\documentclass{amsart}
\usepackage{amssymb, latexsym}
\usepackage{enumerate}
\usepackage{footnote}
\usepackage{xcolor}

\newcommand\supp{\operatorname{supp}}

\newcommand \Isep{[ \bar{t}_{j_{0}}, \bar{\bar{t}}_{j_{0}} ]}

\theoremstyle{plain}
\newtheorem{thm}{Theorem}

\newtheorem{rem}[thm]{Remark}
\newtheorem{claim}[thm]{Claim}
\newtheorem{prop}[thm]{Proposition}
\newtheorem{lem}[thm]{Lemma}
\newtheorem{res}[thm]{Result}

\numberwithin{equation}{section} \numberwithin{thm}{section}

\begin{document}

\title[Global Dynamics For The Energy-Critical Klein-Gordon Equation]
{Global Dynamics Above The Ground State Energy For The Energy-Critical Klein-Gordon Equation}

\author{Tristan Roy}
%\thanks{}

\address{American University of Beirut}
\email{tr14@aub.edu.lb}

%\subjclass{35Q55} \keywords{nonlinear Schr\"odinger equation,
%well-posedness}

\vspace{-0.3in}

\begin{abstract}

Consider the focusing energy-critical Klein-Gordon equation in dimension $d  \in \{ 3,4,5 \}$

\begin{equation}
\left\{
\begin{array}{ll}
\partial_{tt} u  - \Delta u  +  u & = |u|^{\frac{4}{d-2}} u \\
u(0,x) & :=  f_{0}(x) \\
\partial_{t} u(0,x) & := f_{1}(x)
\end{array}
\right. \nonumber
\end{equation}
with data $(f_{0},f_{1}) \in \mathcal{H} := H^{1} \times L^{2}$. We describe the global dynamics of real-valued solutions of which the energy is slightly larger than
that of the ground states'. We classify the flows of the solutions that are ejected from a small neighborhood of the ground states or that are away from them. The classification relies upon a modification of the arguments in \cite{paynesatt} to prove blow-up in finite time, and a modification of the arguments in  \cite{kenmer,kenmerwave,ibramasmnak,kriegnakschlagnonrad} to prove scattering as $t \rightarrow \pm \infty$ . There are three main differences between this paper and \cite{kriegnakschlagnonrad}. The first one is the lack of scaling symmetry. The second one appears in the proof of the ejection lemma: one has to control the mass in the ejection process. The third one appears in the proof of the one-pass lemma: in the worst scenario, one cannot use the equipartition of energy
and therefore one has to prove a decay estimate which allows to use an argument in \cite{bourgjams}. %Another difference appears in the proof of scattering. One has to do extra work because of the lack of scaling when we apply the concentration compactness procedure

\end{abstract}

\maketitle

\section{Introduction}

In this paper we consider the energy-critical Klein-Gordon equation on $\mathbb{R}^{d}$

\begin{equation}
\partial_{tt} u  - \Delta u + u   =   |u|^{2^{*}-2} u \\
\label{Eqn:NlkgCrit}
\end{equation}
with data $(u(0),\partial_{t} u(0)):=(f_{0},f_{1}) $ with $ \mathcal{H} := H^{1} \times L^{2} $. Here $2^{*} := \frac{2d}{d-2}$ denotes the critical Sobolev exponent
(with $d$ the number of dimensions of the space), and $H^{1}$ denotes the standard inhomogeneous Sobolev space, i.e $H^{1}$ is the completion of
the Schwartz space with respect to the norm $ \| f \|_{H^{1}} := \| f \|_{L^{2}} + \| \nabla f \|_{L^{2}} $. \\
\\
In this paper we restrict our attention to (strong) real-valued solutions of (\ref{Eqn:NlkgCrit}), i.e real-valued functions
$ \vec{u} := \left( u , \partial_{t} u \right) \in  \mathcal{C} \left( I, H^{1} \right) \times \mathcal{C} \left( I, L^{2} \right) $ that satisfy the
integral equation below on an interval $I$ containing $0$

\begin{equation}
\begin{array}{ll}
u(t) & = \cos{(t \langle D \rangle)} f_{0} - \frac{ \sin{ \left( t \langle D \rangle \right)}}{\langle D \rangle} f_{1} - \int_{0}^{t}
\frac{ \sin{ \left( (t-t') \langle D \rangle \right)}  }  {\langle D \rangle} \left( |u|^{2^{*} -2} (t') u(t') \right) \, dt' \cdot
\end{array}
\label{Eqn:StrongSol}
\end{equation}
It is well-known that (strong) real-valued solutions $\vec{u}$ of (\ref{Eqn:NlkgCrit}) enjoy the following energy conservation law

\begin{equation}
\begin{array}{ll}
E(\vec{u}) & := \frac{1}{2} \int_{\mathbb{R}^{d}} | \partial_{t}u(t,x)|^{2} \, dx +
\frac{1}{2} \int_{\mathbb{R}^{d}} |\nabla u(t,x)|^{2} \, dx +
 \frac{1}{2} \int_{\mathbb{R}^{d}} |u(t,x)|^{2} \, dx
- \frac{1}{2^{*}} \int_{\mathbb{R}^{d}} |u(t,x)|^{2^{*}} \, dx \cdot
\end{array}
\label{Eqn:Nrj}
\end{equation}
The energy-critical Klein-Gordon equation (\ref{Eqn:NlkgCrit}) is closely related to the energy-critical wave equation, that is

\begin{eqnarray*}
\partial_{tt} u - \triangle u & = & |u|^{2^{*}-2} u \cdot
\label{Eqn:CriticalWave}
\end{eqnarray*}
We recall some properties of (\ref{Eqn:CriticalWave}). It is well-known that of (strong) real-valued solutions
$\vec{u}:= (u, \partial_{t} u) $ of (\ref{Eqn:CriticalWave}) \footnote{i.e real-valued functions $ \vec{u} := \left( u , \partial_{t} u \right) \in  \mathcal{C}
\left( I, \dot{H}^{1} \right) \times \mathcal{C} \left( I, L^{2} \right) $ that satisfy on an interval $I$ containing $0$ (\ref{Eqn:StrongSol}), replacing
``$\langle D \rangle $'' with ``$D$'' } satisfy the following energy conservation law

\begin{eqnarray*}
E_{wa}(\vec{u}) & := \frac{1}{2} \int_{\mathbb{R}^{d}} | \partial_{t} u(t,x)|^{2} \, dx + \frac{1}{2} \int_{\mathbb{R}^{d}}
| \nabla u(t,x)|^{2} \, dx - \frac{1}{2^{*}} \int_{\mathbb{R}^{d}} |u(t,x)|^{2^{*}} \, dx \cdot
\end{eqnarray*}
We can write (\ref{Eqn:CriticalWave})in the Hamiltonian form

\begin{equation}
\begin{array}{l}
\partial_{t} \vec{u} = \mathcal{J} E^{'}_{wa}(\vec{u}), \; \text{with} \; \mathcal{J} :=
\left(
\begin{array}{ll}
0 & 1 \\
-1 & 0
\end{array}
\right)
\end{array}
\nonumber
\end{equation}
Here $E^{'}_{wa}(\vec{v})$ defined in the following fashion \footnote{Recall that the natural dot product on $L^{2}$ is given by
$ \langle f,g \rangle  :=  \int_{\mathbb{R}^{d}}  f g \, dx$ }

\begin{equation}
\begin{array}{ll}
\langle E^{'}_{wa}(\vec{v}), \vec{h} \rangle & := \left(
\partial_{\lambda} E_{wa} \left( \vec{v} + \lambda (h_{1},0) \right)_{|_{\lambda =0}} ,
\partial_{\lambda} E_{wa} \left( \vec{v} + \lambda (0,h_{2}) \right)_{|_{\lambda =0}}
\right) \\
& \\
& = \langle \left( -\triangle v_{1} + |v_{1}|^{2^{*}-2} v_{1}  ,v_{2} \right) , \vec{h} \rangle \cdot
\end{array}
\nonumber
\end{equation}
The symplectic form  $\omega$ associated to this hamiltonian system is defined by

\begin{equation}
\begin{array}{ll}
\omega (u, v) & := \langle \mathcal{J} u, v \rangle, \,  \,
\end{array}
\nonumber
\end{equation}
We now recall the definition of some object that appear in the statement of our theorem. The equation (\ref{Eqn:CriticalWave}) admits a family of stationary solutions $\pm \mathcal{S}$  (or, equivalently, ground states) described by the translation parameter $c \in \mathbb{R}^{d}$ and the scaling parameter $\sigma \in \mathbb{R}$:

\begin{equation}
\begin{array}{ll}
\pm \mathcal{S} & := \pm \left\{  T^{c} \vec{S}^{\sigma}  \vec{W}, \, \vec{W} :=
\left(
\begin{array}{l}
W \\
0
\end{array}
\right)
\right\}
\end{array}
\nonumber
\end{equation}
Here $T^{c}$ and $\vec{S}^{\sigma}$ denote the following operators

\begin{equation}
\begin{array}{l}
T^{c} ( \vec{f} ) := \vec{f} (\cdot - c), \; \text{and} \; \\
\\
\vec{S}^{\sigma} (\vec{f}) := \left( S_{-1}^{\sigma} f_{1}, S_{0}^{\sigma} f_{2} \right)
:= \left( e^{\sigma \left( \frac{d}{2} -1 \right) } f_{1} (e^{\sigma}  \cdot),
e^{\frac{\sigma d}{2} } f_{2} (e^{\sigma}  \cdot) \right),
\end{array}
\nonumber
\end{equation}
with $\vec{f}:= \mathbb{R}^{d} \rightarrow \mathbb{R}^{2}: x \rightarrow  \left( f_{1}(x),f_{2}(x) \right)$ a vector-valued function. Here $W$ satisfies $ - \triangle W = W^{2^{*}-1} $. It is well-known that $ W(x) = \frac{1}{ \left( 1 + \frac{|x|^{2}}{d(d-2)} \right)^{\frac{d-2}{2}}} $.  \\
\\
The local well-posedness theory of (\ref{Eqn:NlkgCrit}) is well-known. %has been studied in \cite{pecher1}.
A consequence of this theory is that one can construct a
maximal time interval of existence $ I(u):=(-T_{-}(u), T_{+}(u) )$. The global behavior of solutions of (\ref{Eqn:NlkgCrit}) with small energy
is now well-understood: in particular, global existence (i.e., $T_{+}(u) = T_{-}(u) = \infty$) and scattering (i.e linear asymptotic behavior) of solutions of
(\ref{Eqn:NlkgCrit}) was proved in \cite{pecher1, pecher2}. The next step is to classify the global behavior of
solutions of (\ref{Eqn:NlkgCrit}) with large energy. %It was proved in \cite{paynesatt} that for energies smaller than that of
%the ground states (i.e $E(\vec{u}) < E_{wa}(\vec{W})$), the solutions blows-up in finite time in some cases (i.e $T_{+}(u) < \infty$
%or $T_{-}(u)< \infty$), and, in other cases, the solutions exist globally in time.
It was proved in \cite{ibramasmnak} that for energies smaller than the ground states' (i.e $E(\vec{u}) < E_{wa}(\vec{W})$), the solutions blows-up in finite time in some cases and, in other cases, the solutions exist globally in time and scatter. The classification of the global behaviors of solutions of (\ref{Eqn:NlkgCrit}) for energies that are equal to the ground states' (i.e $E(\vec{u}) =E_{wa}(\vec{W})$) was studied in \cite{ibramasmnakshold}. The goal of this paper is to classify the global behavior of solutions of (\ref{Eqn:NlkgCrit}) for energies that are slightly above the ground states'. To this end we define for $ \epsilon > 0 $
\footnote{For a definition of $E(\vec{\phi})$, see Subsection \ref{Subsec:NotatFunc} }

\begin{equation}
\begin{array}{l}
\mathcal{H}^{\epsilon}  := \left\{ \vec{\phi} \in \mathcal{H}, \; E(\vec{\phi}) < E_{wa}(\vec{W}) + \epsilon^{2} \right\} \cdot
\end{array}
\nonumber
\end{equation}
We also define for $\mu > 0$ the set

\begin{equation}
\begin{array}{l}
\mathcal{B}^{\mu} := \left\{ \vec{\phi} \in \dot{\mathcal{H}}, \; d_{\mathcal{S}}(\vec{\phi}) \leq \mu \right\} \cdot
\end{array}
\nonumber
\end{equation}
Here $d_{\mathcal{S}}(\vec{f})$ denotes the distance of $\vec{f}$ to $\mathcal{S}$, i.e
$ d_{\mathcal{S}}(\vec{f}) := \inf \limits_{(\sigma,c) \in \mathbb{R} \times \mathbb{R}^{N}} \left\| \vec{f}- T^{c} \vec{S}^{\sigma}(\vec{W}) \right\|_{\dot{\mathcal{H}}} $, with
$ \dot{\mathcal{H}} := \dot{H} \times L^{2}$. \\
\\
We now state the main result of this paper:

\begin{thm}
Let $d \in \{ 3,4,5 \}$. There exist $ 0 <  \epsilon_{*} \ll 1  $, $ 0 < \mu = O(\epsilon_{*}) $, a non-empty set $ \mathcal{T} \subset
\mathcal{B}^{\mu} \cap \mathcal{H}^{\epsilon_{*}} $, and a continuous functional  $\Theta: \mathcal{H}^{\epsilon_{*}} \setminus \mathcal{T} \rightarrow  \left\{ \pm 1 \right\}$ such that for any solution $u$ of (\ref{Eqn:NlkgCrit}) with real data $\vec{\phi} \in \mathcal{H}^{\epsilon_{*}}$ that has maximal time interval of existence $I(u)$ the following properties hold:

\begin{enumerate}

\item[$$]

\item[$(1)$] $ I_{0}(u) := \left\{ t \in I(u): \; \vec{u}(t) \in \mathcal{T} \right\}  $ is an interval,
$I_{+}(u) := \left\{ t \in I(u): \; \vec{u}(t) \notin \mathcal{T}, \; \Theta (\vec{u}(t)) = 1 \right\} $ consists of at most two infinite intervals, and
$I_{-}(u) := \left\{ t \in I(u): \; \vec{u}(t) \notin \mathcal{T}, \; \Theta (\vec{u}(t)) = -1 \right\} $ consists of at most two finite intervals;

\item[$$]

\item[$(2)$] $u$ scatters to a solution of the linear Klein-Gordon equation as $t \rightarrow \pm \infty$ if and only if $ \pm t \in I_{+}(u) $ for large $t > 0$,  and, moreover, $\| u \|_{L_{t}^{\frac{2(d+1)}{d-2}} L_{x}^{\frac{2(d+2)}{d-2}} (I_{\pm}(u))} < \infty$ ;

\item [$$]

\item[$(3)$] For each $\sigma_{1}, \sigma_{2} \in \{ \pm \}$,  let $A_{\sigma_{1},\sigma_{2}}$ be the collection of initial
data $(u_{0},u_{1}) \in \mathcal{H}^{\epsilon_{*}}  $ such that for some $t_{1} < 0  < t_{2}$

\begin{equation}
(- \infty , t_{1}) \cap I(u) \subset I_{\sigma_{1}} (u), \qquad (t_{2}, \infty) \cap I(u) \subset I_{\sigma_{2}}(u).
\nonumber
\end{equation}
Then $A_{\sigma_{1},\sigma_{2}} \neq \emptyset$.

\end{enumerate}

\label{Thm:Main}
\end{thm}

\begin{rem}
Observe that if $\vec{\phi}:= (\phi_{1},\phi_{2}) \in \mathcal{T} $ then not only $ d_{\mathcal{S}}(\vec{\phi}) \lesssim \epsilon_{*} $
(in other words $\vec{\phi}$ is close to $\mathcal{S}$ with $O(\epsilon_{*})$ distance in $\dot{\mathcal{H}}$)
but also the mass is small: more precisely $\| \phi_{1} \|_{L^{2}} \lesssim \epsilon_{*} $.
\end{rem}

Now we explain the main ideas of this paper and how it is organized. \\ \\

In order to prove Theorem \ref{Thm:Main}, we use a strategy based upon two lemmas: the
\textit{ejection lemma} and the \textit{one-pass lemma}. This strategy was designed by the Nakanishi and Schlag in \cite{nakschlagkg} in the study of the global dynamics above the ground states of solutions of the focusing Klein-Gordon equation with a focusing and energy-subcritical power-type nonlinearity (namely
the cubic nonlinearity). It was successfully applied to a focusing Schr\"odinger equation with a energy-subcritical cubic nonlinearity \cite{nakaschlagschrod},
to the focusing energy-critical wave equation with radial data \cite{kriegnakschlagrad} and nonradial data \cite{kriegnakschlagnonrad}, and to the focusing
energy-critical Schr\"odinger equation with radial data \cite{nakashtriroyschrod}. \\ \\
The \textit{ejection lemma} (see Proposition \ref{Prop:DynEjection}) aims at describing the dynamics of the solution if it is close to $\mathcal{S}$ in $\dot{\mathcal{H}}$ and if it moves away from it. In the ejection mode, we would like the dynamics be ruled by the unstable eigenmode of the spectral decomposition of the remainder resulting from the linearization  around $\mathcal{S}$. But this can only be done if we can con control the dynamics of the orthogonal component $\vec{\gamma}$. The dynamics of $\vec{\gamma}$ is estimated by that of the quadratic expansion of the energy  along the orthogonal direction of the remainder, provided that some orthogonality conditions are satisfied. In \cite{kriegnakschlagnonrad}, the authors perform the decomposition (\ref{Eqn:Decompu}) in order to deal with the solutions of the energy-critical wave equation close to $\mathcal{S}$. The decomposition involves two parameters (the translation parameter $c$ and the scaling parameter $\sigma$) that evolve as time goes by so that two orthogonality conditions hold, namely (\ref{Eqn:Orth1}) and (\ref{Eqn:Orth2}). In our work, one has to control extra terms (such as $Res(t)$) that appear due to the ``$u$'' term in (\ref{Eqn:NlkgCrit}). We perform this task by relating the dynamics of these terms to that of the mass of the solution, then by controlling the growth of the mass through again the variation of
the  quadratic expansion of the energy along the orthogonal component. At the end of the ejection, we prove that the dynamics of the solution is dominated by the exponential growth of the unstable eigenmode; moreover, a relevant functional (denoted by $K$) grows exponentially and its sign eventually becomes opposite to that of the eigenmode. \\ \\
The \textit{one-pass lemma} (see Proposition \ref{Prop:OnePassLemma}) aims at classifying the flow of the
solution as the nonlinear distance is in the ejection mode. It shows that the nonlinear distance cannot be at two different times
(say $t_{2}$ and $t_{3}$) too small and in the ejection mode. The proof is by contradiction. Assuming that it is false, then we can apply
the ejection lemma whenever the nonlinear distance is small and then variational estimates (see Proposition \ref{Prop:VarEst}) whenever the nonlinear
distance is large. The contradiction appears when we integrate by part the virial identity
(\ref{Eqn:VirialIdentNLKG}). The left-hand side is much smaller than the right-hand side thanks to the exponential growth
of the functional $K$ in the ejection mode and variational estimates far from
$\mathcal{S}$. The process involves a parameter $m$. We mention the main differences between our work and the previous works. First notice that unlike the
energy-subcritical focusing cubic Klein-Gordon equation there is an additional parameter, i.e the scaling parameter. Notice also that, unlike the energy-critical wave equation (see \cite{kriegnakschlagrad,kriegnakschlagnonrad}), the energy-critical Klein-Gordon equation does not have any scaling property. Because of this lack of symmetry one cannot restrict our analysis to a narrower range of values of this parameter. Consequently one has to perform an analysis taking into account all the possible values of this
parameter. In the case where $\Theta(u) = 1 $, the worst scenario is the following: one has degeneracy of
$K$, i.e the kinetic part is small on average in the region where the nonlinear distance is large. Notice that for the energy-critical wave equation
(see \cite{kriegnakschlagnonrad}) a localized version of the equipartition of energy allowed to prove that this scenario never occurs, more
precisely we always have \footnote{Here $t_{2}$ and $t_{3}$ are defined in Proposition \ref{Prop:OnePassLemma} and its proof}

\begin{equation}
\begin{array}{ll}
\int_{t_{2}}^{t_{3}} \int_{\mathbb{R}^{N}} |\nabla u(t)|^{2} dx \, dt & \gtrsim t_{3}-t_{2}.
\end{array}
\nonumber
\end{equation}
Notice also that if we were to apply this strategy to the energy-critical Klein-Gordon equation, the same argument would yield

\begin{equation}
\begin{array}{ll}
\int_{t_{2}}^{t_{3}} \int_{\mathbb{R}^{d}} |\nabla u(t)|^{2} + |u(t)|^{2} dx \, dt & \gtrsim t_{3}-t_{3},
\end{array}
\nonumber
\end{equation}
which would not be sufficient to prevent this scenario from occurring. Instead, we proceed as follows: we prove a decay estimate (see (\ref{Eqn:DecayEstimate})) that allows to use an argument in \cite{bourgjams}. More precisely, using to our advantage this decay estimate, we prove a propagation estimate involving a frequency localized
part of the solution and dispersive estimates to construct a solution $\tilde{w}$ that is initially close to our solution and such that its energy is smaller than the
energy of the ground states; then, applying the theory for energies below that the ground states  (see \cite{ibramasmnak}), one can prove that $\tilde{w}$
is far from $\mathcal{S}$ and, applying a perturbation argument, we see that the nonlinear distance cannot be too small at $t_{2}$ or $t_{3}$, which is a contradiction. \\ \\
The fate of the solution depends on the value of $\Theta (u)$ once it is ejected.
In the case where $\Theta(u) = -1$ after ejection, we prove the existence of blow-up in finite time by
by modifying an argument due to Payne-Sattinger \cite{paynesatt}. If $\Theta(u) =1$, then we prove that there is scattering: see Section
\ref{Subsection:OnePassLemmaScatt}. The proof is proved by a modification of an approach from \cite{kenmer}. Assuming that scattering fails, then one can find a critical level of energy above which scattering does not hold for solutions that satisfy $\Theta (u) =1$ and such that the nonlinear distance is large. But this means that there exists a sequence of solutions $(u_{n})_{n \geq 1}$ of the energy-critical Klein-Gordon equation (\ref{Eqn:NlkgCrit}) that satisfy the properties that we have just mentioned (in fact, the nonlinear distance can be upgraded from large to very large, by applying to the ejection lemma). Next
we use a linear profile decomposition at an appropriate time (see Proposition \ref{Proposition:ProfileDecomp}) that is proved in
\cite{ibramasmnak}. The process involves a list of sequences of scaling parameters $ \{ h_{n}^{j} \}_{n \in \mathbb{N}, j \in \mathbb{N}}$ that are
either equal to one or tend to zero as $n$ goes to infinity. We then construct for each $j$ a nonlinear profile depending on the
sequence $\{ h_{n}^{j} \}_{n \in \mathbb{N}}$. If $h_{n}^{j}=1$ then the nonlinear profile is obviously a solution of (\ref{Eqn:NlkgCrit}); if
$h_{n}^{j}$ approaches zero as $n$ grows, then we prove that the nonlinear profile is a good approximation of a solution of the energy-critical wave equation. Again, in the latter
case, one has to deal with the lack of scaling property for the energy-critical Klein-Gordon equation. We manage to overcome this issue by proving some estimates that are uniform
throughout the process. finally, the theory below the ground states for the energy-critical wave equation \cite{kenmerwave}, that
for the energy-critical Klein-Gordon equation \cite{ibramasmnak}, and  that just above the ground states for the energy-critical wave equation
\cite{kriegnakschlagnonrad}, one can construct after some work a critical element $U_{c}$ that is a solution of (\ref{Eqn:NlkgCrit}), that does not scatter,  has an energy equal to the critical level of energy, satisfies  $\Theta (U_{c}) > 0 $, and such that its nonlinear distance is large. Moreover its flow is precompact modulo some translation parameter. Finally, the arguments in \cite{ibramasmnak} (see also \cite{kenmerwave}) allow to prove that this critical element does not exist. \\
\\
\textbf{Acknowledgements}:  This project is based upon discussions with Kenji Nakanishi. \\
\\
\textbf{Funding}: This article was partially funded by an URB grant (ID: 3245) from the American University of Beirut. This research was also partially conducted in Japan and funded by a (JSPS) Kakenhi grant [15K17570 to T.R.].

%Part of this work was done in Japan where the author was supported by an JSPS fellowship.

\section{Notation}

In this section, we set up some notation that appear in this paper.

\subsection{General Notation}

If $y$ is a real number, then $y+$, $y-$ is a slightly larger, smaller number respectively. If $x \in \mathbb{R}$ then $sign(x):=1$ if $x \geq 0$ and
$sign(x):= -1$ if $x < 0$. Let $\vec{f} := \left( \begin{array}{l} f_{1} \\ f_{2} \end{array} \right) $ be a vector-valued function that is made of two real-valued functions: $f_{1}$ and $f_{2}$.
If $\vec{u}$ is a real-valued solution of (\ref{Eqn:NlkgCrit}) or a solution of (\ref{Eqn:CriticalWave}) then we write $u_{1} := u $, $u_{2} := \partial_{t} u$, and $\vec{u} :=(u,\partial_{t}u)$. Unless otherwise specified, we do not mention specify in the sequel the spaces to which the functions belong.

%Let $g$ be a function defined on an interval $J$. Given  $t \in J$, we denote by $\partial_{t}^{+} g(t) $, $\partial_{t}^{-} g(t)$ the following
%numbe

%\begin{equation}
%\begin{array}{ll}
%\partial_{t}^{+} g(t) & = \lim_{h \rightarrow 0, h > 0} \frac{g(t+h) - g(t)}{h}, \; \text{and}   \\
%\partial_{t}^{-} g(t) & = \lim_{h \rightarrow 0, h < 0} \frac{g(t+h) - g(t)}{h} \cdot
%\end{array}
%\nonumber
%\end{equation}

\subsection{Paley-Littlewood decomposition}

Some estimates that we establish throughout this paper require the Paley-Littlewood technology. We set it up
now. Let $f$ be a function defined on $\mathbb{R}^{d}$. Let $\phi(\xi)$ be a real, radial, nonincreasing function that is equal to $1$ on the unit ball
$\left\{ \xi \in \mathbb{R}^{N}: |\xi| \leq 1 \right\}$ and that is supported on $\left\{ \xi \in \mathbb{R}^{N}: |\xi| \leq 2  \right\}$.
Let $\psi$ denote the function $\psi(\xi):=  \phi(\xi) - \phi(2 \xi)$. If $ \alpha \in \mathbb{R}^{+}_{*} $ then we define the Paley-Littlewood operators
$ P_{\alpha} f $, $P_{<\alpha} f$, and $ P_{\geq \alpha} f $ in the Fourier domain by

\begin{equation}
\begin{array}{ll}
\widehat{P_{\alpha} f}(\xi) & := \psi \left( \frac{\xi}{\alpha} \right) \widehat{f}(\xi), \\
\widehat{P_{<\alpha} f}(\xi) & := \phi \left( \frac{\xi}{\alpha} \right) \widehat{f}(\xi), \; \text{and} \\
\widehat{P_{\geq \alpha} f}(\xi) & := \widehat{f}(\xi) - \widehat{P_{< \alpha} f}(\xi) \cdot
\end{array}
\nonumber
\end{equation}
We also define $P_{0} f $ by the formula  $ \widehat{P_{0} f}(\xi) := \phi(\xi) \hat{f}(\xi) $. Observe that we can write
$f = P_{0} f  + \sum \limits_{M \in 2^{\mathbb{N}}} P_{M} f  $: this decomposition is referred to as the Paley-Littlewood decomposition of $f$.\\
\\
It is well-known that there exists $\bar{C} > 0$ such that $ \| P_{0} f \|_{L^{2^{*}}} \leq \bar{C} \| f  \|_{L^{2^{*}}}$ and
$\| P_{M} f \|_{L^{2^{*}}} \leq \bar{C} \| f \|_{L^{2^{*}}}$. Hence, by defining the ``normalized'' Paley-Littlewood projectors
$\bar{P}_{0} f := \bar{C}^{-1} f $ and $ \bar{P}_{M} f := \bar{C}^{-1} P_{M} f $, we get $ \| \bar{P}_{0} f \|_{L^{2^{*}}} \leq \| f  \|_{L^{2^{*}}}$ and
$\| \bar{P}_{M} f \|_{L^{2^{*}}} \leq  \| f \|_{L^{2^{*}}}$ \footnote{ We will use these inequalities in Subsection
\ref{Subsection:CondImplBourg} : see the estimate  $I \left( \chi_{R'}(0) \bar{P}_{< M^{'}_{j_{0}}} u(0) \right) \leq I( u(0)) $ }. We also have  $ f = \bar{C} \left( \bar{P}_{0} f  + \sum \limits_{M \in 2^{\mathbb{N}}} \bar{P}_{M} f \right)  $: this decomposition is referred to as the ``normalized'' Paley-Littlewood decomposition of $f$.

\subsection{Operators}

We define some operators that we use throughout this paper. \\
\\
Let $f$ be a function defined on $\mathbb{R}^{d}$. Let $(\sigma, a) \in \mathbb{R}^{2} $  and $S^{\sigma}_{a}$ be the operator defined as follows

\begin{equation}
\begin{array}{ll}
S^{\sigma}_{a} f(x) & := e^{ \left( \frac{d}{2} + a \right) \sigma} f(e^{\sigma} x)
\end{array}
\nonumber
\end{equation}
Let $S^{'}_{a}:= \partial_{\sigma} S^{\sigma}_{a} (\sigma=0)$. We define the operator
$\vec{\Lambda} := \left( S_{1}^{'}, S_{0}^{'} \right) $. Let $(S^{\sigma}_{a})^{*}$  be the adjoint of $S^{\sigma}_{a}$. It is easy to see that
$(S_{a}^{\sigma})^{*}= S_{-a}^{- \sigma}$. By differentiating
$ \langle S^{\sigma}_{a} \phi_{1}, \phi_{2} \rangle = \langle \phi_{1}, (S^{\sigma}_{a})^{*} \phi_{2} \rangle $  with respect to $\sigma$
at $\sigma=0$ we see that

\begin{equation}
\begin{array}{ll}
(S^{'}_{a})^{*} & = - S^{'}_{-a}
\end{array}
\label{Eqn:AdjointS}
\end{equation}
Let $ \vec{\Lambda}^{*} := \left( S^{'}_{-1}, S^{'}_{0} \right) $. \\
\\
We will work with families of operators. Given $(h_{n}^{j}, x^{j}_{n}) \in (0,1] \times \mathbb{R}^{d}$ let $\tau^{j}_{n}$, $T_{n}^{j}$ and
$ \langle D \rangle_{n}^{j} $ denote the scaled time shift, the unitary and self-adjoint operators, that is

\begin{align*}
\tau_{n}^{j}:= - \frac{ t_{n}^{j}}{h_{n}^{j}}, \, T_{n}^{j} f (x) := ( h_{n}^{j} )^{-\frac{d}{2}} f \left( \frac{ x -x_{n}^{j}}{h_{n}^{j}} \right), \; \text{and} \;
\langle D \rangle_{n}^{j} := \sqrt{-\triangle + (h_{n}^{j})^{2}} \cdot
\end{align*}

\subsection{Norms and Function spaces}

With $h_{n}^{j}$ defined above let $\mathcal{H}_{h_{n,j}} := H_{h_{n,j}} \times L^{2} $ where $H_{h_{n}^{j}}$ is  the closure of the Schwartz space with respect to the norm

\begin{equation}
\begin{array}{ll}
\| f \|_{H_{h_{n}^{j}}} := \| \langle D \rangle_{n}^{j} f \|_{L^{2}} \cdot
\end{array}
\nonumber
\end{equation}
Let $\sigma \in \mathbb{R} - \{ 0 \}$  and  $\check{q} \geq 1$. Let $M_{n}^{j}$ be the smallest dyadic number such that
$M_{n}^{j} \geq h_{n}^{j}$. With $h_{n}^{j}$ defined above we define the following quantity

\begin{equation}
\begin{array}{l}
\| f \|_{\bar{B}^{\sigma,h_{n}^{j}}_{\check{q},2}} := (h_{n}^{j})^{\sigma}
\left\| P_{\lesssim h_{n}^{j}} f \right\|_{L^{\check{q}}} +  \left( \sum \limits_{M \in 2^{\mathbb{Z}} : M > M_{n}^{j}}
 M^{ 2 \sigma} \| P_{M} f \|^{2}_{L^{\check{q}} } \right)^{\frac{1}{2}} \cdot
\end{array}
\nonumber
\end{equation}
Observe that if $h_{n,j} = 1 $ then $ \| f \|_{\bar{B}^{\sigma,h_{n,j}}_{\check{q},2}} = \| f \|_{B^{\sigma}_{\check{q},2}}
$, with $B^{\sigma}_{\check{q},2}$ belonging to the class of well-known inhomogeneous Besov spaces.

\subsection{Sobolev-type embedding}

We recall some Sobolev-type embeddings that we use throughout this paper.  \\
\\
Let $f$ be a function depending on $\mathbb{R}^{d}$. Then we have

\begin{equation}
\begin{array}{l}
\| f \|_{L^{2^{*}}} \lesssim C \| f \|_{\dot{H^{1}}} \cdot
\end{array}
\label{Eqn:SobIneq}
\end{equation}
It is known (see \cite{aubin,tal}) that $W$ is an extremizer for the inequality (\ref{Eqn:SobIneq}), i.e

\begin{equation}
\begin{array}{l}
C_{*} := \sup \left\{ \frac{\| f \|_{L^{2^{*}}}}{\| f \|_{\dot{H}^{1}}}, \, f  \in \dot{H}^{1}, \, f \neq 0 \ \right\} = \frac{ \| W \|_{L^{2^{*}}} } { \| W \|_{\dot{H}^{1}}} \cdot
\end{array}
\label{Eqn:SobIneqMax}
\end{equation}

%??????????????

%Let $v$ be a function depending on time and space. Then we have

%\begin{equation}
%\begin{array}{l}
%\| v \|_{(\dot{0} ,J)} = \| v \|_{(0,J)}  \lesssim
%\| v \|_{\left( \dot{A} , J \right)} \lesssim  \| v \|_{(A,J)} \\
%\\
%\| v \|_{(\dot{0} , J)} = \| v \|_{( \dot{0'}, J )} = \| u \|_{(\dot{0}, h_{n}^{j},J)}  \lesssim \| v \|_{ \left( \oneoftwo \cdot, J \right)}
%\lesssim \| v \|_{ \left( \oneoftwo, h_{n}^{j},J \right) } \lesssim \| v \|_{\left( \oneoftwo,J \right)}
%\end{array}
%\label{Eqn:SobIneqStrich}
%\end{equation}

%?????????????

\subsection{Dispersive estimates and nonlinear estimates}

We recall dispersive estimates and nonlinear estimates that we will use in this paper. \\
\\
Let $f$ and $g$ be two functions depending on $\mathbb{R}^{d}$. Let $F(f) := |f|^{2^{*}-2} f$. Let
$
(R, \bar{R})  \in
\left(  ( B^{\frac{1}{2}}_{\frac{2(d+1)}{d+3},2} ,  B^{\frac{1}{2}}_{\frac{2(d+1)}{d-1},2} ),
( \dot{B}^{\frac{1}{2}}_{\frac{2(d+1)}{d+3},2} ,  \dot{B}^{\frac{1}{2}}_{\frac{2(d+1)}{d-1},2} ), \,
( \bar{B}^{\frac{1}{2},h_{n}^{j}}_{\frac{2(d+1)}{d+3},2} ,  \bar{B}^{\frac{1}{2},h_{n}^{j}}_{\frac{2(d+1)}{d-1},2} )
\right) $. Then the nonlinear estimates read as follows

\begin{equation}
\begin{array}{ll}
\left\| F(f) \right\|_{R} & \lesssim
\| f \|_{\bar{R}}  \| f \|^{2^{*}-2}_{L^{\frac{2(d+1)}{d-2}}}, \; \text{and}
\end{array}
\label{Eqn:NonlinearEstGen}
\end{equation}

\begin{equation}
\begin{array}{ll}
\left\| F(f) - F(g) \right\|_{R} &
\lesssim \| f - g \|_{\bar{R}}
\left( \| f \|^{2^{*}-2}_{L^{\frac{2(d+1)}{d-2}}} + \| g \|^{2^{*}-2}_{L^{\frac{2(d+1)}{d-2}}} \right) \\
& + \| f- g \|_{L^{\frac{2(d+1)}{d-2}}} \left(  \| f \|^{2^{*}-3}_{L^{\frac{2(d+1)}{d-2}}} + \| g \|^{2^{*}-3}_{L^{\frac{2(d+1)}{d-2}}} \right) \\
& \left( \| f \|_{\bar{R}} + \| g \|_{\bar{R}} \right) \cdot
\end{array}
\label{Eqn:NonlinearEstGenDiff}
\end{equation}

We will combine the estimates above with the Strichartz-type estimates below with  $(q,r)$ two numbers that satisfy $q > 2$ and $\frac{1}{q} + \frac{d}{r} = \frac{d}{2} - 1 $. \\
\\
Let $w$ be a solution of $ \partial_{tt} w - \triangle w + w = F + G $ on $J=[a, \cdot]$. Recall the Duhamel formula, that is $w(t)= w_{l,a}(t) + w_{nl,a}(t)$ with
$ w_{l,a}(t) := \cos{\left( (t-a)\langle D \rangle \right)} w(a) + \frac{\sin{ \left((t-a) \langle D \rangle \right) }}{ \langle D \rangle} \partial_{t} w(a) $  and $ w_{nl,a}(t) := - \int_{a}^{t} \frac{\sin \left( (t- t^{'}) \langle D \rangle \right)}{\langle D \rangle} \left( F(t') + G(t') \right) \; dt' $. $w_{l,a}$ (resp. $w_{nl,a}$) denotes the linear part (resp. the nonlinear part of $w$ starting from $a$; moreover the following wave-type Strichartz estimate holds

\begin{equation}
\begin{array}{l}
\left\|  \left( w , \partial_{t} w \right) \right\|_{L_{t}^{\infty} \mathcal{H} (J)} + \| w \|_{L_{t}^{q} L_{x}^{r} (J)} +
\| \vec{w} \|_{L_{t}^{\frac{2(d+1)}{d-1}} B^{\frac{1}{2}}_{\frac{2(d+1)}{d-1},2} (J) \times L_{t}^{\frac{2(d+1)}{d-1}} B^{-\frac{1}{2}}_{\frac{2(d+1)}{d-1},2} (J) }  \\
\lesssim \left\| \left( w(a) , \partial_{t} w(a) \right) \right\|_{\mathcal{H}} + \| F \|_{L_{t}^{1} L_{x}^{2}(J)} + \| G \|_{L_{t}^{\frac{2(d+1)}{d-3}} B^{\frac{1}{2}}_{\frac{2(d+1)}{d-3},2} (J)}
\end{array}
\label{Eqn:StrichKg0}
\end{equation}

Let $w$ be a solution of  $ \partial_{tt} w - \triangle w = F + G $ on $J=[a,\cdot]$. Recall the Duhamel formula that is $w(t)= w_{l,a}(t) + w_{nl,a}(t) $ with
$ w_{l,a}(t) := \cos{((t-a)D)} w(a) + \frac{\sin{ \left((t-a) D \right) }}{D} \partial_{t} w(a) $  and $ w_{nl,a}(t) := - \int_{a}^{t} \frac{\sin \left( (t- t^{'}) D \right)}{D} \left( F(t') + G(t') \right) \; dt' $. $w_{l,a}$ (resp. $w_{nl,a}$) denotes the linear part (resp. the nonlinear part of $w$ starting from $a$; moreover the following Strichartz-type estimate holds

\begin{equation}
\begin{array}{l}
\left\|  \left( w (t), \partial_{t} w \right) \right\|_{L_{t}^{\infty} \dot{\mathcal{H}} (J)} + \| w \|_{L_{t}^{q} L_{x}^{r}(J)} +
\left\| (w, \partial_{t} w ) \right\|_{L_{t}^{\frac{2(d+1)}{d-1}} \dot{B}^{\frac{1}{2}}_{\frac{2(d+1)}{d-1},2} (J) \times L_{t}^{\frac{2(d+1)}{d-1}} \dot{B}^{-\frac{1}{2}}_{\frac{2(d+1)}{d-1},2} (J) }
 \\
\lesssim \| \left( w(a) , \partial_{t} w(a) \right) \|_{\dot{\mathcal{H}}} + \| F \|_{L_{t}^{1} L_{x}^{2}(J)} + \| G \|_{L_{t}^{\frac{2(d+1)}{d+3}} \dot{B}^{\frac{1}{2}}_{\frac{2(d+1)}{d+3},2} (J)} \cdot
\end{array}
\label{Eqn:StrichWave0}
\end{equation}

If $w$ satisfies $ \partial_{tt} w - \triangle w + (h_{n,j})^{2} w  = F + G  $ on $ J=[a,\cdot] $ then we can write similarly $w(t)= w_{l,a}(t) + w_{nl,a}(t)$ with
$ w_{l,a}(t) := \cos{((t-a) \langle D \rangle_{n}^{j})} w(a) + \frac{\sin{ \left( (t-a) \langle D \rangle_{n}^{j} \right)}}{\langle D \rangle_{n}^{j}}
\partial_{t} w(a)  $ and  $ w_{nl,a}(t) := - \int_{a}^{t} \frac{\sin \left( (t- t^{'}) \langle D \rangle_{n}^{j}  \right)}{\langle D \rangle_{n}^{j} } \left( F(t') + G(t') \right) \; dt' $. Moreover we have

\begin{equation}
\begin{array}{l}
\left\|  \left( w (t), \partial_{t} w \right) \right\|_{L_{t}^{\infty} \mathcal{H}_{h_{n}^{j}} (J)} + \| w \|_{L_{t}^{q} L_{x}^{r}(J)}
+ \left\| (w,\partial_{t} w) \right\|_{L_{t}^{\frac{2(d+1)}{d-1}} \bar{B}^{\frac{1}{2},h_{n,j}}_{\frac{2(d+1)}{d-1},2} (J)
\times L_{t}^{\frac{2(d+1)}{d-1}} \bar{B}^{-\frac{1}{2},h_{n,j}}_{\frac{2(d+1)}{d-1},2} (J)}  \\
\lesssim \left\|  \left( w(a), \partial_{t} w(a) \right)  \right\|_{\mathcal{H}_{h_{n}^{j}}} + \| F \|_{L_{t}^{1} L_{x}^{2}(J)} + \| G \|_{L_{t}^{\frac{2(d+1)}{d-3}} \bar{B}^{\frac{1}{2},h_{n}^{j}}_{\frac{2(d+1)}{d-3},2} (J)} \cdot
\end{array}
\label{Eqn:StrichWave0Sc}
\end{equation}

\begin{rem}
The estimates (\ref{Eqn:NonlinearEstGen}) and (\ref{Eqn:NonlinearEstGenDiff}) for $(R, \bar{R}) \neq ( \bar{B}^{\frac{1}{2},h_{n}^{j}}_{\frac{2(d+1)}{d+3},2} ,  \bar{B}^{\frac{1}{2},h_{n}^{j}}_{\frac{2(d+1)}{d-1},2} )$ can be proved by using finite differences characterizations of homogeneous Besov spaces:
see e.g. \cite{bulutlipav}. In order to prove (\ref{Eqn:NonlinearEstGen}) and (\ref{Eqn:NonlinearEstGenDiff}), we apply
(\ref{Eqn:NonlinearEstGen}) and (\ref{Eqn:NonlinearEstGenDiff}) for $(R, \bar{R}) = \left( B^{\frac{1}{2}}_{\frac{2(d+1)}{d+3},2} ,
\bar{B}^{\frac{1}{2}}_{\frac{2(d+1)}{d-1},2} \right)$ to $ h_{n}^{j} \left( T_{n}^{j} f, T_{n}^{j} g \right)$ instead of $(f,g)$, taking into account
$ (*): \, P_{M } \left( f \left( \frac{\cdot}{\lambda} \right) \right) = P_{\lambda M} f \left( \frac{\cdot}{\lambda} \right) \,
\left(  \text{resp.} \,  P_{ \leq M  } \left( f \left( \frac{\cdot}{\lambda} \right) \right) = P_{\leq \lambda M} f \left( \frac{\cdot}{\lambda} \right) \right) $.
The estimates (\ref{Eqn:StrichKg0}) and (\ref{Eqn:StrichWave0}) are well-known. For papers dealing with Strichartz-type estimates for the wave equation and
the Klein-Gordon equation see e.g \cite{keeltao,ginebvelo,linsog,ibramasmnak} and references therein. In order to prove (\ref{Eqn:StrichWave0Sc}) observe that
if $w$ satisfies $ \partial_{tt} w - \triangle w + (h_{n}^{j})^{2} w  = H  $, then
$ \bar{w}(t,x) := w \left( \frac{t}{h_{n}^{j}}, \frac{x}{h_{n}^{j}} \right) $ is a solution of
$\partial_{tt} \bar{w} - \triangle \bar{w} + \bar{w} = (h_{n}^{j})^{-2} H \left( \frac{t}{h_{n}^{j}}, \frac{x}{h_{n}^{j}} \right)$;
then use (\ref{Eqn:StrichKg0}) and $ (*) $ to get (\ref{Eqn:StrichWave0Sc}).
\end{rem}
We will mostly use (\ref{Eqn:StrichKg0}), (\ref{Eqn:StrichWave0}), and (\ref{Eqn:StrichWave0Sc}) with $ (q,r) = \frac{2(d+1)}{d-2} (1,1) $ or $(q,r) := \frac{d+2}{d-2} (1,2) $.

\subsection{Functionals}
\label{Subsec:NotatFunc}

In this paper we will often use the functionals defined below.\\

Let $\vec{g} := (g_{1},g_{2})$ with $g_{1}$ and $g_{2}$ two real-valued functions defined on $\mathbb{R}^{d}$. We denote by $E(\vec{g})$, $E_{wa}(\vec{g})$ the following expressions (for $\vec{g}$ and $(g_{1},0)$ lying in the appropriate spaces)

\begin{equation}
\begin{array}{l}
E(\vec{g}) :=  \frac{1}{2} \int_{\mathbb{R}^{d}} g_{2}^{2} \; dx + \frac{1}{2} \int_{\mathbb{R}^{d}} |\nabla g_{1}|^{2} \; dx
+ \frac{1}{2} \int_{\mathbb{R}^{2}} g_{1}^{2} \; dx - \frac{1}{2^{*}} \int_{\mathbb{R}^{d}} |g_{1}|^{2^{*}} \; dx \\
\\
E_{wa}(\vec{g}) := \frac{1}{2} \int_{\mathbb{R}^{d}} g^{2}_{2} \; dx   + \frac{1}{2} \int_{\mathbb{R}^{d}} |\nabla g_{1}|^{2} \; dx
 - \frac{1}{2^{*}} \int_{\mathbb{R}^{d}} |g_{1}|^{2^{*}} \; dx \cdot \\
\end{array}
\nonumber
\end{equation}
More generally, given $h$ a function on $\mathbb{R}^{d}$, $R>0$ and $c \in \mathbb{R}^{d}$, let

\begin{align*}
E(h,\vec{g})  := \frac{1}{2} \int_{\mathbb{R}^{d}}  h  | g_{2}|^{2} \, dx + \frac{1}{2} \int_{\mathbb{R}^{d}} h  |\nabla g_{1}|^{2} \, dx
+ \frac{1}{2} \int_{\mathbb{R}^{d}} h |g_{1}|^{2} \, dx - \frac{1}{2^{*}} \int_{\mathbb{R}^{d}} h |g_{1}|^{2^{*}} \, dx, \, and \\
\\
E_{R}(c,\vec{g}) := \frac{1}{2} \int_{|x -c| \geq R }  g_{2}^{2} + \frac{1}{2} \int_{|x -c| \geq R} | \nabla g_{1}|^{2} + \frac{1}{2} \int_{|x-c| \geq R}
|g_{1}|^{2} \, dx - \frac{1}{2^{*}} \int_{|x-c| \geq R} |g_{1}|^{2^{*}} \, dx.
\end{align*}
Let $f$ be a function defined on $\mathbb{R}^{d}$. We denote by $K(f)$, $I(f)$, and $G(f)$ the following numbers

\begin{equation}
\begin{array}{l}
K(f) := \| \nabla f \|^{2}_{L^{2}} - \| f \|^{2^{*}}_{L^{2^{*}}}, \;
I(f) :=  \frac{1}{d} \| f \|^{2^{*}}_{L^{2^{*}}}, \; \text{and} \;
G(f) := \frac{1}{d} \| f \|^{2}_{\dot{H}^{1}} \cdot
\end{array}
\nonumber
\end{equation}
Observe that

\begin{equation}
\begin{array}{l}
E(\vec{g}) - \frac{K(g_{1})}{2} = \frac{1}{2} \| g_{2} \|^{2}_{L^{2}} + \frac{1}{2} \| g_{1} \|^{2}_{L^{2}} + I(g_{1}), \;
E_{wa}(\vec{g}) - \frac{K_(g_{1})}{2} = \frac{1}{2} \| g_{2} \|^{2}_{L^{2}} + I(g_{1}),
\end{array}
\label{Eqn:RelI}
\end{equation}
%\label{Eqn:DefI} \label{Eqn:RelNrjK}

\begin{equation}
\begin{array}{l}
E(\vec{g}) - \frac{K(g_{1})}{2^{*}}  = \frac{1}{2}  \left( \|  g_{2} \|^{2}_{L^{2}} + \| g_{1} \|^{2}_{L^{2}} \right) + G(g_{1}), \; \text{and} \\
E(\vec{g}) - \frac{K(g_{1})}{2^{*}} = \frac{1}{2}  \|  g_{2} \|^{2}_{L^{2}} + G(g_{1}) \cdot
\end{array}
\label{Eqn:RelG}
\end{equation}
It follows from (\ref{Eqn:SobIneqMax}) that

\begin{equation}
\begin{array}{ll}
E_{wa}(\vec{W}) & = \inf \{ E_{wa}((g_{1},0)), \, g_{1} \neq 0, \, g_{1} \in \dot{H}^{1}, K(g_{1})=0 \}
\end{array}
\label{Eqn:Charac2GdState}
\end{equation}
We also have

\begin{equation}
\begin{array}{ll}
E_{wa} (\vec{W}) & = \inf \{ G_(f), f \neq 0, \, f \in \dot{H}^{1}, K(f) \leq 0 \}
\end{array}
\label{Eqn:Charac3GdState}
\end{equation}
and

\begin{equation}
\begin{array}{ll}
E_{wa}(\vec{W}) & = \inf \{ I(f), f \neq 0, \, f \in \dot{H}^{1}, K(f) \leq 0 \} \cdot
\end{array}
\label{Eqn:Charac4GdState}
\end{equation}
%Let

%\begin{align*}
%\vec{P}(\vec{g}) := \int_{\mathbb{R}^{d}} g_{2}  \nabla g_{1} \, dx.
%\end{align*}
%It is well-known that if $\vec{u}$ is a solution of (\ref{Eqn:NlkgCrit}), then $\vec{P}(\vec{u})$ is nothing but the momentum and it is conserved. \\ \\

%\footnote{ Here $ p:= \frac{2(n+2)}{n-2}+$,  $ q:=\frac{2(n+2)}{n-2} -$
%denotes a slightly larger, smaller number than $ \frac{2(n+2)}{n-2} $ such that $\frac{1}{p} + \frac{n}{2q} = \frac{n}{4}$  }

\subsection{Linearized operators}

In this paper, we constantly use the linearized operator $L_{+}$ defined by

\begin{align*}
L_{+}   := - \triangle - (2^{*}-1) W^{2^{*}-2}
\end{align*}
We recall some spectral properties of $ L_{+}$:

\begin{itemize}

\item the discrete spectrum consists of a unique negative eigenvalue (denoted by $-k^{2}$ ) and
there exists a unique smooth, positive, exponentially decaying eigenfunction $\rho$ such that

\begin{align}
L_{+} \rho = - k^{2} \rho, \, \rho >0, \, \| \rho \|_{L^{2}} =1
\label{Eqn:DefDiscSpect}
\end{align}

\item the essential spectrum  is $[0, \infty)$

\item the threshold $0$ of the essential spectrum of $L_{+}$ is an eigenvalue: indeed
$L_{+} \left(  \partial_{x_{i}} W \right) = 0 $ for all $1 \leq i \leq d $. Moreover
$ L_{+} \left( \Lambda_{-1} W \right) =0 $.

\end{itemize}

\section{Proof of Theorem \ref{Thm:Main}}
\label{Section:ThmMain}

The proof of Theorem \ref{Thm:Main} relies upon some propositions that we state below. \\
In the first proposition, proved in \cite{kriegnakschlagnonrad}, we perform a decomposition of $u$ close to the ground states, taking into
account the symmetries of the equation and some constraints (the so called orthogonality conditions)

\begin{prop}%{\textbf{`` Orthogonal Decomposition Of $u$ ''}}
There is a small constant $ 0 <  \delta_{s} \ll 1 $  such that for all  $\vec{f} \in \mathcal{V} := \left\{  \vec{g} \in \dot{\mathcal{H}}: \; d_{\mathcal{S}} (\vec{g} ) < \delta_{s} \right\} $ there exists a $\mathcal{C}^{1}-$ function

\begin{equation}
\begin{array}{ll}
\mathcal{V} & \rightarrow \mathbb{R}^{d} \times \mathbb{R} \times \dot{\mathcal{H}}  \\
\vec{f} & \rightarrow ( c, \sigma, \vec{v}) := \left( c( \vec{f} ), \sigma( \vec{f} ), \vec{v}(f) \right)
\end{array}
\nonumber
\end{equation}
such that

\begin{equation}
\vec{f} := T^{c} \vec{S}^{\sigma} ( \vec{W} + \vec{v} ),
\label{Eqn:Decompu}
\end{equation}
and such that $\vec{v}$ satisfies two orthogonality conditions, namely

\begin{equation}
\langle v_{1}, S'_{0} \rho \rangle  = 0, \qquad and
\label{Eqn:Orth1}
\end{equation}

\begin{equation}
\langle v_{1}, \nabla \rho \rangle = 0.
\label{Eqn:Orth2}
\end{equation}

\label{Prop:OrthDecompu}
\end{prop}

\begin{rem}

Let  $ \mathcal{L}  :=  \left(
\begin{array}{ll}
L_{+} & 0 \\
0 & 1
\end{array}
\right)   $. We can decompose the vector $\vec{v}$ as a linear combination of (generalized) eigenvectors of $\mathcal{J} \mathcal{L}$ plus a remainder term.
More precisely let $\vec{g}^{^{\pm}}:= \frac{ (1, \pm k) \rho} { \sqrt{2k}}$. Then $ \mathcal{J} \mathcal{L} \vec{g}^{^{\pm}} = \pm k \vec{g}^{^{\pm}}$. Let $\vec{\gamma}$ be such that

\begin{align}
\vec{v} = \lambda_{+} \vec{g}^{^{+}} + \lambda_{-} \vec{g}^{^{-}} + \vec{\gamma},
\label{Eqn:Decompv}
\end{align}
with $ \lambda_{+} := \omega ( \vec{v}, \vec{g}^{^{-}} ) $ and  $ \lambda_{-} := - \omega ( \vec{v} , \vec{g}^{^{+}} )$. Let $\lambda_{1}:= \frac{1}{\sqrt{2k}} \left( \lambda_{+} + \lambda_{-}  \right)$  and $ \lambda_{2}:= \frac{\sqrt{k}}{2} \frac{\lambda_{+} - \lambda_{-}}{2}$. Then

\begin{align}
v_{i} = \lambda_{i} \rho  + \gamma_{i}, & \qquad \, i \in [1..2] \cdot
\label{Eqn:Decompvi}
\end{align}
Observe that $ \lambda_{i} : = \langle v_{i}, \rho \rangle $ and that $ \langle \gamma_{i}, \rho \rangle  = 0$.

\label{Rem:Decompv}
\end{rem}

\begin{rem}

We will apply the decomposition (\ref{Eqn:Decompu}) to $\vec{u}$, solution of (\ref{Eqn:NlkgCrit}) on an interval $J$ such that
$ \vec{u}(t) \in \mathcal{V}, \; t \in J $. From Proposition \ref{Prop:OrthDecompu} we see that there exists a $\mathcal{C}^{1}(J)-$ function $ t \rightarrow
\left( c(t), \sigma(t), \vec{v}(t) \right) $ such that for all $t \in J$

\begin{equation}
\vec{u}(t) = T^{c(t)} \vec{S}^{\sigma(t)} \left( \vec{W} + \vec{v}(t) \right),
\label{Eqn:Decomput}
\end{equation}
$\langle v_{1}(t), S'_{0} \rho \rangle = 0 $ , and $ \langle v_{1}(t), \nabla \rho \rangle = 0 $.

\end{rem}

The next proposition, proved in Section \ref{Section:LinearParamW}, aims at describing the dynamics of the solution near the ground states, using the
decomposition in Proposition \ref{Prop:OrthDecompu}:

\begin{prop}%{\textbf{`` Dynamics Around Ground States ''}}
Let $\vec{u}$ be a solution of (\ref{Eqn:NlkgCrit}) on an interval $I:= [0,.] $ such that for all $t \in I$ we have $ \vec{u}(t) \in \mathcal{V}$.
Consider the decomposition (\ref{Eqn:Decomput}). Let $\tau$ be such that $\partial_{t} \tau:= e^{\sigma (t)}$ and $\tau(0)=0$. Then

\begin{enumerate}

\item[$(1)$]

\begin{equation}
\partial_{\tau} \vec{v}   = \mathcal{J} \mathcal{L} \vec{v} +  \left( e^{\sigma} \partial_{\tau} c \cdot \nabla  - \partial_{\tau} \sigma  \vec{\Lambda} \right)
( \vec{W} + \vec{v} ) + \underline{N} (\vec{v} )
- e^{- 2 \sigma}
\left(
\begin{array}{l}
0 \\
W  + v_{1}
\end{array}
\right)
\label{Eqn:Linearization}
\end{equation}
Here $ \underline{N}(\vec{v}) := \left( \begin{array}{l} 0 \\ N(v_{1})   \end{array} \right) $ with $N(f) := |W + f|^{2^{*} -2} (W + f) - W^{2^{*} -1} - (2^{*} -1) W^{2^{*} -1} f$.

\item[$$]

\item[$(2)$] We now apply the decomposition (\ref{Eqn:Decompv}) to $v(t)$. We have

\begin{equation}
\partial_{\tau} \lambda_{\pm} = \pm k \lambda_{\pm} \mp \omega \left( \vec{v},  \nabla \vec{g}^{^{\mp}} \right) \cdot \partial_{\tau} c e^{\sigma}
\mp \omega \left( \vec{v}, \vec{\Lambda}^{*} \vec{g}^{^{\mp}} \right) \partial_{\tau} \sigma  \mp \omega \left( \underline{N}(\vec{v}), \vec{g}^{^{\mp}} \right) +
e^{- 2 \sigma }
\omega \left(
\left(
\begin{array}{l}
0 \\
W + v_{1}
\end{array}
\right), \vec{g}^{^{\mp}}
\right),
\label{Eqn:EstDynLambdaPlus}
\end{equation}

\begin{equation}
\begin{array}{l}
\partial_{\tau} \lambda_{1} =  \lambda_{2}   + \frac{1}{\sqrt{2k}} e^{- 2 \sigma }
\omega \left( \left(
\begin{array}{l}
0 \\
W + v_{1}
\end{array}
\right), \vec{g}^{^{-}} + \vec{g}^{^{+}}
\right), \\
\partial_{\tau} \lambda_{2}  =  k^{2} \lambda_{1} - \partial_{\tau} c e^{\sigma} \cdot \langle v_{2}, \nabla \rho  \rangle
+ \partial_{\tau} \sigma \langle v_{2}, S_{-1}^{'} \rho \rangle + \langle N(v_{1}), \rho \rangle, \; \text{and}
\end{array}
\label{Eqn:EstDynLambdaOne}
\end{equation}

\begin{align}
| \partial_{\tau} \sigma | + |\partial_{\tau} c| e^{\sigma} \lesssim  \| \vec{\gamma} \|_{\dot{\mathcal{H}}} \cdot
\label{Eqn:EvolPremSigma}
\end{align}

\end{enumerate}

\label{Prop:LinearParamW}
\end{prop}

The next proposition proved in Section \ref{Section:LocalExistGround} shows that a solution $\vec{u}$ of (\ref{Eqn:NlkgCrit}) with data $\vec{u}(0)$ that is closed enough to $\mathcal{S}$ has a maximal time interval of existence is larger or equal to an interval of size roughly equal to $e^{-\sigma(0)}$; moreover, on this interval,
the distance to $\mathcal{S}$ does not vary much:

\begin{prop}

Let $ 0 \leq \delta_{l} \ll \delta_{s} $ and $ 0 < c_{l} \ll 1 $. Let $\vec{u}$ be a solution of (\ref{Eqn:NlkgCrit}) such that
 $ E ( \vec{u} ) \leq E_{wa}(\vec{W}) + c_{l} d^{2}_{\mathcal{S}} \left( \vec{u}(0) \right)  $ and such that $d_{\mathcal{S}} \left( \vec{u}(0) \right) \leq \delta_{l} $. There exists $1 \gg \bar{c}_{l} > 0$ such that $ \left[ -\bar{c}_{l} e^{-\sigma(0)} , \bar{c}_{l} e^{-\sigma(0)} \right] \subset I(u) $  and there exists a constant $C_{l} \gtrsim 1$ such that

\begin{equation}
\begin{array}{l}
t \in \left[ - \bar{c}_{l} e^{-\sigma(0)}  ,  \bar{c}_{l} e^{-\sigma(0)} \right]: \; \frac{1}{C_{l}} d_{\mathcal{S}} \left( \vec{u}(0) \right) \leq  d_{\mathcal{S}}(\vec{u}(t)) \leq C_{l} d_{\mathcal{S}} \left( \vec{u}(0) \right) \cdot  % \text{and} \; \sigma(t)  = \sigma(0) + O(\delta) \cdot
\end{array}
\nonumber
\end{equation}

\label{Prop:LocalExistGround}
\end{prop}

%\begin{rem}
%Let $\vec{u}$ be a solution of (\ref{Eqn:NlkgCrit} ??? DEF SOL ???. Assume that for some $ T_{+}(u) > T > 0$ we have
%$d_{\mathcal{S}} \left( \vec{u}(t) \right) \leq \delta_{l} $ for $t \in [0, T]$. Hence we see from Proposition
%\ref{Prop:LocalExistGround} that $ T + \bar{c}_{l} e^{-\sigma(T)} < T_{+}(u) $.
%\label{Rem:LocalExistGround}
%\end{rem}

The next proposition, proved in e.g \cite{kriegnakschlagnonrad,nakschlagbook}, establishes some variational estimates:

\begin{prop}%{\textbf{`` Variational Estimates ''} \cite{kriegnakschlagnonrad,nakschlagbook}}
Let $\delta > 0$. Then there exist $\epsilon_{v}:= \epsilon_{v}(\delta) > 0$, $k:= k(\delta) > 0$
and an absolute constant $c>0$ such that if $d_{\mathcal{S}}(\vec{f}) > \delta $, $\bar{\epsilon} \leq \epsilon_{v}$,  and $ E_{wa}(\vec{f}) < E_{wa}(\vec{W}) +  \bar{\epsilon} $, then

\begin{equation}
\begin{array}{ll}
K(f_{1}) & \geq \min \left( k, c \| \nabla f_{1} \|^{2}_{L^{2}} \right)
\end{array}
\nonumber
\end{equation}
or else

\begin{equation}
\begin{array}{ll}
K(f_{1}) & \leq -k
\end{array}
\nonumber
\end{equation}
Here $d_{\mathcal{S}} (\vec{f}) := \inf_{ (\sigma ,c) \in \mathbb{R} \times \mathbb{R}^{N}} \| \vec{f} - T^{c} \vec{S}^{\sigma} \vec{W} \|_{\dot{\mathcal{H}}}$.

\label{Prop:VarEst}
\end{prop}

The next proposition (see \cite{kriegnakschlagnonrad} ) shows that the orthogonal direction $\vec{\gamma}$ of $\vec{v}$ in (\ref{Eqn:Decompv})
can be controlled by the linearized energy $\langle \mathcal{L} \vec{\gamma}, \vec{\gamma} \rangle$:

\begin{prop}%{\textbf{`` Control Of Orthogonal Direction Of $v$ By Linearized Energy''}}

For any $ g \in \dot{H}^{1}$ satisfying $ \langle g, \rho \rangle=0$,  we have

\begin{equation}
\begin{array}{ll}
\| \nabla g \|^{2}_{L^{2}} & \sim \langle  L_{+} g, g \rangle  + |\langle g, S_{0}^{'} \rho \rangle|^{2} +
| \langle  g, \nabla \rho \rangle| ^{2}
\end{array}
\label{Eqn:ControlOrth}
\end{equation}
In particular let $\vec{v} \in \dot{\mathcal{H}}$. Then $ \left\| \vec{\gamma} \right\|_{\dot{\mathcal{H}}}  \sim \langle \mathcal{L} \vec{\gamma},
\vec{\gamma} \rangle $, with $\vec{\gamma}$ defined in the decomposition (\ref{Eqn:Decompvi}).
\label{Prop:ControlOrth}
\end{prop}

In the next proposition, we recall the definition and some properties of the functions $d_{0}$ that was used in the study of the energy-critical
focusing wave equation (see \cite{kriegnakschlagnonrad}):

\begin{prop}

Let $\vec{f} \in \mathcal{V}$. Consider the decomposition (\ref{Eqn:Decompu}) followed by the decomposition (\ref{Eqn:Decompv}). Let $d_{0}$ denote the
following number \footnote{Here $k$ is the number defined in (\ref{Eqn:DefDiscSpect})}:

\begin{equation}
d_{0}^{2}(\vec{f})  := E_{wa}(\vec{f}) - E_{wa}(\vec{W}) + k^{2} \lambda_{1}^{2} \cdot
\nonumber
\end{equation}
Then

\begin{align}
d_{0}^{2}(\vec{f}) \sim  \| \vec{v} \|_{E}^{2} \sim  \| \vec{v} \|^{2}_{\mathcal{H}} \sim d^{2}_{\mathcal{S}} (\vec{f}),
\label{Eqn:Equivdwa}
\end{align}
with $\| \vec{v} \|_{E}$ denoting the linearized energy norm, i.e

\begin{equation}
\| \vec{v} \|^{2}_{E} := \frac{1}{2} (k^{2} \lambda_{1}^{2} + \lambda_{2}^{2}) +
\frac{1}{2} \langle \mathcal{L} \vec{\gamma}, \vec{\gamma} \rangle \cdot
\nonumber
\end{equation}

\label{Prop:DistFuncEigen}
\end{prop}

We now define $\tilde{d}_{\mathcal{S}}$ for the equation (\ref{Eqn:NlkgCrit}) and we recall that for $\vec{\phi}$
such that (\ref{Eqn:CondIneqNrj}) holds, it behaves like the unstable mode:

\begin{prop}

Let $ 0 < \delta^{'}_{s} \ll \delta_{l}$. %Let $ 0 <  \bar{\epsilon} \ll \delta_{l} $.
Let $\chi$ is a smooth, decreasing, and nonnegative function such that $\chi(r) =1 $ if $r \leq 1 $ and $\chi(r)=0$ if
$r \geq 2$.  Let $\vec{\phi} \in \mathcal{H}^{\bar{\epsilon}}$. Then the nonlinear distance function $\tilde{d}_{\mathcal{S}}$ is defined as follows:

\begin{equation}
\tilde{d}_{\mathcal{S}} (\vec{\phi}) := \chi \left(  \frac{d_{\mathcal{S}}(\vec{\phi})}{\delta^{'}_{s}} \right) d_{0}(\vec{\phi}) +
\left( 1 -  \chi \left( \frac{ d_{\mathcal{S}} (\vec{\phi})}{\delta^{'}_{s}} \right) \right) d_{\mathcal{S}}(\vec{\phi}) \cdot
\label{Eqn:NonlinearDistDef}
\end{equation}
Hence if  $ d_{\mathcal{S}}(\vec{\phi}) \ll \delta^{'}_{s} $  and

\begin{equation}
E(\vec{\phi})  < E_{wa}(\vec{W}) + \frac{\tilde{d}^{2}_{\mathcal{S}} (\vec{\phi})}{2},
\label{Eqn:CondIneqNrj}
\end{equation}
then

\begin{equation}
\tilde{d}^{2}_{\mathcal{S}} (\vec{\phi})  \sim  \lambda^{2}_{1} (\vec{\phi})  \cdot
\label{Eqn:DomEigenMode}
\end{equation}

\label{Prop:DistFuncEigenGen}
\end{prop}
The proof of (\ref{Eqn:DomEigenMode}) is short: it follows from (\ref{Eqn:CondIneqNrj}) and the definition of
$d_{0}$.

\begin{rem}
The same conclusion holds if  $ \vec{\phi} \in \dot{\mathcal{H}}^{1}$, $E_{wa}(\vec{\phi}) < E_{wa} (\vec{W}) + \bar{\epsilon}^{2}$, and
$ E_{wa}(\vec{\phi})  < E_{wa}(\vec{W}) + \frac{\tilde{d}^{2}_{\mathcal{S}} (\vec{\phi})}{2} $.
\end{rem}

In the sixth proposition, proved in Section \ref{Section:Ejection}, we study the dynamics of the solution $u$ close to the ground states in the ejection mode:

\begin{prop}%{\textbf{`` Dynamics In The Ejection Mode ''}}

Let $\vec{u}$ be a solution of (\ref{Eqn:NlkgCrit}). Let $t_{0} \in I(u)$. Let $c_{d}$ be such that  $0 < c_{d} \ll c_{l} $.
Let $ R: = \tilde{d}_{\mathcal{S}} \left( \vec{u}(t_{0}) \right)$ and  $ \delta_{f} $ be a positive constant such that $ R \ll \delta_{f} \ll \delta^{'}_{s} $. Assume that $\vec{u}$ satisfies the energy estimate

\begin{equation}
\begin{array}{l}
E(\vec{u}) \leq E_{wa} (\vec{W}) + c_{d} R^{2},
\end{array}
\label{Eqn:EnergyEstInit}
\end{equation}
and the ejection scenario, i.e

\begin{align}
\partial_{t} \tilde{d}^{2}_{\mathcal{S}} \left( \vec{u} \right)(t_{0}) & \geq 0 \cdot
\label{Eqn:EjecScplus}
\end{align}
Then there exists $t_{f}$ such that the following properties hold: $ [t_{0},t_{f}] \subset I(u)$; $ \tilde{d}_{\mathcal{S}} \left( \vec{u}(t_{f}) \right) = \delta_{f}$;
$\tilde{d}_{\mathcal{S}}(\vec{u})$ is increasing on $[t_{0},t_{f}]$; $sign(\lambda_{1})$ is constant on $[t_{0},t_{f}]$;

\begin{equation}
\begin{array}{l}
t \in [t_{0},t_{f}]: \;  \tilde{d}_{\mathcal{S}} \left( \vec{u}(t) \right)  \sim \left| \lambda_{1}(t) \right| \sim R e^{k \tau};
\end{array}
\label{Eqn:Dyndq}
\end{equation}

\begin{equation}
\begin{array}{l}
t \in [t_{0},t_{f}]: \;  \| \vec{\gamma}(t)  \|_{\dot{\mathcal{H}}} + \| u(t) \|_{L^{2}}  \lesssim R +  R^{2} e^{2 k \tau} ;
\end{array}
\label{Eqn:Dynvperp}
\end{equation}

\begin{equation}
\begin{array}{l}
t \in [t_{0},t_{f}]: \; \left| \left( \sigma(t) - \sigma(t_{0}), e^{\sigma(t_{0})} \left( c(t) - c(t_{0}) \right) \right) \right|  \ll 1; \; \text{and}
\end{array}
\label{Eqn:EstVarParam}
\end{equation}
there exists a constant  $C_{K} > 0$ such that

\begin{equation}
\begin{array}{l}
t \in [t_{0},t_{f}]: \; -sign(\lambda_{1}(t)) K(u(t)) \gtrsim \left( e^{k \tau} - C_{K}  \right) R  \cdot
\end{array}
\label{Eqn:DynK}
\end{equation}
Here we consider the decomposition (\ref{Eqn:Decomput}) followed by the decomposition
(\ref{Eqn:Decompv}) applied to $\vec{v}(t)$; we denote by  $\tilde{d}_{\mathcal{S}} (\vec{u} ) $ the map
$ \tilde{d}_{\mathcal{S}} \left( \vec{u} \right):= t \rightarrow \tilde{d}_{\mathcal{S}} \left( \vec{u}(t) \right)$;
we define $\tau:=\tau(t)$ by  $\frac {d \tau}{d t} = e^{\sigma(t)}$ and $\tau(t_{0}):= 0$.

\label{Prop:DynEjection}
\end{prop}

\begin{rem}

By time-reversal invariance \footnote{i.e $ (t,x) \rightarrow \vec{u}(t,x)$ is a solution of (\ref{Eqn:NlkgCrit})
then $(t,x) \rightarrow \vec{u}(-t,x)$} and time-translation invariance
\footnote{if $a \in \mathbb{R}$  and  $ (t,x) \rightarrow \vec{u}(t,x)$ is a solution of (\ref{Eqn:NlkgCrit}) then $(t,x) \rightarrow \vec{u}(t-a,x)$ }
a similar result holds if

\begin{equation}
\begin{array}{ll}
\partial_{t} \tilde{d}^{2}_{\mathcal{S}} \left( \vec{u} \right) (t_{0}) & \leq 0
\end{array}
\label{Eqn:EjecScmin}
\end{equation}
More precisely three exists $t_{f} \leq t_{0}$ such that the following properties hold: $[t_{0},t_{f}] \subset I(u)$;
$ \tilde{d}_{\mathcal{S}} (\vec{u}(t_{f})) = \delta_{f}$; $\tilde{d}_{\mathcal{S}}(\vec{u})$ is decreasing on $[t_{f},t_{0}]$;
$sign(\lambda_{1})$ is constant on $[t_{f},t_{0}]$; (\ref{Eqn:Dyndq}), (\ref{Eqn:Dynvperp}), (\ref{Eqn:EstVarParam}),  and (\ref{Eqn:DynK}) hold with
``$\tau(t)$'' ( resp. ``$[t_{0},t_{f}]$'' ) substituted with ``$-\tau(t)$'' ( resp. ``$[t_{f},t_{0}]$'' ).

\label{Rem:DynEjection}
\end{rem}

In the next proposition, we define a continuous sign function that will decide the fate of the solution of (\ref{Eqn:NlkgCrit}) as
$t \rightarrow T_{+}(u)$ , $t < T_{+}(u)$, or $t \rightarrow T_{-}(u)$, $t > T_{-}(u)$

\begin{prop}%{\textbf{`` Sign function ''}}

Let $\delta_{b} > 0$ be a positive number such that $ \delta_{b} \ll \delta_{f}$. Let $ 0 <  \bar{\epsilon} \ll  \epsilon_{v} (\delta_{b}) $. Let

\begin{equation}
\begin{array}{ll}
\tilde{\mathcal{H}}^{\bar{\epsilon}} := \left\{ \vec{\phi} \in \mathcal{H}^{\bar{\epsilon}}, \,
E(\vec{\phi}) < E_{wa}(\vec{W}) + \frac{\tilde{d}^{2}_{\mathcal{S}} (\vec{\phi})}{2}
\right\}
\end{array}
\nonumber
\end{equation}
We define $\Theta$ in the following fashion:

\begin{equation}
\begin{array}{ll}
\Theta : \tilde{\mathcal{H}}^{\bar{\epsilon}}  &  \rightarrow \left\{ -1, 1 \right\} \\
\vec{\phi} &  \rightarrow \Theta(\vec{\phi}) :=
\left\{
\begin{array}{l}
- sign \left( \lambda_{1} (\vec{\phi}) \right), \,  \tilde{d}_{\mathcal{S}} (\vec{\phi}) \leq \delta_{f}  \\
  sign \left( K(\phi_{1}) \right), \,  \tilde{d}_{\mathcal{S}}(\vec{\phi}) \geq \delta_{b}
\end{array}
\right.
\end{array}
\nonumber
\end{equation}
Then $\Theta$ is continuous. Here in the region $ \left\{ \tilde{d}_{\mathcal{S}} (\vec{\phi}) \leq \delta_{f} \right\}$  we apply the decomposition (\ref{Eqn:Decompvi})
to $\vec{\phi}$.

\label{Prop:SignProp}
\end{prop}

The proof is essentially well-known: see \cite{nakschlagkg} in which the global dynamics above the ground state (energy) of an energy-subcritical
nonlinear Klein-Gordon equation with radial data was studied. See in particular Section entitled ``Sign function away from the ground states'' in this paper.  We write down the proof for convenience of the reader in Section \ref{Sec:SignTheta}.

%\begin{prop}%{\textbf{`` Sign function ''}}
%Let $\delta_{1}$ and $\delta_{2}$  be two positive numbers such that $\delta_{1} < \delta_{2} \ll \delta_{f}$. Let $\epsilon \ll \epsilon_{v} (\delta_{1})$.
%Let

%\begin{equation}
%\begin{array}{ll}
%\tilde{\mathcal{H}_{\epsilon}} := \left\{ \vec{\phi} \in \mathcal{H}^{\epsilon}, \, E(\vec{\phi}) < E_{wa}(\vec{W}) + \frac{ \tilde{d}^{2}_{\mathcal{S}} (\vec{\phi})}{2}
%\right\}
%\end{array}
%\nonumber
%\end{equation}
%We define $\Theta_{\delta_{1}, \delta_{2}}$ in the following fashion:

%\begin{equation}
%\begin{array}{ll}
%\Theta_{\delta_{1}, \delta_{2}} : \tilde{\mathcal{H}_{\epsilon}}  &  \rightarrow \left\{ -1, 1 \right\} \\
%\phi &  \rightarrow \Theta_{\delta_{1}, \delta_{2}} (\vec{\phi}) :=
%\left\{
%\begin{array}{l}
%-sign\left( \lambda_{1} (\vec{\phi}) \right), \,  \tilde{d}_{\mathcal{S}} (\vec{\phi}) \leq \delta_{2}  \\
%sign \left( K(\phi) \right), \,  \tilde{d}_{\mathcal{S}}(\vec{\phi}) \geq \delta_{1}
%\end{array}
%\right.
%\end{array}
%\end{equation}
%Then $\Theta_{\delta_{1},\delta_{2}}$ is continuous.

%\label{Prop:SignProprPrev}
%\end{prop}
%The proof is a straightforward modification of that in \cite{nakashtriroyschrod} and therefore omitted.\\
%Theorem \ref{Thm:Main} $(1)$ is a direct consequence of Proposition \ref{Prop:OnePassLemma}, proved in Section \ref{Section:OnePassLemma}, Remark \ref{Rem:OnePassLemma}, %Proposition \ref{Prop:FarFromGd} proved in Section \ref{Section:FarFromGd}, and Remark \ref{Rem:FarFromGd}: we
%refer to \cite{nakashtriroyschrod} for more details.

\begin{prop}

Let $u$ be a solution of (\ref{Eqn:NlkgCrit}) such that $\vec{u}(0) \in \mathcal{H}^{\bar{\epsilon}}$ with

\begin{equation}
\begin{array}{l}
\bar{\epsilon} \ll \min \left(  \delta_{b}, \epsilon_{v}(\delta_{V}), k(\delta_{V}) \right) \cdot
\end{array}
\label{Eqn:BarEpsProp}
\end{equation}
Assume that there exist $t_{1}, t_{2} \in I(u)$ such that $t_{1} < t_{2}$,
$\tilde{d}_{\mathcal{S}} \left( \vec{u}(t_{1}) \right) < \delta_{b} = \tilde{d}_{\mathcal{S}} \left( \vec{u}(t_{2})\right) $. Then for all
$ t \in I(u) \cap [t_{2}, \infty)$ we have $ \tilde{d}_{\mathcal{S}} \left( \vec{u}(t) \right) \geq \delta_{b}$.
\label{Prop:OnePassLemma}
\end{prop}

\begin{rem}
By time-reversal symmetry the following statement is also true: if there exist $t_{1}, t_{2} \in I(u)$ such that $t_{2} < t_{1}$ and
$\tilde{d}_{\mathcal{S}} \left( \vec{u}(t_{1}) \right) < \delta_{b} = \tilde{d}_{\mathcal{S}} \left( \vec{u}(t_{2})\right) $, then
for all $t \in I(u) \cap (-\infty, t_{2} ]$ we have $ \tilde{d}_{\mathcal{S}} \left( \vec{u}(t) \right) \geq \delta_{b}$.
\label{Rem:OnePassLemma}
\end{rem}

\begin{prop}

Let $u$ be a solution of (\ref{Eqn:NlkgCrit}) such that $\vec{u}(0) \in \mathcal{H}^{\bar{\epsilon}}$ with $\bar{\epsilon}$ satisfying (\ref{Eqn:BarEpsProp}) and
$\bar{\epsilon} \ll \epsilon_{v} (\delta_{b})$. Assume that there exists $t_{0} \in I(u)$ such that $ \tilde{d}_{\mathcal{S}} \left(\vec{u}(t) \right) \geq \delta_{b}$ for all
$t \in [t_{0}, \infty) \cap I(u)$. If $ \Theta \left( u(t_{0}) \right) = 1 $ then $T_{+}(u) = \infty $ and $u$ scatters to a free solution as $ t \rightarrow \infty $;  if $ \Theta \left( u(t_{0}) \right) = -1$ then $ T_{+}(u) < \infty $.

\label{Prop:FarFromGdPr}
\end{prop}

\begin{rem}
By time-reversal symmetry the following statement is also true: assume that there exists $t_{0} \in I(u)$ such that
$ \tilde{d}_{\mathcal{S}} \left(\vec{u}(t) \right) \geq \delta_{b}$ for $t \in (-\infty, t_{0}] \cap I(u)$ ; if $\Theta \left( u(t_{0}) \right) = 1 $
then $ T_{-}(u) = \infty $ and $u$ scatters to a free solution as $t \rightarrow -\infty$; if
$\Theta \left( u(t_{0}) \right) = -1 $ then $ T_{-}(u)= \infty $.
\label{Rem:FarFromGd}
\end{rem}

With these propositions in mind, we can now prove Theorem \ref{Thm:Main}. \\
\\
Let $0 < \epsilon_{*} \ll \min \left(\delta_{b} , k(\delta_{V}), \epsilon_{v}(\delta_{b}) \right)$. \\
Let  $ \mathcal{T} := \left\{ \vec{\psi} \in \mathcal{H}^{\epsilon_{*}}: \; E(\vec{\psi}) \geq E_{wa} (\vec{W}) + c_{d}
\tilde{d}_{\mathcal{S}}^{2} \left( \vec{\psi} \right) \right\} $ \footnote{There are other choices for $\mathcal{T}$: one could also have chosen
for example $\mathcal{T} := \left\{ \vec{\psi} \in \mathcal{H}^{\epsilon_{*}}: \; \tilde{d}_{\mathcal{S}} \left( \vec{\psi} \right) \leq (c_{d})^{-\frac{1}{2}}
\epsilon_{*}   \right\} $}.  \\
Proposition \ref{Prop:DistFuncEigen} implies that  $\tilde{d}_{\mathcal{S}}^{2} \left( \vec{\psi}  \right) \sim d_{\mathcal{S}}^{2} \left( \vec{\psi} \right)$: hence
$ \mathcal{T} \subset \mathcal{B}^{O(\epsilon_{*})} \cap \mathcal{H}^{\epsilon_{*}} $. \\
If  $\vec{u}(t) \in \mathcal{T}$ for all $t \in I(u)$ then the conclusions of Theorem \ref{Thm:Main}, $(1)$ and $(2)$, hold. \\
If not there exists $ t^{'} \in I(u) $ such that  $ E(\vec{u}) \leq E_{wa} (\vec{W}) +  c_{d} \tilde{d}_{\mathcal{S}}^{2} \left( \vec{u}(t^{'}) \right) $. Assume that
$\partial_{t} \tilde{d}_{\mathcal{S}} \left( \vec{u}(t') \right) \geq 0 $ (resp. $\partial_{t} \tilde{d}_{\mathcal{S}} \left( \vec{u}(t') \right) < 0 $);
then we see by applying Proposition \ref{Prop:DynEjection} that $\tilde{d}_{\mathcal{S}} (\vec{u}) $ increases (resp. decreases) from $t'$ (resp. $t^{''}$) to
 $t^{''}$ (resp. $t^{'}$) with $t^{''}$ such that $ \tilde{d}_{\mathcal{S}} \left( \vec{u}(t^{''} ) \right) = \delta_{b} $; hence
by applying Proposition \ref{Prop:OnePassLemma} and Proposition \ref{Prop:SignProp} we see that either $(f+)$ or $(f-)$ holds, with
$(f \pm): \, \Theta (\vec{u}(t)) = \pm  1 \, \text{for all} \, t \in [t^{'}, \infty) \cap I(u)$ (resp. $(b+)$ or $(b-)$ holds, with
$(b \pm): \, \Theta (\vec{u}(t)) = \pm  1 \, \text{for all} \, t \in (-\infty, t^{'}] \cap I(u) $); moreover by Proposition \ref{Prop:FarFromGdPr} we see that if
$(f+)$ holds  (resp. $(f-)$ holds ) then $u$ scatters to a free solution as  $t \rightarrow \infty $
(resp. $t \rightarrow -\infty$) and if $(b+)$ holds (resp. $(b-)$ holds) then
 $ [t^{'}, \infty) \cap I(u) $ (resp. $ (-\infty, t^{'}] \cap I(u)$) is finite. Hence the conclusions of Theorem \ref{Thm:Main}, $(1)$ and $(2)$, also hold. \\
Theorem \ref{Thm:Main}, $(3)$ is proved in Section \ref{Section:ProofThmMain2}.

\section{Proof Of Proposition \ref{Prop:LinearParamW}}
\label{Section:LinearParamW}

In this section we prove Proposition \ref{Prop:LinearParamW}. \\

\begin{align*}
\partial_{t}  \vec{u} & = \mathcal{J} D \vec{u} +
\left(
\begin{array}{l}
0 \\
-u
\end{array}
\right) +
\left(
\begin{array}{l}
0 \\
|u|^{2^{*}-2} u
\end{array}
\right)
\end{align*}

Noting that (with $\vec{h}:= (h_{1},h_{2})$)

\begin{align*}
\partial_{t} (T_{c} \vec{h}) = T_{c} ( \partial_{t} \vec{h} - \partial_{t} c \cdot \nabla  \vec{h}), & \qquad
\partial_{t} ( \vec{S}^{\sigma} \vec{h}) = \vec{S}^{\sigma} (\partial_{t} \vec{h} +  \partial_{t}  \sigma \vec{\Lambda} \vec{h} ) \\
\triangle S_{-1}^{\sigma} h_{1}  = e^{2 \sigma} S_{-1}^{\sigma} \triangle h_{1}, & \qquad | S_{-1}^{\sigma} h_{1} |^{2^{*} -2}  S_{-1}^{\sigma} h_{1}  =
e^{2 \sigma} S_{-1}^{\sigma}
( |h_{1}|^{2^{*} -2} h_{1} ) \end{align*}
we see ( after composition by $T^{-c} \vec{S}^{-\sigma} $)

\begin{equation*}
\partial_{t} \vec{v}   =  \left( e^{\sigma} \partial_{t} c \cdot \nabla  - \partial_{t} \sigma  \vec{\Lambda} \right) ( \vec{W} + \vec{v} )
+ e^{\sigma} \left( \mathcal{J} \mathcal{L} \vec{v} + \underline{N}(\vec{v})  \right)
-e^{- \sigma}  \left(
\begin{array}{l}
0 \\
W + v_{1}
\end{array}
\right)
\label{Eqn:Linearizationaux}
\end{equation*}
Therefore (\ref{Eqn:Linearization}) holds. \\
Hence, differentiating with respect to $\tau$ the relation $\lambda_{\pm} := \omega \left( \vec{v}, \vec{g}^{^{\mp}} \right) $, one sees that (\ref{Eqn:EstDynLambdaPlus}) holds.
Hence, differentiating with respect to $\tau$ (\ref{Eqn:Orth1}) and (\ref{Eqn:Orth2}), we get, by integration by part and
(\ref{Eqn:AdjointS}), (\ref{Eqn:EvolPremSigma}): see \cite{kriegnakschlagnonrad} for more details.

%Notice also from (\ref{Eqn:Const}) and (\ref{Eqn:EstRes}) that

%\begin{align}
%e^{-2 \sigma} \lesssim \| u \|^{2}_{L^{2}} \cdot
%\label{Eqn:BoundSigma}
%\end{align}s

\section{Proof of Proposition \ref{Prop:LocalExistGround}}
\label{Section:LocalExistGround}

In this section we prove Proposition  \ref{Prop:LocalExistGround}. \\
\\
By using the time-reversal invariance \footnote{i.e if $u$ is a solution of (\ref{Eqn:NlkgCrit}) then $(t,x) \rightarrow u(-t,x)$ is also a solution
of (\ref{Eqn:NlkgCrit})} it suffices to prove the proposition  by replacing `` $ \left[ -\bar{c}_{l} e^{-\sigma(0)}, \bar{c}_{l} e^{-\sigma(0}) \right] $ '' with
`` $ \left[ 0, \bar{c}_{l} e^{-\sigma(0)} \right] $ '' in its statement. \\
\\
Let $\mathcal{O}$ be the largest interval of the form $ [0, a )$  such that $\vec{u}(t) \in \mathcal{V}$, $t \in [0,a)$. We apply the decomposition (\ref{Eqn:Decomput}) for all $t \in \mathcal{O}$. Let $C^{'}_{l} \gtrsim 1  $ be a constant such that all the statements below and in this section are true. Let $ \vec{w} (t) := T^{c(t)} \vec{S}^{\sigma(t)} \vec{v}(t) = \vec{u}(t) - T^{c(t)} \vec{S}^{\sigma(t)}  \vec{W} $. We write $\vec{w}(t) := \left( w_{1}(t), w_{2}(t) \right)$. If $J$ is an interval let $ Q(J) := \| \vec{w} \|_{ L_{t}^{\infty} \dot{\mathcal{H}} (J) } + \| w_{1} \|_{L_{t}^{\frac{d+2}{d-2}} L_{x}^{\frac{2(d+2)}{d-2}} (J)} $.  Let $\mu_{max} :=
\max \left\{ \mu \in \mathbb{R}^{+}: \, Q \left( \left[ 0, \mu e^{-\sigma(0)} \right] \right) \leq C^{'}_{l} d_{\mathcal{S}} \left( \vec{u}(0) \right)  \right\} $.  \\
\\
Assume that $\mu_{max} < \bar{c}_{l}$. We will show that it is impossible. \\
\\
Let $ J := [0,.] \subset \left[ 0, \mu_{max} e^{-\sigma(0)} \right] $. Then $\bar{v}$ satisfies

\begin{equation}
\begin{array}{l}
\partial_{tt} w_{1} - \triangle w_{1} = - u + \left| \bar{v} + T^{c(t)}S_{-1}^{\sigma(t)} W \right|^{2^{*} - 2 }
\left( w_{1} + T^{c(t)} S_{-1}^{\sigma(t)} W \right) -
| T^{c(t)} S_{-1}^{\sigma(t)} W|^{2^{*}-2} T^{c(t)} S_{-1}^{\sigma(t)} W \cdot
\end{array}
\nonumber
\end{equation}
We have

\begin{equation}
\begin{array}{l}
Q(J) \lesssim  \| \overrightarrow{\bar{v}}(0) \|_{\dot{\mathcal{H}}} +  \| u \|_{L_{t}^{1} L_{x}^{2} (J)} + \| X \|_{L_{t}^{1} L_{x}^{2} (J)}, \; \text{with}  \\
\end{array}
\nonumber
\end{equation}
$ X := |\bar{v}(t) + S_{-1}^{\sigma(t)} W|^{2^{*} -2}  \left( \bar{v}(t) + S_{-1}^{\sigma(t)} W \right)  - |S_{-1}^{\sigma(t)} W|^{2^{*}-2} S_{-1}^{\sigma(t)} W $. \\
\\
The fundamental theorem of calculus yields  $ |X| \lesssim  \left( | S_{-1}^{\sigma(t)} W |^{2^{*}-2}  + | \bar{v}(t)|^{2^{*} -2} \right)
\left| \bar{v}(t) \right| $ . Hence $ \| X \|_{L_{t}^{1} L_{x}^{2} (J)} \lesssim  \| A \|_{L_{t}^{1} L_{x}^{2}(J)}  + \| B \|_{L_{t}^{1} L_{x}^{2}(J)} $, with \\ $ (A,B) := \left( |\bar{v}(t)|^{2^{*} -2} \bar{v}(t),  \left| S_{-1}^{\sigma(t)} W  \right|^{2^{*} -2} \bar{v}(t) \right)  $.  Clearly
$ \| A \|_{L_{t}^{1} L_{x}^{2}(J)}  \lesssim \| \bar{v} \|^{2^{*} -1}_{L_{t}^{\frac{d+2}{d-2}} L_{x}^{\frac{2(d+2)}{d-2}} (J)} \lesssim
d_{\mathcal{S}}^{2^{*}-1} \left( \vec{u}(0) \right) $. We also get
from (\ref{Eqn:EvolPremSigma})

\begin{equation}
\begin{array}{ll}
\| B \|_{L_{t}^{1} L_{x}^{2}(J)} \lesssim  \left\| S_{-1}^{\sigma(t)} W  \right\|^{2^{*}-2}_{L_{t}^{\frac{d+2}{d-2}} L_{x}^{\frac{2(d+2)}{d-2}} (J)}
\| \bar{v} \|_{L_{t}^{\frac{d+2}{d-2}} L_{x}^{\frac{2(d+2)}{d-2}} (J)}  \ll d_{\mathcal{S}}^{2^{*}-1} \left( \vec{u}(0) \right)  \cdot
\end{array}
\nonumber
\end{equation}
We then estimate $\| u \|_{L_{t}^{1} L_{x}^{2}(J)}$. Let $t \in J$. We have $ E_{wa} \left( \vec{u}(t) \right) = E_{wa}  \left( T^{c(t)} \vec{S}^{\sigma(t)} (\vec{W} + v(t)) \right) = E_{wa} \left( \vec{W} + \vec{v}(t) \right) $. Hence, expanding around $W$, we get from (\ref{Eqn:DefDiscSpect})

\begin{equation}
\begin{array}{ll}
E_{wa} \left( \vec{u}(t) \right) & = E_{wa}(\vec{W}) + \frac{1}{2} \langle \mathcal{L} \vec{v}(t), \vec{v}(t) \rangle  - C \left( \vec{v}(t) \right) \\
& = E_{wa}(\vec{W}) + \frac{1}{2}  \left( \lambda_{2}^{2} - k^{2} \lambda_{1}^{2} \right) + \frac{1}{2} \langle \mathcal{L} \vec{\gamma}(t), \vec{\gamma}(t) \rangle
-C \left( \vec{v}(t) \right) \cdot
\end{array}
\label{Eqn:ExpNrjArW}
\end{equation}
Here

\begin{align}
C (\vec{f})  := \int_{\mathbb{R}^{N}} \frac{ |W + f_{1}|^{2^{*}} - W^{2^{*}} } {2^{*}} - W^{2^{*}-1} f_{1} - \frac{2^{*} -1 }{2} W^{2^{*}-2}
f_{1}^{2} \, dx \cdot
\label{Eqn:DefCf}
\end{align}
Applying elementary estimates for $ \left| v_{1}(t) \right| \gtrsim |W| $, the Taylor formula for $|W| \ll \left| v_{1}(t) \right|$, and the embedding $\dot{H}^{1} \hookrightarrow L^{2^{*}}$, we get $ \left| C(\vec{v}(t)) \right| \lesssim \left\| \vec{v}(t) \right\|^{3}_{\mathcal{H}} + \left\| \vec{v}(t) \right\|^{2^{*}}_{\mathcal{H}}
\lesssim d_{\mathcal{S}}^{3} \left( \vec{u}(0) \right) $. Hence
$E_{wa} \left( \vec{u}(t) \right) = E_{wa}(\vec{W}) + O \left( d_{\mathcal{S}}^{2} \left( \vec{u}(0) \right) \right) $. Since
$ E  \left( \vec{u} \right) \leq  E_{wa} \left( \vec{W} \right) + c_{l} d^{2}_{\mathcal{S}} \left( \vec{u}(0) \right) $, we deduce that

\begin{equation}
\begin{array}{l}
\left\| u(t) \right\|_{L^{2}} \lesssim d_{\mathcal{S}} \left( \vec{u}(0) \right) \cdot
\end{array}
\nonumber
\end{equation}
Let $\chi$ be a radial function such that $\chi(r) \geq 0 $, $\chi(r) = 1 $ if $0 \leq r \leq 1$, and $\chi(r) =0 $ if $r \geq 2$. We get from H\"older inequality

\begin{equation}
\begin{array}{l}
\left( \int_{\mathbb{R}^{d}} |W + v_{1}(t)|^{2} \; dx \right)^{\frac{1}{2}}  \gtrsim \int_{\mathbb{R}^{d}} \left( W + v_{1}(t) \right) \chi(r) \; dx  \gtrsim 1,
\end{array}
\nonumber
\end{equation}
using again at the last step H\"older inequality followed by $\| v_{1}(t) \|_{L^{2^{*}}} \lesssim Q(J)$. The above estimates and
(\ref{Eqn:EvolPremSigma}) yield $  d_{\mathcal{S}} \left( \vec{u}(0) \right) \gtrsim \| u(t) \|_{L^{2}} \sim e^{- \sigma(0)} \left( \int_{\mathbb{R}^{d}} \left| W + v_{1}(t) \right|^{2} \; dx \right)^{\frac{1}{2}} $. Hence $ e^{-\sigma(0)} \lesssim  d_{\mathcal{S}} \left( \vec{u}(0) \right) $ and

\begin{equation}
\begin{array}{l}
\| u \|_{L_{t}^{1} L_{x}^{2} (J)} \lesssim e^{-\sigma(0)} \| u \|_{L_{t}^{\infty} L_{x}^{2}(J)} \lesssim  d_{\mathcal{S}} \left( \vec{u}(0) \right)   \cdot
\end{array}
\nonumber
\end{equation}
The estimates above and a continuity argument show that  $ \left\| \overrightarrow{\bar{v}}  ( \mu_{max} e^{- \sigma(0)} ) \right\|_{\dot{\mathcal{H}}} \leq 2
d_{\mathcal{S}} \left( \vec{u}(0) \right) $, which contradicts $ \left\| \overrightarrow{\bar{v}}  ( \mu_{max} e^{- \sigma(0)} ) \right\|_{\dot{\mathcal{H}}} = C^{'}_{l} d_{\mathcal{S}} \left( \vec{u}(0) \right) $. \\
\\
Hence using Proposition \ref{Prop:DistFuncEigen}, we see that $d_{\mathcal{S}} \left( \vec{u}(t) \right) \lesssim
d_{\mathcal{S}} \left( \vec{u}(0) \right) $. \\
\\
Let $t_{0} \in [0, c_{l} e^{-\sigma(0)} ]$. Since $ v(t) := u ( t_{0} - t ) $ is a solution of (\ref{Eqn:NlkgCrit}) that satisfies
$d_{\mathcal{S}} \left( \vec{v}(0) \right) \leq C_{l} \delta_{l} $ we see from applying the above result to $v$ that  $ d_{\mathcal{S}} \left( \vec{u} (t_{0}) \right) \lesssim d_{\mathcal{S}} \left( \vec{u}(0) \right) $.

\section{Proof Of Proposition \ref{Prop:DynEjection}}
\label{Section:Ejection}

In this section we prove Proposition \ref{Prop:DynEjection} . \\
\\
Let $\tilde{C} > 0$ be a fixed constant large enough such that all the statements below are true. \\
Let

\begin{equation}
\bar{t} := \sup \left\{ t \geq t_{0}: \;
\begin{array}{l}
\frac{1}{\tilde{C}} \lambda_{1}(t_{0}) e^{ k \tau} \leq  \lambda_{1}(t) \leq \tilde{C} \lambda_{1}(t_{0}) e^{k \tau} \; \; (A), \\
\\
\tilde{d}_{\mathcal{S}}^{2}  \left( \vec{u} \right) \, \text{is} \, \nearrow \, \text{on} \, [t_{0},t] \; \; (B),  \\
\\
\| u(t) \|_{L^{2}} + \left\| \vec{\gamma}(t) \right\|_{\dot{\mathcal{H}}} \leq \tilde{C} \left( R + R^{2} e^{2 k \tau} \right) \; \; (C) \cdot
\end{array}
\right\} \cdot
\nonumber
\end{equation}
Let $\tilde{t}$ be such that $\delta^{'}_{s} \gg R e^{ \tilde{C}^{-1}  e^{\sigma(t_{0})} (\tilde{t} - t_{0}) } \gg \delta_{f} $. \\
\\
\underline{Claim}: $ \bar{t} > t_{0} + \tilde{t} $. \\
\\
Assuming that the claim holds, then by Proposition \ref{Prop:DistFuncEigen}  and (\ref{Eqn:EvolPremSigma}), we see that (\ref{Eqn:Dyndq}), (\ref{Eqn:Dynvperp}), and (\ref{Eqn:EstVarParam}) hold. From the invariance of $K$ by the scaling transform $S_{-1}^{\sigma(t)}$ and by the translation transform $T^{c(t)}$  we get from plugging (\ref{Eqn:Decomput}) and Taylor expansion in the region $\{ |v_{1}(t)| \ll |W| \}$

\begin{align*}
K(u(t)) = - (2^{*} -2) \langle W^{2^{*}-1}, v_{1}(t) \rangle + O (\| v_{1} (t)\|^{2}_{\dot{H}^{1}})
\end{align*}
Since $v_{1}(t) = \lambda_{1}(t) \rho + \gamma_{1}(t)$, we see from (\ref{Eqn:Dyndq}) and  (\ref{Eqn:Dynvperp}) that
(\ref{Eqn:DynK}) holds. \\
\\
It remains to prove the claim. Let $\bar{\tau}$ be such that $ \bar{\tau} := \int_{t_{0}}^{\bar{t}} e^{\sigma(t)} \, dt $. Assume that it fails. Let $ Res(t) := \frac{1}{\sqrt{2k}} e^{- 2 \sigma(t)}  \omega
\left(
\left(
\begin{array}{l}
0 \\
W + v_{1}
\end{array}
\right), \vec{g}^{^{+}} + \vec{g}^{^{-}}
\right)
$. Let $t \in [t_{0}, \bar{t})$. Observe by positivity of $\rho$ and $W$, by (\ref{Eqn:DomEigenMode}), and the embedding $ \dot{H}^{1} \hookrightarrow L^{2^{*}} $,
that

\begin{equation}
\begin{array}{l}
\int_{\mathbb{R}^{d}} \left| W + v_{1}(t) \right| \rho \; dx  \gtrsim  \int_{\mathbb{R}^{d}} W \rho \; dx - \| v_{1}(t) \|_{L^{2^{*}}} \| \rho \|_{L^{(2^{*})^{'}}}
\gtrsim 1 \cdot
\end{array}
\nonumber
\end{equation}
Hence, using also the Cauchy-Schwarz inequality and (\ref{Eqn:Decomput}) we get

\begin{equation}
\begin{array}{l}
| Res(t) | \lesssim e^{-2 \sigma(t)} \left(  \int_{\mathbb{R}^{d}} |W + v_{1}(t)| \rho \; dx \right)^{2} \lesssim \| u(t) \|^{2}_{L^{2}} \lesssim
R^{2} + R^{4} e^{4 k \tau} \cdot
\end{array}
\label{Eqn:EstRes}
\end{equation}
We also see from Proposition \ref{Prop:DistFuncEigen} and Proposition \ref{Prop:DistFuncEigenGen}  that $ \tilde{d}_{\mathcal{S}} (\vec{u})(\tau)  \sim |\lambda_{1}(t_{0})| \sim
\| \overrightarrow{v}(\tau) \|_{\mathcal{H}}$. A computation of the first derivative and the second derivative of $ \tau \rightarrow \tilde{d}_{\mathcal{S}}^{2}(\vec{u}) (\tau) $ combined with (\ref{Eqn:Linearization}), (\ref{Eqn:EstDynLambdaOne}), (\ref{Eqn:EvolPremSigma}), and the above estimates shows that

\begin{equation}
\begin{array}{l}
\partial^{2}_{\tau} \tilde{d}_{\mathcal{S}}^{2}(\vec{u}) =   k^{2}  \lambda^{2}_{1} + o \left( \lambda_{1}^{2} \right)
\end{array}
\cdot
\label{Eqn:EstSndDerivDist}
\end{equation}
In the expression above we also used the estimate $|N(f)| \lesssim |f|^{2} + |f|^{2^{*}-1}$ \footnote{
in the region $ \{ |f| \ll |W| \} $ we use the Taylor-Lagrange formula} combined with the H\"older inequality and
the embedding $\dot{H} \hookrightarrow L^{2^{*}}$. Hence we get for $\tau \in [0, \bar{\tau} )$ $ \partial_{t} \tilde{d}_{\mathcal{S}} (\vec{u}(t)) > 0 $. \\
We see from (\ref{Eqn:EjecScplus}) that  $ sign \left( \lambda_{1} \right) \lambda_{2} (0 ) \gtrsim  - \lambda_{1}^{2} (0)  $. Consequently
$ \lambda_{+}(0) \sim \lambda_{1}(0) $ . Integrating (\ref{Eqn:EstDynLambdaPlus}) and using the above estimates combined with the H\"older inequality
and the embedding $ \dot{H}^{1} \hookrightarrow L^{2^{*}} $, we get

\begin{align}
\left| \lambda_{\pm}(\tau) - e^{\pm k \tau } \lambda_{\pm} (0) \right| \lesssim
\int_{\bar{\tau}}^{\tau} e^{k (\tau -s)} R^{2} e^{ 2 \mu s } \, ds \lesssim R^{2} e^{2 k \tau}
\nonumber
\end{align}
Collecting the estimates above one sees that $sign(\lambda_{1})$ is constant on $[0, \bar{\tau})$  and $(A')$ holds with  $ (A'): \; \lambda_{1} (\tau) \sim \lambda_{1}(t_{0}) e^{k \tau}$. \\
Next we modify an argument used in Lemma $4.3$ of \cite{nakschlagbook} ( see also \cite{nakschlagkg,kriegnakschlagnonrad} ). More precisely we consider the expansion of $\vec{v}$, ignoring the $\gamma$ term. i.e

\begin{equation}
\vec{v}_{\lambda} := \lambda_{+} \vec{g}^{^{+}} + \lambda_{-} \vec{g}^{^{-}} = \vec{v} - \vec{\gamma}
\nonumber
\end{equation}
Expanding around $\vec{W}$ we get

\begin{align*}
E_{wa}(\vec{W} + \vec{v}_{\lambda} ) & =  E_{wa}(\vec{W}) +
\frac{1}{2} ( \lambda_{2}^{2} - k^{2} \lambda^{2}_{1} )  - C( \vec{v}_{\lambda} )
\end{align*}
with $C(\vec{f})$ defined in (\ref{Eqn:DefCf}). Hence using also (\ref{Eqn:ExpNrjArW})

\begin{equation}
E (\vec{u}) - E_{wa}(\vec{W} + \vec{v}_{\lambda}) = \frac{1}{2} \| u \|^{2}_{L^{2}} +
\frac{1}{2} \langle \mathcal{L} \vec{\gamma}, \vec{\gamma} \rangle
+ C (\vec{v}_{\lambda}) - C(\vec{v})
\label{Eqn:ExpNrjExpOrth}
\end{equation}
Let $ X := \| \vec{\gamma} \|_{\dot{\mathcal{H}}} + \| u \|_{L^{2}} $. We compute ( using (\ref{Eqn:EstDynLambdaOne}) and (\ref{Eqn:EvolPremSigma}) )

\begin{equation}
\begin{array}{ll}
\partial_{\tau} E_{wa} ( \vec{W} + \vec{v}_{\lambda} )  & =  \lambda_{2} \partial_{\tau} \lambda_{2} - k^{2} \lambda_{1} \partial_{\tau} \lambda_{1}
- \langle N( v_{\lambda,1}| \partial_{\tau} v_{\lambda,1} \rangle \\
& = \lambda_{2} \left(  k^{2} \lambda_{1} +  O ( \| \vec{\gamma} \|_{\dot{\mathcal{H}}} \| \vec{v} \|_{\mathcal{H}})  +
\langle N(v_{1}) | \rho \rangle  \right)
- k^{2} \lambda_{1} \left(  \lambda_{2} + Res \right)  \\
& - \langle N(v_{\lambda,1}) | \rho  \rangle ( \lambda_{2} + Res )   \\
& = O \left( (R e^{k \tau})^{2} X + R e^{k \tau}  X^{2} \right)
\end{array}
\label{Eqn:Ewapar}
\end{equation}
Hence, using also Proposition \ref{Prop:ControlOrth} and (\ref{Eqn:ExpNrjExpOrth}), we see that for $t' \in [t_{0},t]$

\begin{align*}
| E (\vec{u} (t')) - E_{wa}(\vec{W} + \vec{v}_{\lambda}(t'))|  & \lesssim  | E( \vec{u}(t_{0})) - E_{wa} (\vec{W} + \vec{v}_{\lambda}(t_{0}) ) |
+ | E_{wa}( \vec{W} + \vec{v}_{\lambda} (t')) - E_{wa} ( \vec{W} + \vec{v}_{\lambda} (t_{0}) ) | \\
& \lesssim  R^{2} + (R e^{k \tau})^{2} \| X \|_{L_{t}^{\infty} ([t_{0},t])} + R e^{k \tau} \| X \|^{2}_{L_{t}^{\infty} ([t_{0},t])} \cdot
\end{align*}
Therefore

\begin{align*}
\| X \|_{L_{t}^{\infty} ([t_{0},t])}^{2} & \lesssim \sup_{t' \in [t_{0},t]}  \left( | E(\vec{u}(t')) - E_{wa}(\vec{W} + \vec{v}_{\lambda}(t')) | +
| C (\vec{v}_{\lambda}(t')) - C(\vec{v}(t')) | \right)  \\
& \lesssim  R^{2} + (R e^{k \tau})^{2} \| X \|_{L_{t}^{\infty} ([t_{0},t])} +  R e^{k \tau}  \| X \|_{L_{t}^{\infty} ([t_{0},t])}^{2}
\end{align*}
Hence $(B')$ holds with $ (B'): \;  \| u(t) \|_{L^{2}} + \left\| \vec{\gamma}(t) \right\|_{\dot{\mathcal{H}}} \lesssim  R + R^{2} e^{2 k \tau} $.  \\
Now by applying Proposition \ref{Prop:LocalExistGround} close (from the left) to $\bar{t}$ we see that
$\bar{t} \in I(u)$. Hence  $\bar{t}$ is in fact the max such that $(A)$, $(B)$, and $(C)$ hold; moreover we see from
(\ref{Eqn:EstSndDerivDist}) that $ \partial_{t} \tilde{d}_{\mathcal{S}} (\vec{u}(\bar{t})) > 0 $. Hence contradiction.

% COMMENTS: in previous version text just after ``...one sees that  $\lambda_{1} (\tau) \sim \lambda_{1}(t_{0}) e^{k \tau}$.''
%Integrating (\ref{Eqn:EstDynLambdaOne}) one sees that

%\begin{align*}
%sign(\lambda_{1}) \left( \lambda_{2}(\tau) - \lambda_{2}(\bar{t}) \right) & \gtrsim k^{2} R  (e^{k \tau} - 1 ) + O (R^{2} e^{ 2 k \tau}) \cdot
%\end{align*}
%Hence, using also (\ref{Eqn:ComputDerivd0Star}) and $(Cd1)$, we see that $(Cd2)$ holds on $ [ t_{0},t_{0} + \tilde{t}]$. \\
%It remains to prove that $(Cd3)$ holds on $ [ t_{0},t_{0} + \tilde{t}]$.

\section{Proof of Proposition \ref{Prop:SignProp}}
\label{Sec:SignTheta}

In this section we prove Proposition \ref{Prop:SignProp}. \\
\\
We claim that $ sign \left( \lambda_{1} \right) $ is continuous at  $ \vec{\phi} \in \tilde{\mathcal{H}}^{\bar{\epsilon}} $ such that
$ \tilde{d}_{\mathcal{S}}(\vec{\phi}) \leq \delta_{f} $. Indeed $\lambda_{1}$ is clearly continuous at $\vec{\phi}$: this follows from
Proposition \ref{Prop:OrthDecompu} and the formula of $\lambda_{1}$ in Remark \ref{Rem:Decompv}. Moreover
$ \lambda_{1}(\vec{\phi}) \neq 0 $: if not Proposition \ref{Prop:DistFuncEigenGen} implies that $\tilde{d}_{\mathcal{S}}(\vec{\phi}) = 0$ (and hence
$E(\vec{\phi})= E_{wa}(\vec{W})$: see Proposition \ref{Prop:DistFuncEigen}), which contradicts $\vec{\phi} \in \tilde{\mathcal{H}}^{\bar{\epsilon}} $. \\
We also claim that $ sign( K \circ P  ) $ \footnote{Here $P$ is the projection defined by $P(\vec{\phi}) = \phi_{1}$ } is continuous at $ \vec{\phi} \in \tilde{\mathcal{H}}^{\bar{\epsilon}} $ such that $ \tilde{d}_{\mathcal{S}}(\vec{\phi}) \geq \delta_{b} $. Indeed $K$ is continuous at
$ \vec{\phi} \in \tilde{\mathcal{H}}^{\bar{\epsilon}} $ such that $ \tilde{d}_{\mathcal{S}}(\vec{\phi}) \geq \delta_{b} $. Moreover
$ K(\phi_{1}) \neq 0 $: this follows from Proposition \ref{Prop:VarEst} and Proposition \ref{Prop:DistFuncEigen}. \\
It remains to show that $ sign  \left( \lambda_{1} (\vec{\phi}) \right) = - sign \left( K(\phi_{1}) \right) $ for all
$ \vec{\phi} \in \tilde{\mathcal{H}}^{\bar{\epsilon}} $ such that $ \delta_{b} \leq \tilde{d}_{\mathcal{S}} (\vec{\phi}) \leq \delta_{f} $. Let $\mathcal{O}$ be the
connected component of $ \left\{ \vec{\psi} \in \tilde{\mathcal{H}}^{\bar{\epsilon}}: \; \delta_{b} \leq \tilde{d}_{\mathcal{S}}
(\vec{\psi}) \leq \delta_{f} \right\} $ that contains $\vec{\phi}$. By continuity $sign(\lambda_{1})$ and $sign(K \circ P)$ are constant on $\mathcal{O}$. Let $\vec{u}$ be a solution of (\ref{Eqn:NlkgCrit}) such that $\vec{u}(0) := \vec{\phi}$ and such that (\ref{Eqn:EjecScplus}) holds. Then we can apply
Proposition \ref{Prop:DynEjection} to get $ sign \left(  K(u(t)) \right) s=  - sign \left( \lambda_{1}(t) \right) $ for
$ t \approx e^{- \sigma(0)}   $. We have $ \vec{u}(t) \in \mathcal{O} $ since $ \left\{ \vec{u}(t'), t' \in [0,t] \right\}$ is connected as the image
of a connected set by the continuous function $\vec{u}$. Hence $ sign \left( \lambda_{1}(\vec{\phi}) \right) =
-sign(K(\phi_{1}))$.

\section{ Proof of Proposition \ref{Prop:OnePassLemma}}
\label{Section:OnePassLemma}

In this section we prove Proposition \ref{Prop:OnePassLemma}. \\

\vspace{4mm}

\textit{Setting} \\

By decreasing the value of $t_{2}$ if necessary we may assume WLOG that $\partial_{t} \tilde{d}_{\mathcal{S}} \left( \vec{u}(t_{2}) \right) \geq 0$.
Assume towards a contradiction that $\tilde{d}_{\mathcal{S}} (\vec{u}(t)) \geq  \delta_{b}$ does not hold for all $t \in I(u) \cap [t_{2}, \infty)$. This means that
that there exists $t_{3} > t_{2}$ such that $ \tilde{d}_{\mathcal{S}} \left( \vec{u}(t_{3}) \right) = \delta_{b}$, $\partial_{t} \tilde{d}_{\mathcal{S}} \vec{u}(t_{3}) \leq 0$,
and $\tilde{d}_{\mathcal{S}} \left( \vec{u}(t) \right) \geq \delta_{b}$ for all $t \in [t_{2},t_{3}]$. \\
By space translation invariance, we may assume WLOG that

\begin{eqnarray*}
c(t_{2}) & = & 0
\end{eqnarray*}

\textit{Hyperbolic region, Variational region} \\
\\
We now define the hyperbolic region and the variational region. \\
Let $ \delta_{b} \ll \delta_{V} \ll \delta_{m} \ll \delta_{f}$. Let $s$ be a local minimizer of the
function $t \rightarrow \tilde{d} \left( \vec{u}(t) \right)$ on $ \left\{ t \in [t_{2}, t_{3}]: \; \tilde{d}_{\mathcal{S}} \left( \vec{u}(t) \right) < \delta_{V} \right\}$.
By Proposition \ref{Prop:DynEjection} there exists an interval $I[s] \subset [t_{2},t_{3}] $ such that \footnote{Here $\tau$ is defined by $\frac{d \tau}{dt} := e^{\sigma(t)} $  and $\tau(s) :=0 $ }

\begin{itemize}

\item[$1.$] $\tilde{d}_{\mathcal{S}}  \left( \vec{u}(t) \right) \sim  \tilde{d}_{\mathcal{S}} \left( \vec{u}(t_{s}) \right) e^{k \tau} $

\item[$2.$] $ - sign \left( \lambda_{1}(\tau) \right) K(u(t)) \gtrsim  \left( e^{k \tau} - C_{K} \right) \delta_{b}$

\item[$3.$]  $\tilde{d}_{\mathcal{S}} \left( \vec{u} \left( \partial I [s] \cap (t_{2},t_{3}) \right) \right) = \delta_{m}$.

\end{itemize}

Let $\mathcal{L}$ be the set of these $s$. Note that $\mathcal{L}$ is finite: this follows from the uniform continuity of
$ \tilde{d}_{\mathcal{S}} ( \vec{u} ) $ on $[ t_{2},t_{3} ]$, along with the fact that
$\left[ \delta_{V}, \delta_{f} \right]$ is included in the range of $\tilde{d}_{\mathcal{S}} \left( \vec{u} \right)$
between two consecutive $s$. The hyperbolic interval $I_{H}$ and the variational interval are defined by the formulas below:

\begin{equation}
I_{H} := \bigcup_{s \in \mathcal{L}} I[s], \;  I_{V} := [t_{2},t_{3}] / I_{H} \cdot
\nonumber
\end{equation}
Taking into account (\ref{Eqn:EstVarParam}) we see that

\begin{equation}
\begin{array}{l}
s \neq \{ t_{2},t_{3} \}: \; \left| \, I[s] \,  \right| \gtrsim e^{-\sigma(s)} \log \left( \frac{\delta_{m}}{ \delta_{V} } \right) \\
s \in \left\{ t_{2}, t_{3} \right\}: \;  \left| \, I[s] \, \right| \approx e^{-\sigma(s)} \log \left( \frac{\delta_{m}}{\delta_{b}} \right) \cdot
\end{array}
\nonumber
\end{equation}
Let $ m_{H} := \max \limits_{s \in \mathcal{L}} e^{-\sigma(s)} $: this number will be used in the sequel when we use the virial identity.  \\
\\
Let  $s \neq \{ t_{2},t_{3} \} $. By still applying Proposition \ref{Prop:DynEjection} we see that there exist intervals $I_{l}[s] \subset I_{V}$ and
$I_{r}[s] \subset I_{V}$ such that $I_{l}[s], I_{r}[s] \cap \partial I [s] \neq \{ \emptyset \}$, $I_{l}[s]$ (resp. $I_{r}[s]$) is located at the left
(resp. right) of $I[s]$ and such that

\begin{equation}
\begin{array}{l}
s \neq \{ t_{2},t_{3} \}: \; \left| \, I_{l}[s] \,  \right| + \left| \, I_{r}[s] \,  \right|  \approx e^{-\sigma(s)} \log \left( \frac{\delta_{f}}{ \delta_{m} } \right) \\
s = t_{2}: \;  \left| \, I_{r}[s] \, \right| \approx e^{-\sigma(s)} \log \left( \frac{\delta_{f}}{\delta_{m}} \right) \\
s = t_{3}: \;  \left| \, I_{l}[s] \, \right| \approx e^{-\sigma(s)} \log \left( \frac{\delta_{f}}{\delta_{m}} \right)
\end{array}
\nonumber
\end{equation}
We prove the following estimate: \\
\\
\underline{Result}: Define $\bar{I}[s]$ as follows: $\bar{I}[s] \in \left\{  I_{l}[s], I_{r}[s] \right\} $ if $s \notin \{ t_{2}, t_{3} \} $;
$\bar{I}[s] := I_{r}[s]$ if $s=t_{2}$, and $\bar{I}[s] := I_{l}[s]$ if $s = t_{3}$. Then

\begin{align}
\| u \|_{L_{t}^{\frac{2(d+1)}{d-2}} L_{x}^{\frac{2(d+1)}{d-2}}(\bar{I}[s])} \gtrsim 1  \\
\label{Eqn:LowerBdStrichGlob}
\end{align}

\begin{proof}

We shall only prove that $\| u \|_{L_{t}^{\frac{2(d+1)}{d-2}} L_{x}^{\frac{2(d+1)}{d-2}} (I_{r}[t_{2}])  } \gtrsim 1$: the other cases are treated similarly.
Let $ 0 < \delta \ll \bar{c}_{l} $ be a small constant. Let $\tilde{t} \in I_{r}[t_{2}] $. Observe from ( \ref{Eqn:EvolPremSigma} ) that $ \sigma(\tilde{t}) \approx \sigma(t_{2}) $. Hence, by (\ref{Eqn:StrichWave0}), (\ref{Eqn:Decomput}), and (\ref{Eqn:EvolPremSigma}) we see that there exists a constant $C > 0 $ such that \footnote{Here we also use the fact that $ e^{it D} \left( S_{-1}^{\mu} f \right)  = S_{-1}^{\mu}  \left( e ^{itD} f \right) $ for $\mu \in \mathbb{R}$ and that $ e^{it D} \left( T^{\vec{\mu}} f \right)  = T^{\vec{\mu}}  \left( e^{itD} f \right) $ for $\vec{\mu} \in \mathbb{R}^{d} $ }

\begin{align*}
\left\| \left( e^{i (t- \tilde{t}) D} u(\tilde{t}) , \frac{ e^{i (t- \tilde{t}) D}}{D} \partial_{t} u (\tilde{t}) \right)
\right\|_{ L_{t}^{\frac{2(d+1)}{d-2}}  L_{x}^{\frac{2(d+1)}{d-2}}  \left( [\tilde{t}, \tilde{t} + \delta e^{-\sigma(t_{2})}] \right)} & \lesssim
\| e^{i t D} W \|_{L_{t}^{\frac{2(d+1)}{d-2}}  L_{x}^{\frac{2(d+1)}{d-2}} \left( [0, C \delta] \right)  } +
\| \vec{v} (\tilde{t}) \|_{\dot{\mathcal{H}}}
\\
& \ll 1,
\end{align*}
Moreover, since $\vec{u}(t) \in \mathcal{H}^{\bar{\epsilon}}$ for $t \in  [\tilde{t}, \tilde{t} + \delta e^{-\sigma(t_{2})}]$, we see (using Proposition \ref{Prop:LocalExistGround} and Proposition \ref{Prop:DistFuncEigen} ) that  $ E(\vec{u}) \leq E_{wa}(\vec{W}) + c_{l} d^{2}_{\mathcal{S}} \left(  \vec{u}(\tilde{t}) \right) $. Hence by using similar arguments as those from (\ref{Eqn:ExpNrjArW}) until the end of Section \ref{Section:LocalExistGround} we see that
$ e^{- \sigma(t_{2})}  \sup_{t \in [ \tilde{t}, \tilde{t} + \delta e^{-\sigma(t_{2})} ] } \| u(t) \|_{L^{2}} \lesssim \delta_{f} $.   Hence

\begin{align}
\| u \|_{ L_{t}^{1} L_{x}^{2} \left( [ \tilde{t}, \tilde{t} + \delta e^{-\sigma(t_{2})} ] \right) }  \lesssim e^{- \sigma(t_{2})}  \log \left( \frac{\delta_{f}}{\delta_{m}} \right)
\sup_{t \in [ \tilde{t}, \tilde{t} + \delta e^{-\sigma(t_{2})} ] } \| u(t) \|_{L^{2}}  \lesssim \delta^{0+}_{f}
\label{BoundLt1Lx2u}
\end{align}
Let $I := [\tilde{t}, . ] \subset [ \tilde{t}, \tilde{t} + \delta e^{-\sigma(t_{2})} ]$. By (\ref{Eqn:StrichWave0}) and (\ref{Eqn:NonlinearEstGen})
we get

\begin{equation}
\begin{array}{l}
\| u \|_{ L_{t}^{\frac{2(d+1)}{d-2}} L_{x}^{\frac{2(d+1)}{d-2}} (I)} \\
\lesssim \left\|  u_{l,\bar{t}} \right\|_{ L_{t}^{\frac{2(d+1)}{d-2}} L_{x}^{\frac{2(d+1)}{d-2}} (I) } + \| u \|_{L_{t}^{1} L_{x}^{2} (I) } +
\| u \|^{2^{*}-2}_{L_{t}^{\frac{2(d+1)}{d-1}} L_{x}^{\frac{2(d+1)}{d-1}} (I) }
\| u \|_{ L_{t}^{\frac{2(d+1)}{d-1}} \dot{B}^{\frac{1}{2}}_{\frac{2(d+1)}{d-1},2} (I)}, \; \text{and}
\end{array}
\nonumber
\end{equation}

\begin{equation}
\begin{array}{l}
\| u \|_{ L_{t}^{\frac{2(d+1)}{d-1}} \dot{B}^{\frac{1}{2}}_{\frac{2(d+1)}{d-1},2} (I)} +
\| u \|_{L_{t}^{\frac{2(d+1)}{d-2}-} L_{x}^{\frac{2(d+1)}{d-2}+} (I)}  \\
\lesssim  \| u(\tilde{t}) \|_{\dot{H}^{1}} + \| u \|_{L_{t}^{1} L_{x}^{2}(I)} +
\| u \|^{2^{*}-2}_{L_{t}^{\frac{2(d+1)}{d-1}} L_{x}^{\frac{2(d+1)}{d-1}} (I) } \| u \|_{ L_{t}^{\frac{2(d+1)}{d-1}} \dot{B}^{\frac{1}{2}}_{\frac{2(d+1)}{d-1},2} (I)}  \cdot
\end{array}
\nonumber
\end{equation}
Moreover $ \| u(\tilde{t}) \|_{\dot{H}^{1}} \lesssim \| W \|_{\dot{H}^{1}} \lesssim  1 $: this follows from Proposition \ref{Prop:DistFuncEigen}. A continuity argument shows that $ \| u \|_{ L_{t}^{\frac{2(d+1)}{d-2}-} L_{x}^{\frac{2(d+1)}{d-2}+} \left( [ \tilde{t}, \tilde{t} + \delta e^{-\sigma(t_{2})} ] \right)} \ll 1 $ and that
$ \| u \|_{ L_{t}^{\frac{2(d+1)}{d-1}} \dot{B}^{\frac{1}{2}}_{\frac{2(d+1)}{d-1},2} \left( [ \tilde{t}, \tilde{t} + \delta e^{-\sigma(t_{2})} ] \right) } \lesssim 1 $.
Iterating we get

\begin{align*}
 \| u \|_{ L_{t}^{\frac{2(d+1)}{d-2}-} L_{x}^{\frac{2(d+1)}{d-2}+} \left( I_{r}[t_{2}] \right) } \lesssim
 \log \left( \frac{\delta_{f}}{\delta_{m}} \right) \cdot
\end{align*}
Computations show that $ \| T^{c(t)} S_{-1}^{\sigma(t)} W \|_{ L_{t}^{\frac{2(d+1)}{d-2}-}  L_{x}^{\frac{2(d+1)}{d-2}+} \left( I_{r}[t_{2}] \right)} \sim
\log \left(  \frac{\delta_{f}}{\delta_{m}} \right) $. Hence

\begin{align*}
\| \tilde{v}_{1} \|_{ L_{t}^{\frac{2(d+1)}{d-2}-} L_{x}^{\frac{2(d+1)}{d-2}+} \left( I_{r}[t_{2}] \right) } \lesssim \log \left(  \frac{\delta_{f}}{\delta_{m}} \right), \qquad \tilde{v}_{1} := T^{c} S_{-1}^{\sigma} v_{1} \cdot
\end{align*}
Interpolating with

\begin{align*}
\| \tilde{v}_{1} \|_{L_{t}^{\infty} L_{x}^{1_{2}^{*}} \left( I_{r}[t_{2}] \right)} & \lesssim  \| v_{1}  \|_{L_{t}^{\infty} \dot{H}^{1} \left( I_{r}[t_{2}] \right)} \lesssim \delta_{f},
\end{align*}
we see that $ \| \tilde{v}_{1} \|_{ L_{t}^{\frac{2(d+1)}{d-2}} L_{x}^{\frac{2(d+1)}{d-2}} \left( I_{r}[t_{2}] \right) }  \lesssim  \delta^{0+}_{f}$. Computations show that $ \| T^{c(t)} S_{-1}^{\sigma(t)} W \|_{L_{t}^{\frac{2(d+1)}{d-2}} L_{x}^{\frac{2(d+1)}{d-2}} \left( I_{r}[t_{2}] \right) } \sim  \log \left( \frac{\delta_{f}}{\delta_{m}} \right)$. Hence (\ref{Eqn:LowerBdStrichGlob}) holds.

\end{proof}

The result allows to divide $I_{V}$ into subintervals into subintervals $(I_{j}=[a_{j},b_{j}])_{ 1 \leq j \leq J}$ such that

\begin{align}
\| u \|_{ L_{t}^{\frac{2(d+1)}{d-2}} L_{x}^{\frac{2(d+1)}{d-2}} (I_{j})} \sim  \eta \cdot
\label{Eqn:ConcStrich}
\end{align}
with $\eta$ a constant such that $ 0 <  \eta \ll 1$. \\
 % \label{Eqn:Reformdist}

\vspace{4mm}

\textit{Virial identity}

Before proceeding, we recall the well-known virial identity %proved in ??

\begin{equation}
\begin{array}{lll}
\partial_{t} \langle x \cdot \nabla u  + \frac{d}{2} u , w \partial_{t} u \rangle & = & -K(u(t)) + E_{ext}(u(t)) \\
& = & RHS,
\end{array}
\label{Eqn:VirialIdentNLKG}
\end{equation}
with $w(t,x):= \chi \left( \frac{|x|}{t- t_{2} + m}  \right)$ and

\begin{align*}
E_{ext}(\vec{u}(t)) & =  \int_{|x| \geq t - t_{2} + m } |\partial_{t} u(t)|^{2}  + | \nabla u (t)|^{2} + |u(t)|^{2} + |u(t)|^{2^{*}} \, dx \cdot
\nonumber
\end{align*}
Here $m > 0$ is to be determined;  $\chi$ is a smooth, decreasing, and nonnegative function defined on $\mathbb{R}^{+}$ such that
$\chi(r) = 1$ if $r \leq 1$ and $\chi(r)=0$ if $r \geq 2$. Let $ t^{'} \in \{ t_{2}, t_{3} \} $. We place
$ \nabla \left( T^{c(t')} S_{-1}^{\sigma(t')} W \right) $,  $ \nabla \left( T^{c(t')} S_{-1}^{\sigma(t')} v_{1}  \right)$,
and $ T^{c(t')} S_{0}^{\sigma(t')} v_{2}  $ in $ L^{2} $. We also place $  T^{c(t')} S_{-1}^{\sigma(t')} W $ and
$T^{c(t')} S_{-1}^{\sigma(t')} v_{1}$ in  $L^{2^{*}}$. Since $ \vec{u}(t') \in \mathcal{H}^{\bar{\epsilon}} $ we see that
$ \| u(t') \|_{L^{2}} \ll \delta_{b} $. Hence, using (\ref{Eqn:Decomput}) and (\ref{Eqn:Dynvperp}) we get

\begin{equation}
\left| \left[ \langle w \partial_{t} u , x \cdot \nabla u  + \frac{d}{2} u  \rangle \right]_{t_{2}}^{t_{3}} \right|  \lesssim
\delta_{b} (t_{3} - t_{2} + m)
\label{Eqn:EstBound}
\end{equation}
and

\begin{eqnarray*}
E_{ext}(\vec{u}(t_{2})) & \lesssim & (m e^{\sigma(t_{2})})^{-1} + \delta^{2}_{b}
\end{eqnarray*}
Then by finite speed propagation and the local well-posedness theory  for small %[??]%
data we see that if $ m e^{\sigma(t_{2})} \gg 1 $  then

\begin{equation}
\begin{array}{l}
E_{ext} (\vec{u}(t)) \lesssim  E_{ext}(\vec{u}(t_{2})) \lesssim  (m e^{\sigma(t_{2})})^{-1} + \delta^{2}_{b}
\end{array}
\label{Eqn:EstExtt}
\end{equation}
By (\ref{Eqn:DynK}),

\begin{equation}
\begin{array}{ll}
\left| \int_{I[t_{2}]} -K(u) + E_{ext}(\vec{u}) \, dt \right|  & \gtrsim e^{-\sigma(t_{2})}
\left(
\delta_{m} - \log \left( \frac{\delta_{m}}{\delta_{b}} \right)
\left( (m e^{ \sigma(t_{2})} ) ^{-1} +  \delta^{2}_{b} \right)  \right) \\
& \gtrsim e^{-\sigma(t_{2})} \delta_{m}, \qquad m e^{\sigma(t_{2})} \gg \delta_{m}^{-1}  \left( \frac{\delta_{m}}{ \delta_{b}} \right)^{0+}
\end{array}
\label{Eqn:EstLow1maxa}
\end{equation}
Similarly if $s \notin  \{ t_{2}, t_{3} \}$ then

\begin{align}
\left| \int_{I[s]}  -K(u) + E_{ext}(\vec{u}) \, dt \right|  & \gtrsim e^{- \sigma (t_{i})} \delta_{m}, \qquad
m e^{\sigma(s)} \gg \delta_{m}^{-1} \left( \frac{\delta_{m}} {\delta_{b,s}} \right)^{0+} \cdot
\label{Eqn:EstLow1maxb}
\end{align}
( Here $\delta_{b,s} := \tilde{d}_{\mathcal{S}} \left( \vec{u}(s) \right)$). We also have

\begin{equation}
\begin{array}{ll}
\left| \int_{I[t_{3}]} -K (\vec{u}) + E_{ext}(\vec{u}) \, dt \right| & \gtrsim e^{- \sigma (t_{3})} \delta_{m}, \qquad  m e^{\sigma(t_{3})} \gg
\delta_{m}^{-1} \left( \frac{\delta_{m}} {\delta_{b}} \right)^{0+} \cdot
\end{array}
\label{Eqn:EstLow1maxc}
\end{equation}
%In the proof of the one-pass lemma, we say that we play with the parameter $\delta_{b}$ if we make it small enough so that all the
%statements are true. Notice that we can always play with $\delta_{b}$ by making the value of $\epsilon$ arbitrarily small: this allows
%to apply Proposition \ref{Prop:DynEjection} and therefore apply the estimates (\ref{Eqn:EstLow1maxa}), (\ref{Eqn:EstLow1maxb}) and
%(\ref{Eqn:EstLow1maxc}).

\subsection{$\mathbf{\Theta(u) = -1}$}
\label{Subsection:OnePassLemmaBlow}

Choose $m$ such thats

\begin{equation}
\frac{\delta_{m}}{ \delta_{b}}    m_{H}  \gg  m \gg  \delta_{m}^{-1}     \left( \frac{\delta_{m}} {\delta_{b}} \right)^{0+}  m_{H}
\label{Eqn:Choicem}
\end{equation}
From Proposition \ref{Prop:VarEst} we get for $t \in I_{V}$ (with $k:=k(\delta_{V})$)

\begin{eqnarray*}
K(u(t)) & \leq & - k,
\end{eqnarray*}
Hence

\begin{align}
\left| \int_{t_{2}}^{t_{3}}  -K(u) + E_{ext}(\vec{u}) \, dt \right| & \gtrsim \frac{1}{\delta_{b}} \delta_{b} |I_{V}|
+ \sum \limits_{\substack{s \in \mathcal{L} \\ s \neq {t_{2},t_{3}}}} \frac{ \delta_{m} }{ \log \left( \frac{\delta_{m}}{ \delta_{b,s}} \right) } |I[s]|
+ \frac{ \delta_{m} }{ \log \left( \frac{\delta_{m}}{ \delta_{b}} \right) }  \left( |I[t_{1}]| + |I[t_{2}]| \right) \\
& + \delta_{m} m_{H} \nonumber \\
& \gg
\min{ \left( \frac{ \frac{\delta_{m}}{\delta_{b}}}
{ \log{ \left( \frac{\delta_{m}}{\delta_{b}}
\right) } }, \, \frac{1}{\delta_{b}}    \right)}
(t_{3} - t_{2}) \delta_{b} + \delta_{m} m_{H} \nonumber \\
& \gg  (t_{3} - t_{3} +  m) \delta_{b},
\label{Eqn:EstLowerBdK}
\end{align}
This is a contradiction.

\subsection{$\mathbf{\Theta(u) = 1}$}
\label{Subsection:OnePassLemmaScatt}

By (\ref{Eqn:RelG})
\begin{equation}
\| \vec{u} (t) \|^{2}_{\mathcal{H}} + \| u(t) \|^{2}_{L^{2}} + \| u(t) \|^{2^{*}}_{L^{2^{*}}} \lesssim 1
\label{Eqn:ApBound}
\end{equation}
Let $  0 < \beta \ll 1$ be a small constant. There are two cases: \\
\\
\textbf{Case $1$}:

\begin{align*}
\int_{I_{V}} \int_{\mathbb{R}^{d}} |\nabla u|^{2} \, dx \, dt & \geq  \beta |I_{V}|
\label{Eqn:LowerBdKinet}
\end{align*}
%By (\ref{Eqn:EstLow1maxa}), (\ref{Eqn:EstLow1maxb}) and (\ref{Eqn:EstLow1maxc})

%\begin{align*}
%\int_{I_{H,a}} K (u) - |E_{ext}(u)| \, dt &  \gtrsim  \frac{ \frac{\delta_{f_{med}}}{\delta_{bd}} }
%{ \log{ \left( \frac{\delta_{f_{med}}}{\delta_{bd}}  \right)   }} \delta_{bd} |I_{H,a}| \\
%\int_{I_{H,i}} K (u) - |E_{ext}(u)| \, dt &  \gtrsim  \frac{ \frac{\delta_{f_{med}}}{\delta_{bd}} }
%{ \log{ \left( \frac{\delta_{f_{med}}}{\delta_{bd}}  \right)   }} \delta_{bd} |I_{H,i}|, \qquad 1 \leq i \leq k \;  and \\
%\int_{I_{H,b}} K (u) - |E_{ext}(u)| \, dt &  \gtrsim  \frac{ \frac{\delta_{f_{med}}}{\delta_{bd}} }
%{ \log{ \left( \frac{\delta_{f_{med}}}{\delta_{bd}}  \right)   }} \delta_{bd} |I_{H,b}|.
%\nonumber
%\end{align*}
Choose $m$ that satisfies  (\ref{Eqn:Choicem}).
From  Proposition \ref{Prop:VarEst} and (\ref{Eqn:ApBound}) we have for $t \in I_{V}$ (with $k:= k(\delta_{V})$)

\begin{align*}
K(u) \gtrsim  \min{( k, c \| \nabla u \|^{2}_{L^{2}} ) }  \gtrsim  \| \nabla u \|^{2}_{L^{2}}
\end{align*}
Hence by % playing again with $\delta_{b}$ and
using  (\ref{Eqn:EstExtt})

\begin{align*}
\int_{I_{V}} -K(u) + E_{ext}(\vec{u}) \, dt & \lesssim  - \beta |I_{V}|
\end{align*}
Hence by a similar argument used in (\ref{Eqn:EstLowerBdK})

\begin{align*}
\left| \int_{t_{2}}^{t_{3}} -K(u) + E_{ext}(\vec{u}) \, dt \right|  & \gg
(t_{3} - t_{2} + m) \delta_{b}
\end{align*}
This is a contradiction. \\
\\
\textbf{Case $2$}:

\begin{eqnarray*}
\int_{I_{V}} \int_{\mathbb{R}^{N}} | \nabla u|^{2} \, dx \, dt & \leq & \beta |I_{V}|
\end{eqnarray*}
Since $K \geq 0$, we can upgrade this inequality to

\begin{equation}
\int_{I_{V}} \int_{\mathbb{R}^{N}} |u|^{2^{*}} \, dx \, dt \leq \beta |I_{V}|
\label{Eqn:DecayEstimate}
\end{equation}
This allows to use an argument in \cite{bourgjams}. Following \cite{nakashimrn}, we see from (\ref{Eqn:ConcStrich}) that there exists
$ M_{j} \in  2^{\mathbb{N}} $ such that for $M \geq M_{j}$ there exist $J_{j} \subset I_{j} $ and $R_{j} \sim \frac{1}{M}$ that satisfy
$ |J_{j}| \sim \frac{1}{M} $ \footnote{Let $C^{'}$ be a large constant such that all the statements below are true. It is showed in \cite{nakashimrn} that there exist $M_{j} \in 2^{\mathbb{N}}$ such that (say) $\| \bar{P}_{\frac{M_{j}}{4}} u \|_{L_{t}^{\infty} L_{x}^{\infty} (I_{j}) } \geq \eta^{C} $  with $ M_{j} |I_{j}| \lesssim 1 $. Observe that the former inequality implies that $\| \bar{P}_{< M} u \|_{L_{t}^{\infty} L_{x}^{\infty} (I_{j}) } \geq  \eta^{C'} $ for $M \geq M_{j}$ since
$ \bar{P}_{< M} \bar{P}_{\frac{M_{j}}{4}}  = \bar{P}_{\frac{M_{j}}{4}} $. Once this observation made, the arguments to get
to get (\ref{Eqn:ConcKinetNrj}) and (\ref{Eqn:ConcPotNrj}) are almost the same as those in \cite{nakashimrn}: therefore they are omitted.}

\begin{equation}
\begin{array}{lll}
\int_{|x-c_{j}| \leq R_{j}}  | \nabla \bar{P}_{< M} u(t) |^{2} + | \bar{P}_{< M} u (t)|^{2} \, dx & \geq & \eta^{C},
\end{array}
\label{Eqn:ConcKinetNrj}
\end{equation}
 and

\begin{equation}
\begin{array}{lll}
\int_{|x-c_{j}| \leq R_{j}} | \bar{P}_{< M} u(t) |^{2^{*}} \, dx & \geq & \eta^{C}
\end{array}
\label{Eqn:ConcPotNrj}
\end{equation}
for $t \in J_{j}=[\tilde{t}_{j},.]$ Here $C \gg 1$ is a well-chosen large positive constant. Hence using (\ref{Eqn:DecayEstimate}) we see that there exists $j_{0} \in [1,..,J]$ such that

\begin{equation}
|J_{j_{0}}| \lesssim  \beta  |I_{j_{0}}|
\label{Eqn:Smallness}
\end{equation}
Without loss of generality we may assume that $\tilde{t}_{j_{0}} - a_{j_{0}} \leq b_{j_{0}} - \tilde{t}_{j_{0}}$.
Let $K := [a_{j}, t ] \subset I_{j} $. Notice that on each $K$, we have, by (\ref{Eqn:ApBound}), (\ref{Eqn:StrichKg0}), and
(\ref{Eqn:NonlinearEstGen})

\begin{equation}
\begin{array}{l}
\| \vec{u} \|_{L_{t}^{\frac{2(d+1)}{d-1}} B^{\frac{1}{2}}_{\frac{2(d+1)}{d-1},2} (K) \times L_{t}^{\frac{2(d+1)}{d-1}} B^{-\frac{1}{2}}_{\frac{2(d+1)}{d-1},2} (K) } \\
\lesssim  \left\| \vec{u}(a_{j}) \right\|_{\mathcal{H}} + \| u \|_{L_{t}^{\frac{2(d+1)}{d-1}} B^{\frac{1}{2}}_{\frac{2(d+1)}{d-1},2} (K)}
\| u \|^{2^{*}-2}_{L_{t}^{\frac{2(d+1)}{d-2}} L_{x}^{\frac{2(d+1)}{d-2}} (K)} \\
\lesssim 1 + \eta^{2^{*}-2}  \| u \|_{L_{t}^{\frac{2(d+1)}{d-1}} B^{\frac{1}{2}}_{\frac{2(d+1)}{d-1},2} (K)}
\end{array}
\nonumber
\end{equation}
Hence a continuity argument shows that

\begin{equation}
\begin{array}{l}
\| \vec{u} \|_{ L_{t}^{\frac{2(d+1)}{d-1}} B^{\frac{1}{2}}_{\frac{2(d+1)}{d-1},2} (I_{j}) \times L_{t}^{\frac{2(d+1)}{d-1}} B^{-\frac{1}{2}}_{\frac{2(d+1)}{d-1},2} (I_{j})}
\lesssim 1 \cdot
\end{array}
\label{Eqn:ControlAuxNorm}
\end{equation}
Hence

\begin{equation}
\begin{array}{l}
\left\| |u|^{2^{*}-2} u \right\|_{L_{t}^{\frac{2(d+1)}{d+3}} L_{x}^{\frac{2(d+1)}{d+3}}(I_{j}) } \lesssim  \eta^{2^{*}-2}
\| u \|_{L_{t}^{\frac{2(d+1)}{d-1}} B^{\frac{1}{2}}_{\frac{2(d+1)}{d-1},2} (I_{j})} \lesssim 1
\end{array}
\label{Eqn:ControlAuxNorm2}
\end{equation}
%\begin{align}
%\| u \|_{(0',I_{j})} + \| u \|_{\left( B , I_{j} \right)} + \| \partial_{t} u \|_{\left( C, I
Next we prove the following estimates:

\begin{lem}
 Let $ 0 < \eta_{2} \ll \eta^{C} $ (with $C$ defined just below (\ref{Eqn:ConcKinetNrj})) and let $B \gg 1 $. Then there exist $ M^{'}_{j_{0}} \sim M_{j_{0}} $ and
 $ \alpha := \alpha(B) \gg 1$ such that for some $ \bar{R}^{'} \leq \frac{\alpha(B)}{M_{j_{0}}} $, $\bar{t}_{j_{0}} \in I_{j_{0}}$,  and $\bar{\bar{t}}_{j_{0}} \in I_{j_{0}}$
 we have

\begin{equation}
\begin{array}{ll}
(A_{1}): \;  \left| \bar{\bar{t}}_{j_{0}} - \bar{t}_{j_{0}}  \right| \geq  B \alpha^{B}(B) \bar{R}^{'}; & (A_{2}): \; \| \bar{P}_{< M^{'}_{j_{0}}} \vec{u}(\bar{t}_{j_{0}}) \|_{\mathcal{H}(|x- c_{j_{0}}| < \bar{R}^{'})}  \gtrsim \eta^{C}; \\
(A_{3}): \;  \| u \|_{ L_{t}^{\frac{2(d+1)}{d-2}} L_{x}^{\frac{2(d+1)}{d-2}} \left( \Isep  \right) } \leq  \eta_{2}; &
(A_{4}): \; \left\| \frac{ \bar{P}_{ < M^{'}_{j_{0}} }  u  (\bar{t}_{j_{0}})}{x- c_{j_{0}}} \right\|_{L^{2}(|x-c_{j_{0}}| \sim \bar{R}^{'})}  \ll  \eta^{C} \cdot
% | \bar{\bar{t}}_{j_{0}} - \bar{t}_{j_{0}} | & \geq  B \langle M \rangle^{ \frac{ \frac{d+1}{d} + \frac{1}{p} }{ \frac{d-1}{N} - \frac{1}{p}} } ( \bar{R}^{'})^{\frac{2}{ %\frac{N-1}{N} - \frac{1}{p}  }} \label{Cond1} \\
\end{array}
\nonumber
\end{equation}

\label{Lem:CondImplBourg}
\end{lem}
The proof of Lemma \ref{Lem:CondImplBourg} is postponed to Subsection \ref{Subsection:CondImplBourg}. We define $\vec{v}:= (v, \partial_{t} v) $ solution of the free Klein-Gordon equation with initial data

\begin{align*}
\vec{v}(\bar{t}_{j_{0}}) := \chi_{\bar{R}^{'}} \bar{P}_{< M^{'}_{j_{0}}} \vec{u}(\bar{t}_{j_{0}}), \qquad
\chi_{\bar{R}'}(x):= \chi \left( \frac{x -c_{j_{0}}}{\bar{R}^{'}} \right)
\end{align*}
(Recall that $\chi$ is defined below (\ref{Eqn:VirialIdentNLKG})). Assume that $ \bar{\bar{t}}_{j_{0}} > \bar{t}_{j_{0}}$ (resp. $ \bar{\bar{t}}_{j_{0}} < \bar{t}_{j_{0}}$ )
Using the well-known dispersive estimate (see e.g \cite{ginebvelo}) we see that there exists $\gamma > 0$ such that for $ t \geq  \bar{\bar{t}}_{j_{0}} $
(resp. $ t \leq  \bar{\bar{t}}_{j_{0}} $ )

\begin{equation}
\begin{array}{l}
\left\| \cos{ \left( (t - \bar{t}_{j_{0}}) \langle D \rangle \right)} v(\bar{t}_{j_{0}}) +
\frac{ \sin{ \left( (t -\bar{t}_{j_{0}}) \langle D \rangle \right)}}{ \langle D \rangle} \partial_{t}  v(\bar{t}_{j_{0}}) \right\|_{B^{0}_{2^{*}, 2}} \\
\lesssim \frac{1}{ | t - \bar{t}_{j_{0}} |^{\frac{d-1}{d}} }
\left(
\left\| v( \bar{t}_{j_{0}} ) \right\|_{B^{\frac{d+1}{d}}_{ (2^{*})^{'} , 2 }}
+ \left\| \partial_{t} v ( \bar{t}_{j_{0}} ) \right\|_{B^{\frac{1}{d}}_{(2^{*})^{'},2}}
\right) \\
\\
\lesssim \frac{1}{ | t - \bar{t}_{j_{0}} |^{\frac{d-1}{d}} }
\left(
\begin{array}{l}
 \| \chi_{\bar{R}^{'}} \|_{B^{\frac{d+1}{d}}_{\frac{d}{2} ,2} }     \| \bar{P}_{< M^{'}_{j_{0}}} u(\bar{t}_{j_{0}}) \|_{L^{2^{*}}}
 + \| \chi_{\bar{R}^{'}} \|_{L^{\frac{d}{2}}}   \| \bar{P}_{< M^{'}_{j_{0}}} u( \bar{t}_{j_{0}} ) \|_{ B^{\frac{d+1}{d}}_{2^{*},2} } \\
 + \| \chi_{\bar{R}^{'}} \|_{B^{\frac{1}{d}}_{\frac{d}{2} ,2} } \| \bar{P}_{< M^{'}_{j_{0}}} \partial_{t} u(\bar{t}_{j_{0}}) \|_{L^{2^{*}}}
 + \| \chi_{\bar{R}^{'}} \|_{L^{\frac{d}{2}}}  \| \bar{P}_{< M^{'}_{j_{0}}} \partial_{t} u (\bar{t}_{j_{0}}) \|_{B^{\frac{1}{d}}_{2^{*},2}}
\end{array}
 \right) \\
\\
\lesssim  B^{-\gamma},
\end{array}
\label{Eqn:DispEst}
\end{equation}
Here at the last step of (\ref{Eqn:DispEst}) we used the definition of the Besov norms that uses norms of the normalized Paley-Littlewood projectors and we combined (\ref{Eqn:ApBound}) with the estimates below (derived from Bernstein-type inequalities)

\begin{equation}
\begin{array}{l}
M \in \{ 2^{\mathbb{N}},0 \},  M \lesssim (\bar{R}{'})^{-1}: \; \| \bar{P}_{M} \chi_{\bar{R}'} \|_{ B^{r}_{\frac{d}{2},2}}  \lesssim M^{r} \| \chi_{\bar{R}'} \|_{L^{\frac{d}{2}}} \lesssim M^{r} (\bar{R}^{'})^{2}, \\
(M,r) \in 2^{\mathbb{N}} \times \left\{ \frac{d+1}{d}, \frac{d}{2} \right\} : \; M \gg (\bar{R}{'})^{-1} : \; \| \bar{P}_{M} \chi_{\bar{R}'} \|_{ B^{r}_{\frac{d}{2},2}}  \lesssim M^{r-3} \| \nabla^{3} \chi_{\bar{R}'} \|_{L^{\frac{d}{2}}} \lesssim M^{r-3} (\bar{R}^{'})^{-1}, \\
 \left\| \bar{P}_{< M^{'}_{j_{0}}} u \left( \bar{t}_{j_{0}} \right) \right\|_{L^{2^{*}}} \lesssim \left\| u \left( \bar{t}_{j_{0}} \right) \right\|_{H^{1}}, \;
 \| \bar{P}_{< M^{'}_{j_{0}}} u( \bar{t}_{j_{0}} ) \|_{ B^{\frac{d+1}{d}}_{2^{*},2} } \lesssim M_{j_{0}}^{\frac{d+1}{d}} \left\| u \left( \bar{t}_{j_{0}} \right) \right\|_{H^{1}}, \\
\left\| \bar{P}_{< M^{'}_{j_{0}}} \partial_{t} u \left( \bar{t}_{j_{0}} \right) \right\|_{L^{2^{*}}} \lesssim  M_{j_{0}}  \left\| \partial_{t} u \left( \bar{t}_{j_{0}} \right) \right\|_{L^{2}}, \; \text{and} \;  \| \bar{P}_{< M^{'}_{j_{0}}} \partial_{t} u (\bar{t}_{j_{0}}) \|_{B^{\frac{1}{d}}_{2^{*},2}} \lesssim
M_{j_{0}}^{\frac{d+1}{d}}  \left\| \partial_{t} u \left( \bar{t}_{j_{0}} \right) \right\|_{L^{2}} \cdot
\end{array}
\nonumber
\end{equation}
Assume  that $ \bar{\bar{t}}_{j_{0}} > \bar{t}_{j_{0}}$ (resp. $ \bar{\bar{t}}_{j_{0}} < \bar{t}_{j_{0}}$ ). A similar argument that uses again well-known dispersive estimates (see e.g \cite{ginebvelo}) and integration w.r.t time shows that if $ J \subset [\bar{\bar{t}}_{j_{0}}, \infty ) $ (resp. $ J \subset (-\infty, \bar{\bar{t}}_{j_{0}}] $ ) then there exists a constant (that we still denote by  $\gamma$) such that for $ t \geq \bar{\bar{t}}_{j_{0}} $ (resp.
$t \leq \bar{\bar{t}}_{j_{0}}$ )

\begin{equation}
\begin{array}{l}
\left\| \cos{ \left( (t - \bar{t}_{j_{0}}) \langle D \rangle \right)} v(\bar{t}_{j_{0}}) +
\frac{ \sin{ \left( (t - \bar{t}_{j_{0}}) \langle D \rangle \right)}}{ \langle D \rangle} \partial_{t}  v(\bar{t}_{j_{0}})
\right\|_{A(J)} \lesssim B^{-\gamma}
\end{array}
\label{Eqn:DispEst2}
\end{equation}
Here $A(J) \in \left\{ L_{t}^{\frac{2(d+1)}{d-2}} L_{x}^{\frac{2(d+1)}{d-2}} (J), L_{t}^{\frac{2(d+1)}{d-1}} B^{\frac{1}{2}}_{\frac{2(d+1)}{d-1},2}(J) \right\} $.
Let $\vec{w} := \vec{u} - \vec{v}$. By (\ref{Eqn:StrichKg0}) there exists $c > 0$ such that

\begin{equation}
\begin{array}{l}
\left\| \vec{w} (\bar{\bar{t}}_{j_{0}}) \right\|_{\mathcal{H}} \lesssim
\left\| \vec{w}(\bar{t}_{j_{0}}) \right\|_{\mathcal{H}} + \| u \|^{2^{*} - 2}_{ L_{t}^{\frac{2(d+1)}{d-2}} L_{x}^{\frac{2(d+1)}{d-2}}
( [\bar{t}_{j_{0}},\bar{\bar{t}}_{j_{0}}] ) } \| u \|_{ L_{t}^{\frac{2(d+1)}{d-1}} B^{\frac{1}{2}}_{\frac{2(d+1)}{d-1},2} (I_{j_{0}}) }
\lesssim  \left\| \vec{w} (\bar{t}_{j_{0}}) \right\|_{\mathcal{H}} + \eta_{2}^{c} \nonumber
\end{array}
\label{Eqn:Estwdoublebar}
\end{equation}
The Plancherel theorem \footnote{More precisely we use $ \int_{\mathbb{R}^{d}} f(x) \bar{g}(x) \; dx = \int_{\mathbb{R}^{d}} \hat{f}(\xi)
\bar{\hat{g}}(\xi) \; d \xi $  with $ f := |u(t)|^{2^{*}-2} u(t) $ and $ g := \partial_{t} u(t) $ } and a Paley-Littlewood decomposition
for $f$ and $g$ show that

\begin{equation}
\begin{array}{l}
\left| \left\| u(\bar{\bar{t}}_{j_{0}}) \right\|^{2^{*}}_{L^{2^{*}}} - \left\| u(\bar{t}_{j_{0}}) \right\|^{2^{*}}_{L^{2^{*}}} \right| \\
=  \left| 2^{*} \int_{\bar{t}_{j_{0}}}^{\bar{\bar{t}}_{j_{0}}} \int_{\mathbb{R}^{d}} |u(t)|^{2^{*} -2} u(t) \partial_{t} u  \, dx \, dt \right| \nonumber \\
\lesssim \| u \|^{2^{*} -2}_{L_{t}^{\frac{2(d+1)}{d-2}} L_{x}^{\frac{2(d+1)}{d-2}} \left( [ \bar{t}_{j_{0}} , \bar{\bar{t}}_{j_{0}} ] \right)}
\| u \|_{L_{t}^{\frac{2(d+1)}{d-1}} B^{\frac{1}{2}}_{\frac{2(d+1)}{d-1},2}  \left( [ \bar{t}_{j_{0}}, \bar{\bar{t}}_{j_{0}} ] \right) }
\| \partial_{t} u \|_{L_{t}^{\frac{2(d+1)}{d-1}} B^{-\frac{1}{2}}_{\frac{2(d+1)}{d-1},2}  \left( [ \bar{t}_{j_{0}}, \bar{\bar{t}}_{j_{0}} ] \right)} \nonumber \\
\lesssim \eta^{c}_{2} \nonumber
\end{array}
\nonumber
\end{equation}
Hence using also $(A_{2})$, $(A_{4})$, and (\ref{Eqn:DispEst}) we get

\begin{equation}
\begin{array}{lll}
E(\vec{w}(\bar{\bar{t}}_{j_{0}})) & \leq  \frac{1}{2} \left\|  \vec{w} ( \bar{\bar{t}}_{j_{0}} )  \right\|^{2}_{\mathcal{H}} -
\frac{1}{2^{*}} \left\| u (\bar{\bar{t}}_{j_{0}}) \right\|_{L^{2^{*}}}^{2^{*}} + o (\eta^{C}) \\
& \\
& \leq  \frac{1}{2} \left\|  \vec{w} ( \bar{t}_{j_{0}} )  \right\|^{2}_{\mathcal{H}}
- \frac{1}{2^{*}} \left\| u (\bar{t}_{j_{0}}) \right\|_{L^{2^{*}}}^{2^{*}} + o(\eta^{C})  \\
& \\
& \leq E(\vec{u}) - \frac{\eta^{C}}{2}
\end{array}
\label{Eqn:EstNrjw}
\end{equation}
Let $\vec{\tilde{w}}$ be the solution of (\ref{Eqn:NlkgCrit}) with data
$ \vec{\tilde{w}}( \bar{\bar{t}}_{j_{0}} ) := \vec{w}( \bar{\bar{t}}_{j_{0}})$.
We claim that $ K(\tilde{w}(\bar{\bar{t}}_{j_{0}})) = K(w(\bar{\bar{t}}_{j_{0}} )) > 0 $. Clearly  from $(A_{2})$ and $(A_{4})$  we see that there exists $c > 0$ such
that $ \| \vec{\tilde{w}} (\bar{t}_{j_{0}}) \|^{2}_{\mathcal{H}}  \leq \| \vec{u} (\bar{t}_{j_{0}}) \|^{2}_{\mathcal{H}} - c \eta^{C}
$. Then, from $\vec{u}(0) \in \mathcal{H}^{\bar{\epsilon}}$, the conservation of energy, $K \left( u(\bar{t}_{j_{0}} ) \right) \geq 0$  and (\ref{Eqn:RelG}), we see that $ \left\| \vec{u}(\bar{t}_{j_{0}}) \right\|^{2}_{\mathcal{H}} <  \| \nabla W \|^{2}_{L^{2}} + \bar{\epsilon}^{2}$. Therefore using also (\ref{Eqn:Estwdoublebar}) we get $ \| \vec{\tilde{w}} (\bar{\bar{t}}_{j_{0}+1}) \|^{2}_{\mathcal{H}} < \| \nabla W \|^{2}_{L^{2}} $ and the claim holds by (\ref{Eqn:Charac3GdState}). Next we write a perturbation lemma:

\begin{lem}
Let $\overrightarrow{w_{1}}$ and $\overrightarrow{w_{2}}$ solutions of (\ref{Eqn:NlkgCrit}). Let $ a \in (- T_{-}(w_{1}), T_{+}(w_{1}))$ and
$I:=[a,b] \subset ( -T_{-}(w_{2}), T_{+}(w_{2})) $. Let $ \nu > 0 $. Assume that there exists a constant $C_{0} < \infty$ such that

\begin{equation}
\begin{array}{l}
\left\| \vec{w}_{2}  \right\|_{\mathcal{H}} + \| w_{2} \|_{L_{t}^{\frac{2(d+1)}{d-2}} L_{x}^{\frac{2(d+1)}{d-2}} (I) }  \leq C_{0} \cdot
\end{array}
\label{Eqn:Boundw2W}
\end{equation}
There exists $ 0 < \epsilon_{0}:= \epsilon_{0}(C_{0},\nu) \ll 1$ such that if for all $j \in \{ 1,2,3 \}$

\begin{equation}
\begin{array}{ll}
\left\| \cos{ \left( (t-a) \langle D \rangle \right)} \left( w_{1}(a) - w_{2}(a) \right)
+ \frac{ \sin{ \left( (t-a) \langle D \rangle  \right) }}{\langle D \rangle} \left( \partial_{t} w_{1}(a) - \partial_{t} w_{2}(a) \right)
\right\|_{X_{j}(I)} & \leq \epsilon_{0},
\end{array}
\label{Eqn:DispFar1}
\end{equation}
then $I \subset \left(-T_{-}(w_{1}), T_{+}(w_{1}) \right)$ and

\begin{equation}
\begin{array}{ll}
\| w_{1} - w_{2} \|_{L_{t}^{\infty} L_{x}^{2^{*}} (I)}  & \leq \nu \cdot
\end{array}
\label{Eqn:Diffw1w2}
\end{equation}
Here $X_{1}(I) := L_{t}^{\frac{2(d+1)}{d-2}} L_{x}^{\frac{2(d+1)}{d-2}} (I) $,  $ X_{2}(I):= L_{t}^{\frac{2(d+1)}{d-1}} B^{\frac{1}{2}}_{\frac{2(d+1)}{d-1},2} (I) $, and
$X_{3}(I) :=  L_{t}^{\infty} B^{0}_{2^{*},2} (I) $.

\label{lem:PerturbBourg}
\end{lem}
The proof of Lemma  \ref{lem:PerturbBourg} is postponed to Subsection \ref{Subsection:Perturb}. \\
\\
Let $\vec{w}_{2}:= \vec{\tilde{w}}$. By (\ref{Eqn:EstNrjw}),
$ K( w_{2}( \bar{\bar{t}}_{j_{0}}) ) = K ( w ( \bar{\bar{t}}_{j_{0}})  ) > 0$ , we see from
\cite{ibramasmnak} that $ \| w_{2} \|_{ L_{t}^{\frac{2(d+1)}{d-2}} L_{x}^{\frac{2(d+1)}{d-2}} (\mathbb{R})} < \infty$. Hence, using also
(\ref{Eqn:RelI}), the Payne-Sattinger argument \footnote{ that is $K(w_{2}(t)) >0$, $t \in I(w_{2}) = \mathbb{R}$ }, and (\ref{Eqn:EstNrjw}) ) that
there exists $c > 0$ such that

\begin{equation}
\begin{array}{ll}
\| w_{2} (t) \|^{2^{*}}_{L^{2^{*}}} & \leq \| W \|^{2^{*}}_{L^{2^{*}}}  - c \eta^{3}
\end{array}
\label{Eqn:Estw2GroundState}
\end{equation}
Notice also that (\ref{Eqn:DispFar1}) is satisfied on $I:=[ \bar{\bar{t}}_{j_{0}} , t_{3}] $ by  (\ref{Eqn:DispEst}) and
(\ref{Eqn:DispEst2}). Consequently (\ref{Eqn:Diffw1w2}) holds. But then we see from (\ref{Eqn:Estw2GroundState}),
(\ref{Eqn:SobIneq}), and Proposition \ref{Prop:DistFuncEigen}, that it contradicts
$\tilde{d}_{\mathcal{S}} (\vec{u}(t_{3})) = \delta_{b}$.

%We turn to estimating $ \left[ \langle \phi_{m} u , i \partial_{r} u \rangle \right]_{t_{1}}^{t_{2}} $ .
%This estimate is more delicate to establish than that in the blow-up case for low dimensions (i.e $N=3,4$) because
%of the lack of integrability of $W$.
%First, notice that by conservation of mass and (\ref{Eqn:Decompu}), we have

%\begin{equation}
%\begin{array}{ll}
%\| W + \Re{(v)} \|_{L^{2}}, \, \| \Im{(v)} \|_{L^{2}} & \lesssim 1
%\end{array}
%\label{Eqn:ImwL2}
%\end{equation}
%By interpolation of this inequality with (\ref{Eqn:SobIneq}) and by (\ref{Eqn:ControlOfH1}) we see that

%\begin{equation}
%\begin{array}{ll}
%\| \Im{(v(t))} \|_{L^{2^{+}}} & \lesssim R^{+}
%\end{array}
%\label{Eqn:ImwL2plus}
%\end{equation}
%if $t=t_{1},t_{2}$.

%\begin{equation}
%\begin{array}{l}
%k - C m^{-p}   > 0 \\
%e^{- 2 \sigma(t_{1})} \left( \delta_{f} - \log{ \left( \frac{\delta_{f}}{ \delta_{V}} \right) } ( m^{-1} + (m e^{\sigma(t_{1})})^{-1})  \right)   \gg
%\max { \left(  m \delta_{V}, m \delta_{V} e^{- \sigma(t_{i})},
%\min{ \left(
%\begin{array}{l}
%m^{\frac{n}{\infty -}} \delta_{V}^{+}  e^{- \sigma(t_{i})},  \\
%m^{1-} e^{- \left( \frac{n-2}{2}- \right) \sigma (t_{i})} \delta_{V}^{+}
%\end{array}
%\right) }
%\right) }
%\end{array}
%\nonumber
%\end{equation}
%if $i= \{1,2 \}$. By (\ref{Eqn:EstSizeI}), (\ref{Eqn:LowerBdK}), (\ref{Eqn:OgawaTsutM}), (\ref{Eqn:ControlOfH1}), (\ref{Eqn:f4mbound}), (\ref{Eqn:ExpK}), %(\ref{Eqn:LowerBdKVarScat}), (\ref{Eqn:X12Est}), (\ref{Eqn:X13Est}) and (\ref{Eqn:X31Est}),
%we see that this leads to a contradiction.

\subsection{Proof of Lemma \ref{Lem:CondImplBourg} }
\label{Subsection:CondImplBourg}

In this subsection we prove Lemma \ref{Lem:CondImplBourg}. \\
\\
We may assume WLOG that $c_{j_{0}} =0$ and $\tilde{t}_{j_{0}} =0$. \\
Given $R' > 0 $, let $\chi_{R^{'}}(t,x) := \chi \left( \frac{|x|}{t + R^{'}} \right)$ with $\chi$ defined below
(\ref{Eqn:VirialIdentNLKG}). Let $\alpha \geq 0$ and let $M_{\alpha} := 2^{\alpha} M_{j_{0}}$.
Let $C^{'} \gg 1$  be a constant large enough such that all the estimates in the sequel of this proof are true. \\
\\
We first prove the following result: \\
\\
\underline{Result}: Let

\begin{equation}
\begin{array}{ll}
X_{1} \left( M_{\alpha} \right) & := \int_{0}^{b_{j_{0}}} \int  \left| \partial_{t} \bar{P}_{< M_{\alpha}} u \right|
\left| |u|^{2^{*}-2} u - \bar{P}_{ < M_{\alpha}} ( |u|^{2^{*}-2} u) \right|  \, dx \, dt, \; \text{and} \\
& \\
X_{2} \left( M_{\alpha} \right) & := \int_{0}^{b_{j_{0}}} \int  \left| \partial_{t} \bar{P}_{< M_{\alpha}} u \right|
\left| |u|^{2^{*}-2} u - | \bar{P}_{< M_{\alpha}} u|^{2^{*} -2} \bar{P}_{< M_{\alpha}} u  \right| \, dx \, dt  \cdot
\end{array}
\nonumber
\end{equation}
Then there exists $2^{C^{'} \eta^{-C}}  M_{j_{0}} \geq M^{'}_{j_{0}} \geq M_{j_{0}} $ such that

\begin{equation}
\begin{array}{l}
X_{1} \left( M^{'}_{j_{0}} \right) + X_{2} \left(  M^{'}_{j_{0}} \right) \ll \eta^{C} \cdot
\end{array}
\label{Eqn:ResX1X2}
\end{equation}

\begin{proof}

 We use a trick in the spirit of \cite{bourgjams}. Let $\bar{\alpha} \geq 1$. Then (\ref{Eqn:ControlAuxNorm}), (\ref{Eqn:ControlAuxNorm2}), and the fundamental theorem of calculus yield

\begin{equation}
\begin{array}{l}
\sum \limits_{\alpha=1}^{\bar{\alpha}}  X_{1} \left( M_{\alpha} \right)   \\
\\
\lesssim \sum \limits_{\alpha=1}^{\bar{\alpha}} \| \partial_{t} \bar{P}_{ < M_{\alpha}} u \|_{ L_{t}^{\frac{2(d+1)}{d-1}} L_{x}^{\frac{2(d+1)}{d-1}} ([0,b_{j_{0}}]) }
\|  \bar{P}_{ \geq M_{\alpha}} (|u|^{2^{*}-2} u)  \|_{ L_{t}^{\frac{2(d+1)}{d+3}}L_{x}^{\frac{2(d+1)}{d+3}} ([0,b_{j_{0}}]) } \\
\\
\lesssim \sum \limits_{\alpha=1}^{\bar{\alpha}}
\sum \limits_{\substack{ (L_{1},L_{2}) \in 2^{\mathbb{N}} \times 2^{\mathbb{N}} \\ L_{1} < M_{\alpha} \leq L_{2} }}
\left( \frac{  \langle L_{1} \rangle }{ \langle L_{2} \rangle} \right)^{\frac{1}{2}}
\langle L_{1} \rangle^{-\frac{1}{2}} \| \bar{P}_{L_{1}} u \|_{ L_{t}^{\frac{2(d+1)}{d-1}} L_{x}^{\frac{2(d+1)}{d-1}} ([0,b_{j_{0}}]) }
\langle L_{2} \rangle^{\frac{1}{2}} \| \bar{P}_{L_{2}} (|u|^{2^{*}-2} u) \|_{ L_{t}^{\frac{2(d+1)}{d+3}} L_{x}^{\frac{2(d+1)}{d+3}} ([0,b_{j_{0}}]) } \\
\\
\lesssim \sum \limits_{ \substack{ (L_{1}, L_{2}) \in  2^{\mathbb{N}} \times 2^{\mathbb{N}} \\ L_{1} \leq L_{2}}} \left( \frac{\langle L_{1} \rangle}{ \langle L_{2} \rangle} \right)^{\frac{1}{2}-}
\langle L_{1} \rangle^{-\frac{1}{2}}
\| \partial_{t} \bar{P}_{L_{1}} u \|_{ L_{t}^{\frac{2(d+1)}{d-1}} L_{x}^{\frac{2(d+1)}{d-1}} ([0,b_{j_{0}}]) }
 \langle L_{2} \rangle^{\frac{1}{2}} \| \bar{P}_{L_{2}} (|u|^{2^{*}-2} u) \|_{ L_{t}^{\frac{2(d+1)}{d+3}} L_{x}^{\frac{2(d+1)}{d+3}} ([0,b_{j_{0}}])} \\
\\
\lesssim \| \partial_{t} u \|_{L_{t}^{\frac{2(d+1)}{d-1}} B^{-\frac{1}{2}}_{\frac{2(d+1)}{d-1}} ([0,b_{j_{0}}])  }
\| |u|^{2^{*}-2}  u \|_{ L_{t}^{\frac{2(d+1)}{d+3}} B^{\frac{1}{2}}_{\frac{2(d+1)}{d+3},2} ([0,b_{j_{0}}])} \\
\\
\lesssim 1, \; \text{and}
\end{array}
\label{Eqn:X1Control}
\end{equation}

%COMMENT
%\begin{align*}
%\left\| \int \int K(x-y) f(y) dy g(x) dx \right\| & \lesssim \left\| \int K(x-y) f(y) dy \right\|_{L^{2}}  \| g \|_{L^{2}} \\
%& \lesssim  \| K \|_{L^{1}}  \| f \|_{L^{2}} \| g \|_{L^{2}}
%\end{align*}

\begin{equation}
\begin{array}{l}
\sum \limits_{\alpha=1}^{\bar{\alpha}} X_{2} \left( M_{\alpha} \right) \\
\\
\lesssim \sum \limits_{\alpha=1}^{\bar{\alpha}} \| \partial_{t} \bar{P}_{< M_{\alpha}} u \|_{ L_{t}^{\frac{2(d+1)}{d-1}} L_{x}^{\frac{2(d+1)}{d-1}} ([0,b_{j_{0}}]) }
\| |u|^{2^{*}-2} u - | P_{ < M_{\alpha}} u|^{2^{*} -2} \bar{P}_{< M_{\alpha}} u  \|_{ L_{t}^{\frac{2(d+1)}{d+3}} L_{x}^{\frac{2(d+1)}{d+3}} ([0,b_{j_{0}}]) } \\
\\
\lesssim \sum \limits_{\alpha=1}^{\bar{\alpha}}  \| \partial_{t} \bar{P}_{< M_{\alpha}} u \|_{L_{t}^{\frac{2(d+1)}{d-1}} L_{x}^{\frac{2(d+1)}{d-1}} ([0,b_{j_{0}}]) }
\| \bar{P}_{ \geq M_{\alpha}} u \|_{L_{t}^{\frac{2(d+1)}{d-1}} L_{x}^{\frac{2(d+1)}{d-1}}  ([0,b_{j_{0}}])}
\| u \|^{2^{*} -2}_{L_{t}^{\frac{2(d+1)}{d-2}} L_{x}^{\frac{2(d+1)}{d-2}} ([0,b_{j_{0}}])} \\
\\
\lesssim \sum \limits_{\alpha=1}^{\bar{\alpha}} \sum \limits_{\substack{ (L_{1},L_{2}) \in 2^{\mathbb{N}} \times 2^{\mathbb{N}} \\ L_{1} < M_{\alpha} \leq L_{2} }}    \left( \frac{ \langle L_{1} \rangle}{ \langle L_{2} \rangle} \right)^{\frac{1}{2}}
\langle L_{1} \rangle^{- \frac{1}{2}} \|  \bar{P}_{L_{1}} u \|_{L_{t}^{\frac{2(d+1)}{d-1}} L_{x}^{\frac{2(d+1)}{d-1}} ([0,b_{j_{0}}])}
\langle L_{2} \rangle^{\frac{1}{2}}  \| \bar{P}_{L_{2}} u \|_{L_{t}^{\frac{2(d+1)}{d-1}} L_{x}^{\frac{2(d+1)}{d-1}}([0,b_{j_{0}}])} \\
\\
\lesssim  \sum \limits_{\substack{ (L_{1},L_{2}) \in 2^{\mathbb{N}} \times 2^{\mathbb{N}} \\ L_{1} \leq L_{2} }}     \left( \frac{ \langle L_{1} \rangle}{ \langle L_{2} \rangle} \right)^{\frac{1}{2}-}
 \langle L_{1} \rangle^{-\frac{1}{2}} \| \partial_{t} \bar{P}_{L_{1}} u \|_{L_{t}^{\frac{2(d+1)}{d-1}} L_{x}^{\frac{2(d+1)}{d-1}} ([0,b_{j_{0}}])}
\langle L_{2} \rangle^{\frac{1}{2}} \| \bar{P}_{L_{2}} u \|_{L_{t}^{\frac{2(d+1)}{d-1}} L_{x}^{\frac{2(d+1)}{d-1}} ([0,b_{j_{0}}]) } \\
\\
\lesssim \| \partial_{t} u \|_{L_{t}^{\frac{2(d+1)}{d-1}} B^{-\frac{1}{2}}_{\frac{2(d+1)}{d-1},2}([0,T_{1}])} \| u \|_{ L_{t}^{\frac{2(d+1)}{d-1}} B^{\frac{1}{2}}_{\frac{2(d+1)}{d-1} ,2} ([0,b_{j_{0}}]) } \\
\lesssim 1 \cdot
\end{array}
\label{Eqn:X2Control}
\end{equation}
Here we also applied the Young inequality for sequences at the last line of (\ref{Eqn:X1Control}) and (\ref{Eqn:X2Control}). The pigeonhole principle implies that there exists $ 2^{C^{'} \eta^{-C}} M_{j_{0}} \geq  M_{\alpha} \geq M_{j_{0}} $ such that  (\ref{Eqn:ResX1X2}) holds with $M_{j_{0}}^{'} := M_{\alpha}$.

\end{proof}

Next, given $0 \leq T_{1} \leq b_{j_{0}} $, we prove some estimates, namely (\ref{Eqn:EstProp}). \\

First we show that  $K \left( \chi_{R'}(0) \bar{P}_{ < M^{'}_{j_{0}}} u(0) \right) \geq 0 $ . \\
Indeed we may assume WLOG that $k \ll 1$. First assume that  $ \left\| \nabla u(0)  \right\|_{L^{2}} \ll k $ with $k$ defined just below (\ref{Eqn:Choicem}). Then by expansion of the gradient, the Hardy inequality, and Plancherel equality we get $ \| \nabla ( \chi_{R'}(0) \bar{P}_{< M^{'}_{j_{0}}} u(0) ) \|_{L^{2}} \lesssim \frac{ \| \bar{P}_{< M^{'}_{j_{0}}} u(0) \|_{L^{2}(|x| \sim R^{'})}}{R^{'}} +\left\| \nabla \bar{P}_{< M^{'}_{j_{0}}} u(0) \right\|_{L^{2}} \ll k $. Hence we see from (\ref{Eqn:SobIneq}) that
 $ K \left( \chi_{R'}(0) \bar{P}_{< M^{'}_{j_{0}}}  u(0) \right)  \geq  k \| \nabla ( \chi_{R'}(0) \bar{P}_{< M^{'}_{j_{0}}} u(0) ) \|^{2}_{L^{2}} $. Now assume that
 $ \| \nabla u(0)  \|_{L^{2}} \gtrsim k $. The estimate $K \left( u(0) \right) \geq k $  that follows from Proposition \ref{Prop:VarEst},
 (\ref{Eqn:RelI}), the equality  $\frac{1}{2} \int_{\mathbb{R}^{d}} |\nabla W|^{2} \; dx = \frac{1}{2} \int_{\mathbb{R}^{d}} |W|^{2^{*}} \; dx $ that comes from $ -\triangle W = W^{2^{*}-1} $, and the estimate $I \left( \chi_{R'}(0) \bar{P}_{< M^{'}_{j_{0}}} u(0) \right) \leq I \left( u(0) \right) $ yield \\
 $ I \left( \chi_{R'}(0)  \bar{P}_{ < M^{'}_{j_{0}}} u(0)  \right) <  I(W) - \frac{k}{2} $.  Hence, using also (\ref{Eqn:SobIneqMax}), we see that there exists $c :=c (\delta_{V}) > 0$ such that

 \begin{equation}
 \begin{array}{ll}
 K \left( \chi_{R'}(0) \bar{P}_{ < M^{'}_{j_{0}}} u(0) \right) & \geq \left\|  \nabla \left( \chi_{R'}(0) \bar{P}_{< M^{'}_{j_{0}}}  u(0) \right)   \right\|^{2}_{L^{2}}
 \left( 1 - C_{*}^{2} \| \chi_{R'}(0) \bar{P}_{ < M^{'}_{j_{0}} } u(0) \|^{2^{*}-2}_{L^{2^{*}}} \right) \\
 & \geq \left\| \nabla \left( \chi_{R'}(0) \bar{P}_{< M^{'}_{j_{0}}}  u(0) \right)  \right\|^{2}_{L^{2}}
 \left( 1 - C_{*}^{2} \left( \| W \|^{2^{*}}_{L^{2^{*}}} - \frac{d k}{2} \right)^{\frac{2^{*}-2}{2^{*}}} \right)  \\
 & \geq c \left\| \nabla \left(   \chi_{R'}(0) \bar{P}_{< M^{'}_{j_{0}}}  u(0) \right)  \right\|^{2}_{L^{2}} \cdot
 \end{array}
 \nonumber
 \end{equation}
Hence in both cases $K \left( \chi_{R'}(0) \bar{P}_{< M^{'}_{j_{0}}}  u(0) \right)  \geq  0$.  \\
Then observe from (\ref{Eqn:RelG}) that

\begin{align*}
E(\chi_{R^{'}}(0), \bar{P}_{< M^{'}_{j_{0}}} \vec{u} (0)) & \geq  \left( \frac{1}{2} - \frac{1}{2^{*}} \right)
\| \chi_{R'}(0) \bar{P}_{< M^{'}_{j_{0}}} \vec{u}(0) \|^{2}_{\mathcal{H}}
\end{align*}
and consequently

\begin{align*}
E \left( \chi_{R'}(0), \bar{P}_{ < M^{'}_{j_{0}}} \vec{u}(0) \right) & \geq \| \bar{P}_{< M^{'}_{j_{0}}} u(0) \|^{2}_{H^{1} ( |x| \leq R' ) } +
O \left( \int_{|x| \sim R'} |\nabla \bar{P}_{< M^{'}_{j_{0}}} u (0)|^{2} \, dx \right)
\end{align*}
Since  $\sum_{Q \in 2^{\mathbb{N}}} \int_{|x| \sim  Q R'}   |\nabla \bar{P}_{< M^{'}_{j_{0}}} u(0)|^{2} \, dx \lesssim
\| u(0) \|^{2}_{H^{1}} $, one can find a sequence $ \left\{  Q_{1,n}  \right\}_{n \geq 1}$ such that
$Q_{1,n} \leq 2^{C^{'} \eta^{-C} n} $ and $ \int_{|x| \sim Q_{1,n} R_{j_{0}}}  |\nabla \bar{P}_{< M^{'}_{j_{0}}} u(0)|^{2} \, dx  \ll \eta^{C} $. Hence (by (\ref{Eqn:ConcKinetNrj}))

\begin{equation}
\begin{array}{ll}
E \left( \chi_{R_{n}}(0), \bar{P}_{ < M^{'}_{j_{0}}} \vec{u}(0) \right) \gtrsim  \eta^{C} & , \, R_{n}:= Q_{1,n} R_{j_{0}}
\nonumber
\end{array}
\nonumber
\end{equation}
Applying  $\bar{P}_{< M^{'}_{j_{0}}}$ to (\ref{Eqn:NlkgCrit}), multiplying the result by $ \chi_{R_{n}} \partial_{t} \bar{P}_{< M^{'}_{j_{0}}} u $ and integrating on the region
$ [0, T_{1}] \times \mathbb{R}^{d}$ we get

\begin{equation}
\begin{array}{ll}
E \left( \chi_{R_{n}}(T_{1}), \bar{P}_{< M^{'}_{j_{0}}} \vec{u}( T_{1} ) \right) &  =  E \left( \chi_{R_{n}}(0), \bar{P}_{< M^{'}_{j_{0}}} \vec{u} (0) \right) +
O \left( X_{1} ( T_{1} ) \right) + O  \left( X_{2} (T_{1}) \right) \\
& + O \left( X_{3}  \left( T_{1}, R_{n} \right) \right),
\end{array}
\nonumber
\end{equation}
with

\begin{equation}
\begin{array}{ll}
X_{3} \left( T_{1},R^{'} \right) & :=  \int_{0}^{T_{1}} \int_{|x| - t \sim R^{'} } \frac{1}{R_{j_{0}}}
\left( |\partial_{t} \bar{P}_{< M^{'}_{j_{0}}} u|^{2} + | \nabla \bar{P}_{< M^{'}_{j_{0}}} u|^{2} + | \bar{P}_{< M^{'}_{j_{0}}} u|^{2} + | \bar{P}_{< M^{'}_{j_{0}}} u|^{2^{*}} \right) \, dx \, dt \cdot
\end{array}
\nonumber
\end{equation}
Since $ \sum \limits_{n \geq 1} X_{3} \left( T_{1}, R_{n} \right) \lesssim \frac{T_{1}}{R_{j_{0}}} $ one can then find a subsequence
$\left\{ Q_{2,n}:= Q_{2,n}(T_{1}) \right\}_{ n \geq 1}$ of $ \left\{ Q_{1,n} \right\}_{n \geq 1}$ and a
sequence $ \{ b_{n} \}_{n \geq 1} $ such that $ 1 \leq b_{n} \leq 2^{C^{'} \eta^{-C}  \frac{T_{1}}{R_{j_{0}}} n } $,

\begin{equation}
Q_{1,b_{n}} \leq Q_{2,n} \leq Q_{1,b_{n+1}}, \; \text{and} \; X_{3} \left( T_{1}, Q_{2,n} M_{j_{0}}^{'} \right)  \ll \eta^{C} \cdot
\nonumber
\end{equation}
Hence, using also (\ref{Eqn:ResX1X2}) we get

\begin{align*}
\left\| \bar{P}_{< M^{'}_{j_{0}}} \vec{u}(T_{1}) \right\|_{ \mathcal{H} (|x| \leq Q_{2,n} R_{j_{0}} + T_{1} )} & \gtrsim \eta^{C} \cdot
\end{align*}
The Hardy inequality and the Plancherel theorem show that $ \sum \limits_{Q \in 2^{\mathbb{N}}} \int_{|x| - T_{1} \sim Q  R_{j_{0}} }
\frac{ | \bar{P}_{< M^{'}_{j_{0}}} u(T_{1}) |^{2}}{|x|^{2}} \, dx \lesssim \| u(T_{1}) \|^{2}_{H^{1}} $. Hence one can find a subsequence $ \left\{ Q_{3,n} := Q_{3,n}  (T_{1}) \right\}_{ n \geq 1}$  of $ \left\{ Q_{2,n}  \right\}_{n \geq 1}$ and a sequence $\{ b_{n} \}_{n \geq 1 }$ such that $ 1 \leq b_{n} \leq 2^{C^{'} \eta^{-C} n} $,

\begin{equation}
\begin{array}{l}
Q_{2,b_{n}} \leq  Q_{3,n} \leq Q_{2,b_{n+1}}, \; \text{and} \; \left\| \frac{ P_{< M^{'}_{j_{0}}} u(T_{1})}{|x|}
\right\|_{L^{2} (|x| - T_{1} \sim Q_{3,n} R_{j_{0}} )} \ll  \eta^{C}
\end{array}
\nonumber
\end{equation}
Hence there exists $\bar{Q}(T_{1})$ (it suffices to choose $\bar{Q}(T_{1}) := Q_{3,1}(T_{1}) $) such that $ \bar{Q} (T_{1}) \leq (C^{'})^{ C^{'} \eta^{-C} (C^{'})^{C^{'} \eta^{-C} \frac{T_{1}}{R_{j_{0}}} (C^{'})^{C^{'} \eta^{-C}} }}$ and

\begin{equation}
\begin{array}{l}
\left\| \bar{P}_{< M^{'}_{j_{0}}} \vec{u}(T_{1}) \right\|_{ \mathcal{H} (|x| \leq \bar{Q}(T_{1}) R_{j_{0}} + T_{1} )}  \gtrsim \eta^{C} , and
 \left\| \frac{ \bar{P}_{< M^{'}_{j_{0}}} u(T_{1})}{|x|} \right\|_{L^{2} (|x| - T_{1} \sim \bar{Q}(T_{1}) R_{j_{0}} )  } \ll  \eta^{C}
\end{array}
\cdot
\label{Eqn:EstProp}
\end{equation}
Next, given $B \gg 1$,  we write an algorithm that will allow us to define $\bar{t}_{j_{0}}$, $\bar{\bar{t}}_{j_{0}}$, $\alpha:= \alpha(B)$, and
$\bar{R}^{'}$. \\
Let $S := 0 $, $R := R_{j_{0}}$, and $T :=  B  \left( M_{j_{0}} R \right)^{B} R $. Then, as long as
$ \| u \|_{L_{t}^{\frac{2(d+1)}{d-2}} L_{x}^{\frac{2(d+1)}{d-2}} (S,T) } > \eta_{2} $,  do the following:
$R := \bar{Q}(T) R_{j_{0}} + T  $, $ S := T $, and  $T := S +  B(M_{j_{0}} R)^{B} R $. Observe that algorithm has a finite number of steps in view of (\ref{Eqn:Smallness}). Moreover $(A1)$, $(A2)$, $(A3)$, and $(A4)$ hold with $ \bar{t}_{j_{0}}  := S $, $\bar{\bar{t}}_{j_{0}} := T$, and $\bar{R}^{'} := \bar{Q}(T) R_{j_{0}} + T $.

% \color{red}
% COMMENTS: PREVIOUS VERSION
% Recall that there exists $C_{0} > 0$ such that $ R_{j_{0}} \leq \frac{C_{0}}{M_{j_{0}}} $. Let $(S,T) := (0, B C_{0}^{B} R_{j_{0}}) $ and
% $ Var := B C_{0}^{B+1} $. If $ \| u \|_{L_{t}^{\frac{2(d+1)}{d-2}} L_{x}^{\frac{2(d+1)}{d-2}} ([S,T])} \leq \eta_{2} $ then we let $ \left( \bar{t}_{j_{0}}, % % \bar{\bar{t}}_{j_{0}} \right) := (S,T) $ and $ (R^{'}_{j_{0}}, \alpha )  := ( R_{j_{0}}, C_{0}) $, and we stop the algorithm. If not we follow the instruction $(A)$ and then the %instruction $(B)$, with  $(A)$ and $(B)$ defined by \\
%$(A)$: let $ (\kappa,S,T ) := \left( \tilde{Q}_{3,1}\left( \frac{Var}{M_{j_{0}}} \right) C_{0} + \frac{Var}{M_{j_{0}}}, T,  S + B \kappa^{B} \left(  Q_{3,1}(S) R_{j_{0}} + S
%\right)  \right) $ and let $ Var := Var + B \kappa^{B} \left( \tilde{Q}_{3,1} \left( \frac{Var}{M_{j_{0}}} \right) C_{0}  + Var \right) $.  \\
%$(B)$: If $ \| u \|_{L_{t}^{\frac{2(d+1)}{d-2}} L_{x}^{\frac{2(d+1)}{d-2}} ([S,T])} \leq \eta_{2} $ then we let $ \left( \bar{t}_{j_{0}}, \bar{\bar{t}}_{j_{0}} \right) := (S,T)$ %and $ \left( R^{'}_{j_{0}} , \alpha \right) :=  \left( Q_{3,1}(S) R_{j_{0}}  + S, \kappa \right)  $, and we stop the algorithm.  \\
%We then iterate the process $(A) + (B)$: the algorithm clearly ends after a finite of steps in view of (\ref{Eqn:Smallness}). Moreover $(A1)$, $(A2)$, $(A3)$, and $(A4)$ hold.
% \color{black}

\subsection{Proof of Lemma  \ref{lem:PerturbBourg}}
\label{Subsection:Perturb}

Let $F(v) := |v|^{2^{*} -2} v $ and $e := w_{1} - w_{2}$. It is made of four steps:

\begin{enumerate}

\item[$$]

\item $ \| w_{2} \|_{L_{t}^{\frac{2(d+1)}{d-1}} B^{\frac{1}{2}}_{\frac{2(d+1)}{d-1},2} (I) }
+ \left\| \vec{w}_{2} \right\|_{L_{t}^{\infty} \mathcal{H}(I)} \lesssim_{C_{0}} 1$. \\
\\
Indeed let $ 1 \gg \alpha > 0$. We can find a partition of $I$ into subintervals  $\left( J_{l}:=[a_{l},b_{l}] \right)_{1 \leq l \leq \bar{l}}$ such that
$\| w_{2} \|_{L_{t}^{\frac{2(d+1)}{d-2}} L_{x}^{\frac{2(d+1)}{d-2}} (J_{l})} = \alpha$ for all
$1 \leq l < \bar{l}$ , and $ \| w_{2} \|_{L_{t}^{\frac{2(d+1)}{d-2}} L_{x}^{\frac{2(d+1)}{d-2}} (J_{\bar{l}})} = \alpha $.
Let $ K_{1} := [a_{1},.] \subset J_{1} $ and $ Y(K_{1}) := \| w_{2} \|_{L_{t}^{\frac{2(d+1)}{d-1}} B^{\frac{1}{2}}_{\frac{2(d+1)}{d-1},2} (K_{1})}
+ \left\| \vec{w}_{2} \right\|_{L_{t}^{\infty} \mathcal{H} (K_{1})} $. By (\ref{Eqn:StrichKg0}) and
(\ref{Eqn:NonlinearEstGen}) we see that there
exists $C > 0$ such that

\begin{equation}
\begin{array}{ll}
Y(K_{1}) & \leq  C  ( C_{0} +  \alpha^{2^{*}-2} \| w_{2} \|_{L_{t}^{\frac{2(d+1)}{d-1}} B^{\frac{1}{2}}_{\frac{2(d+1)}{d-1},2} (K_{1}) } ) \\
& \leq C ( C_{0} +  \alpha^{2^{*}-2} Y(K_{1}) ) \cdot
\end{array}
\nonumber
\end{equation}
Hence a continuity argument shows that $ Y(J_{1})  \leq 2 C \, C_{0}$. Iterating over $l$ we get $ Y(J_{l})  \leq (2 C \, C_{0})^{l}$. Hence by summing over $l$ we get
$  \lesssim_{C_{0}} 1 $.

\item[$$]

\item Short-time perturbation argument.\\
\\
Let $ 1 \gg \alpha > 0$ and $J=[\bar{a}, \bar{b}] \subset I$ such that
$ \| w_{2} \|_{L_{t}^{\frac{2(d+1)}{d-2}} L_{x}^{\frac{2(d+1)}{d-2}}(J) }
+ \| w_{2} \|_{L_{t}^{\frac{2(d+1)}{d-1}} B^{\frac{1}{2}}_{\frac{2(d+1)}{d-1},2} (J) } \leq \alpha $. Then there exists $ 1 \gg \mu > 0 $ such that
if $\epsilon^{'} \leq \mu $ and for all $j \in \{ 1,2,3  \} $

\begin{equation}
\begin{array}{ll}
\left\|  \cos{ \left( (t-\bar{a}) \langle D \rangle \right)} e(\bar{a}) +
\frac{\sin{ \left( (t-\bar{a}) \langle D \rangle \right) }}{ \langle D \rangle} \partial_{t} e(\bar{a})
\right\|_{X_{j}(J) }  \leq \epsilon^{'}, \; \text{then} \\
 \\
(*) : \; \left\| F(w_{1}) - F(w_{2}) \right\|_{L_{t}^{\frac{2(d+1)}{d+3}} B^{\frac{1}{2}}_{\frac{2(d+1)}{d+3},2} (J) } \lesssim \epsilon^{'} \cdot
\end{array}
\nonumber
\end{equation}
Indeed let $ Z(J): = \| e \|_{L_{t}^{\frac{2(d+1)}{d-1}} B^{\frac{1}{2}}_{\frac{2(d+1)}{d-1},2} (J)}
+ \| e \|_{L_{t}^{\frac{2(d+1)}{d-2}} L_{x}^{\frac{2(d+1)}{d-2}}(J)} $. We get from (\ref{Eqn:DispFar1}) and (\ref{Eqn:StrichKg0}) followed by
(\ref{Eqn:NonlinearEstGenDiff})

\begin{equation}
\begin{array}{ll}
Z(J) & \lesssim \epsilon^{'} + \left\|  F(w_{1}) - F(w_{2}) \right\|_{L_{t}^{\frac{2(d+1)}{d+3}} B^{\frac{1}{2}}_{\frac{2(d+1)}{d+3},2} (J)} \\
 & \lesssim  \epsilon^{'} + \| e \|_{L_{t}^{\frac{2(d+1)}{d-1}} B^{\frac{1}{2}}_{ \frac{2(d+1)}{d-1},2}}
 \left( \| w_{1} \|^{2^{*}-2}_{L_{t}^{\frac{2(d+1)}{d-2}} L_{x}^{\frac{2(d+1)}{d-2}} (J) } +
 \| w_{2} \|^{2^{*}-2}_{L_{t}^{\frac{2(d+1)}{d-2}} L_{x}^{\frac{2(d+1)}{d-2}} (J) }   \right)  \\
 & +  \| e \|_{L_{t}^{\frac{2(d+1)}{d-2}} L_{x}^{\frac{2(d+1)}{d-2}} (J)}
 \left( \| w_{1} \|^{2^{*}-3}_{L_{t}^{\frac{2(d+1)}{d-2}} L_{x}^{\frac{2(d+1)}{d-2}} (J) } +
 \| w_{2} \|^{2^{*}-3}_{L_{t}^{\frac{2(d+1)}{d-2}} L_{x}^{\frac{2(d+1)}{d-2}} (J)}
 \right) \\
 &  \left(  \| w_{1} \|_{L_{t}^{\frac{2(d+1)}{d-1}} B^{\frac{1}{2}}_{ \frac{2(d+1)}{d-1},2}(J)}
 + \| w_{2} \|_{L_{t}^{\frac{2(d+1)}{d-1}} B^{\frac{1}{2}}_{ \frac{2(d+1)}{d-1},2}(J)}
 \right) \\
 & \lesssim  \epsilon^{'} +  Z^{2^{*}-1}(J) + \alpha Z^{2^{*}-2}(J)  + \alpha^{2^{*}-3} Z^{2}(J) + \alpha^{2^{*}-2} Z(J)
\end{array}
\nonumber
\end{equation}
Hence a continuity argument shows that $ Z(J) \lesssim \epsilon^{'}$. Hence $(*)$ holds.

\item[$$]

\item  Long-time perturbation argument. \\
\\
Observe from the first step and  (\ref{Eqn:Boundw2W}) that we can partition $I$ into subintervals $ \left( J_{j}:=[a_{j},b_{j}] \right)_{1 \leq j \leq \bar{j}}$ such that $\bar{j} \lesssim _{C_{0}} 1$, \\
$ \max \left( \| w_{2} \|_{L_{t}^{\frac{2(d+1)}{d-1}} B^{\frac{1}{2}}_{\frac{2(d+1)}{d-1},2} (J_{j})},
\| w_{2} \|_{L_{t}^{\frac{2(d+1)}{d-2}} L_{x}^{\frac{2(d+1)}{d-2}} (J_{j})} \right) = \alpha $  for
$1 \leq j \leq \bar{j}-1 $ and $ \max \left( \| w_{2} \|_{L_{t}^{\frac{2(d+1)}{d-1}} B^{\frac{1}{2}}_{\frac{2(d+1)}{d-1},2} (J_{\bar{j}})} ,
\| w_{2} \|_{L_{t}^{\frac{2(d+1)}{d-2}} L_{x}^{\frac{2(d+1)}{d-2}} (J_{\bar{j}})} \right) \leq \alpha $. We claim that for all $1 \leq j \leq \bar{j}$ there
there exists a positive constant $C(j)$ such that (with $F(x):= |x|^{2^{*}-2} x$)

\begin{equation}
\begin{array}{l}
\left\| F(w_{2}) - F(w_{1}) \right\|_{L_{t}^{\frac{2(d+1)}{d+3}} B^{\frac{1}{2}}_{\frac{2(d+1)}{d+3},2} (J_{j})} \leq C(j) \epsilon_{0} \cdot
\end{array}
\label{Eqn:EstDiffF}
\end{equation}
Indeed by induction assume that (\ref{Eqn:EstDiffF}) holds for all $1 \leq j \leq k$ with $ k \in \{ 1,..., \bar{j} -1 \}$. Then from the Duhamel formula above
(\ref{Eqn:StrichKg0}) we see that

\begin{equation}
\begin{array}{l}
\cos{ \left( (t- a_{k+1}) \langle D \rangle \right)} e( a_{k+1} ) +  \frac{\sin{ \left( (t- a_{k+1})  \langle D \rangle \right) }}{ \langle D \rangle}
\partial_{t} e(a_{k+1})  \\ = \cos{ \left( (t- a) \langle D \rangle \right)} e(a) -
\frac{\sin{ \left( (t- a)  \langle D \rangle \right) }}{ \langle D \rangle} \partial_{t} e(a) - \sum \limits_{j=1}^{k} \int_{a_{j}}^{b_{j}} \frac{\sin \left( (t-t^{'}) \langle D \rangle \right) }{\langle D \rangle} \left( F(w_{1}(t')) - F(w_{2}(t')) \right) \; dt'
\end{array}
\nonumber
\end{equation}
Hence we get from (\ref{Eqn:StrichKg0})

\begin{equation}
\begin{array}{l}
\left\|  \cos{ \left( (t- a_{k+1}) \langle D \rangle \right)} e( a_{k+1} )+ \frac{\sin{ \left( (t- a_{k+1})  \langle D \rangle \right) }}{ \langle D \rangle} \partial_{t} e(a_{k+1}) \right\|_{X_{j}(J_{k+1}) } \\
\lesssim \epsilon_{0} +  \sum \limits_{j=1}^{k} \left\| F(w_{2}) - F(w_{1}) \right\|_{L_{t}^{\frac{2(d+1)}{d+3}} B^{\frac{1}{2}}_{\frac{2(d+1)}{d+3},2} (J_{j})} \lesssim \epsilon_{0} + \sum \limits_{j=1}^{k} C(j) \epsilon_{0} \leq \mu
\end{array}
\nonumber
\end{equation}
Hence an application of the short-time perturbation argument shows that (\ref{Eqn:EstDiffF}) also holds for $j= k+1$. This proves the claim. \\
Hence by summation

\begin{equation}
\begin{array}{l}
\left\| F(w_{1}) - F(w_{2}) \right\|_{L_{t}^{\frac{2(d+1)}{d+3}} B^{\frac{1}{2}}_{\frac{2(d+1)}{d+3},2} (I)} \lesssim
\sum \limits_{j=1}^{\bar{j}} C(j) \epsilon_{0} \cdot
\end{array}
\label{Eqn:EstDiffI}
\end{equation}

\item[$$]

\item Estimate (\ref{Eqn:Diffw1w2}): this follows from (\ref{Eqn:DispFar1}), the embeddings $B^{0}_{2^{*},2} \hookrightarrow L^{2^{*}}$ and  $ \dot{H}^{1} \hookrightarrow L^{2^{*}} $, and (\ref{Eqn:EstDiffI}).

\end{enumerate}

\section{Proof of Proposition \ref{Prop:FarFromGdPr}}
\label{Section:FarFromGdPr}

In this section we prove Proposition \ref{Prop:FarFromGdPr}. \\

By using time translation symmetry of (\ref{Eqn:NlkgCrit}) if necessary we may assume WLOG that $t_{0}=0$. The continuity of the flow, Proposition \ref{Prop:VarEst} and Proposition \ref{Prop:DistFuncEigen} show that $\Theta \left( u ( [0, T_{+}(u)) ) \right) = \pm 1 $ if $ \Theta \left( u(0) \right) = \pm 1 $.

\subsection{$\mathbf{\Theta (u(t_{0})) = -1}$}

We prove that $T_{+}(u) < \infty$. We adapt an argument of Payne-Sattinger
\cite{paynesatt} to energies slightly above that of the ground states. Let $y := \langle u , u  \rangle$. Elementary computations show that

\begin{eqnarray*}
\partial_{t} y & =  & 2 \langle \partial_{t} u, u \rangle
\end{eqnarray*}
and

\begin{equation}
\begin{array}{lll}
\partial^{2}_{t} y & = & 2 \| \partial_{t} u \|^{2}_{L^{2}} + 2 \langle \partial_{tt} u,u \rangle \\
& = & 2 \| \partial_{t} u \|^{2}_{L^{2}} - 2 K(u) - 2 \| u \|^{2}_{L^{2}}
\end{array}
\label{Eqn:DerivSec0}
\end{equation}
Since $\tilde{d}_{\mathcal{S}} (\vec{u}(t)) \geq \delta_{b}$, $t \geq t_{2}$, we get from Proposition \ref{Prop:VarEst}
(with $k:=k(\delta_{b})$)

\begin{equation}
\begin{array}{lll}
\partial^{2}_{t} y & \geq & 2 \| \partial_{t} u \|^{2}_{L^{2}} + 2 k - 2 \| u \|^{2}_{L^{2}}
\end{array}
\label{Eqn:DerivSec1}
\end{equation}
We have

\begin{eqnarray*}
K(\vec{u}) + \| u \|^{2}_{L^{2}} & = & 2^{*} E(\vec{u}) - \frac{2^{*}}{2}  \| \partial_{t} u \|^{2}_{L^{2}} - 2^{*} \left( \frac{1}{2} - \frac{1}{2^{*}}  \right)
( \| \nabla u \|^{2}_{L^{2}} + \| u \|^{2}_{L^{2}} )
\nonumber
\end{eqnarray*}
From $\vec{u}(0) \in \mathcal{H}$, the conservation of energy,  and (\ref{Eqn:Charac3GdState}) we get

\begin{eqnarray*}
E(\vec{u}) & \leq & \left( \frac{1}{2} - \frac{1}{2^{*}} \right) \| \nabla u \|^{2}_{L^{2}} + \bar{\epsilon}^{2}
\label{Eqn:BoundE}
\end{eqnarray*}
Hence collecting these estimates and (\ref{Eqn:DerivSec0})

\begin{equation}
\begin{array} {lll}
\partial^{2}_{t} y & \geq & 2 \left( 1 + \frac{2}{2^{*}}  \right) \| \partial_{t} u \|^{2}_{L^{2}} - 2 \times 2^{*} \bar{\epsilon}^{2}
+ 2 \times 2^{*} \left( \frac{1}{2} - \frac{1}{2^{*}} \right) \| u \|^{2}_{L^{2}}
\end{array}
\label{Eqn:DerivSec2}
\end{equation}
From (\ref{Eqn:DerivSec1}) and (\ref{Eqn:DerivSec2})  we see that
$\partial^{2}_{t} y \geq \bar{\epsilon}^{2}$ . Hence, assuming that $T_{*} = \infty$ this leads to

\begin{equation}
\begin{array}{lll}
y(t) & \rightarrow & \infty, \, \text{as} \, t \rightarrow \infty
\end{array}
\label{Eqn:Behy}
\end{equation}
Since $|\partial_{t} y| \leq 2 \| \partial_{t} u \|_{L^{2}} \| u \|_{L^{2}}$ we get for $t \gg 1$

\begin{eqnarray*}
\partial_{tt} y & \geq & 2 \left( 1 +  \frac{2^{*}}{2} \right) \frac{ (\partial_{t} y)^{2}}{4 y}
\nonumber
\end{eqnarray*}
Hence $\partial^{2}_{t} y^{-\alpha} \leq 0 $ with $\alpha:= \frac{1}{2} \left( \frac{2^{*}}{2} -1 \right)$.
This contradicts (\ref{Eqn:Behy}).

\subsection{$ \mathbf{\Theta (u(t_{0})) = 1} $}

Throughout this section, we work with complex-valued functions so that we can use the results in \cite{ibramasmnak}. Unless previously defined we denote by
$\vec{f} = f_{1} + i f_{2} $ a complex-valued function. Let $E (\vec{f}) := \frac{1}{2} \| \vec{f} \|^{2}_{L^{2}} -
\frac{1}{2^{*}} \left\| Re ( \langle D \rangle^{-1} \vec{f} )  \right\|^{2^{*}}_{L^{2^{*}}}  $ and
$ E_{wa}( \vec{f}) := \frac{1}{2} \| \vec{f} \|^{2}_{L^{2}} - \frac{1}{2^{*}} \left\| \Re ( D^{-1} \vec{f} )  \right\|^{2^{*}}_{L^{2^{*}}} $. If
$\vec{f} := q f_{1}  + i f_{2}  $ with $q \in \{ D, \langle D \rangle \}$ and $f_{1}$, $f_{2}$ real-valued functions then
$ \tilde{d}_{\mathcal{S}}(\vec{f}) := \tilde{d}_{\mathcal{S}} ((f_{1}, f_{2})) $. \\

We recall the following result proved in \cite{ibramasmnak}  and obtained by a concentration compactness
procedure (see e.g \cite{bahger,merlevega}, relying on \cite{lions})

\begin{prop} \cite{ibramasmnak}
Let $ \left( v_{n}(0), \partial_{t} v_{n}(0) \right)$ be a bounded sequence in $\mathcal{H}$. Let
$ \vec{v}_{n}(0)  := \langle D \rangle v_{n}(0) - i \partial_{t} v_{n}(0)$. Let $ \left\{ \vec{v}_{n}(t) := e^{i \langle D \rangle t } \vec{v}_{n}(0) \right\}_{n \in \mathbb{N}}$ be a sequence that is bounded in $L^{2}$. Then there exist $ \left\{ \vec{\phi}^{j} \in L^{2}(\mathbb{R}^{d}) \right\}_{j \in \mathbb{N}}$, $ \left\{ \vec{w}_{n}^{j} \in 
L^{2}(\mathbb{R}^{d}) \right\}_{j \in \mathbb{N}}$ and $ \left\{ (t_{n}^{j}, x_{n}^{j}, h_{n}^{j} ) \subset \mathbb{R} \times \mathbb{R}^{d} \times (0,1] \right\}_{n \in \mathbb{N}, j \in \mathbb{N} }$
such that, up to a subsequence, $h_{n}^{j} \rightarrow 0$ as $n \rightarrow \infty$ or $h_{n}^{j} =1$ for all $n$, and

\begin{equation}
\begin{array}{ll}
k \in \mathbb{N}: \vec{v}_{n}(t) = \sum_{j=0}^{k} \vec{v}_{n}^{j}(t) + e^{i  \langle D \rangle t} \vec{w}_{n}^{k}, \, \vec{v}_{n}^{j} := e^{i \langle D \rangle (t -t_{n}^{j})} T_{n}^{j} \vec{\phi}^{j},
\end{array}
\label{Eqn:Decompvn}
\end{equation}
with

\begin{equation}
\begin{array}{ll}
\lim \limits_{k \rightarrow \infty} \overline{\lim}_{n \rightarrow \infty} \| e^{i \langle D \rangle t} \vec{w}^{k}_{n} \|_{ L_{t}^{\infty} B^{-\frac{d}{2}}_{\infty,\infty}} & =0,
\end{array}
\label{Eqn:Asympwnk1}
\end{equation}

\begin{equation}
\begin{array}{l}
j \neq l: \; \lim \limits_{n \rightarrow \infty} \langle \vec{v}_{n}^{j}(t), \vec{v}_{n}^{l}(t) \rangle = 0, \\ \\
\lim \limits_{n \rightarrow \infty} \langle \vec{v}_{n}^{j}(t),  e^{i \langle D \rangle t} \vec{w}_{n}^{k} \rangle = 0, \; \text{and}
\end{array}
\label{Eqn:Asympwnk2}
\end{equation}

\begin{equation}
\begin{array}{l}
j \neq l:\ ; \lim \limits_{n \rightarrow \infty} \left| \log{ \left( \frac{h_{n}^{l}}{h_{n}^{j}} \right) } \right|
+ \frac{| t_{n}^{l} - t_{n}^{j}| + |x_{n}^{l} - x_{n}^{j}| }{h_{n}^{l}} = \infty \cdot
\end{array}
\label{Eqn:OrthParam}
\end{equation}
In particular

\begin{equation}
\begin{array}{ll}
\lim \limits_{n \rightarrow \infty}  \| \vec{v}_{n}(t) \|^{2}_{L^{2}} - \sum_{j=0}^{k} \| \vec{v}_{n}^{j}(t) \|^{2}_{L^{2}} - \| e^{i \langle D \rangle t} \vec{w}_{n}^{k} \|^{2}_{L^{2}} = 0
\end{array}
\label{Eqn:DecouplKinetNrj}
\end{equation}

\label{Proposition:ProfileDecomp}
\end{prop}

We define the following statement $\mathcal{P}(E)$: if $\vec{u}$ is a solution of (\ref{Eqn:NlkgCrit}) with data
$(u_{0},u_{1}) \in \mathcal{H}$, forward maximal time of existence $T_{+}(u)$, and energy $E$ that satisfies

\begin{enumerate}
\item  $ \tilde{d}_{\mathcal{S}} \left( \vec{u}(t) \right) \geq \delta_{b} $, $t \in [0, T_{+}(u))$
\item  $ \Theta (u (t)) = 1 $, $t \in [ 0, T_{+}(u))$
\item $  E(\vec{u}) \leq E_{wa}(\vec{W}) + \bar{\epsilon}^{2} $
\end{enumerate}
Then $T_{+}(u)= \infty$ and $\| u \|_{L_{t}^{\frac{2(d+1)}{d-2}} L_{x}^{\frac{2(d+1)}{d-2}} (\mathbb{R}^{+})} < \infty$ \footnote{This implies scattering, by a standard
procedure}. We define

\begin{equation}
\begin{array}{ll}
E_{c} & =  \sup \left\{ E > 0, \mathcal{P}(E) \, holds \right\}
\end{array}
\nonumber
\end{equation}
Clearly from \cite{ibramasmnak}, we know that we may assume that $E_{c} \geq E_{wa}(\vec{W})$. Assume toward a contradiction that
$E_{c} = E_{wa}(\vec{W}) + \tilde{\epsilon}^{2}$ with $\tilde{\epsilon} < \bar{\epsilon}$. We prove the following claim:

\begin{claim}
There exists a critical element $ U_{c} $, i.e a solution of (\ref{Eqn:NlkgCrit}), that satisfies $ E(\vec{U}_{c}) = E_{c} $, $ (1) $ and
 $(2)$. Here $\vec{U}_{c} := \langle D \rangle U_{c} - i \partial_{t} U_{c}$. Moreover $U_{c}$ does not scatter, i.e $ \| U_{c} \|_{ L_{t}^{\frac{2(d+1)}{d-2}} L_{x}^{\frac{2(d+1)}{d-2}} \left( [0,T_{+}(U_{c})) \right)} = \infty$.
\label{Claim:CritElem}
\end{claim}

By definition of $E_{c}$ there exists a sequence $(\vec{u}_{n})_{n \geq 1}$ of solutions $\vec{u}_{n}:= \langle D \rangle u_{n} - i \partial_{t} u_{n}$ with forward
maximal time of existence $T_{+}(u_{n})$ that satisfy $(1)$, $(2)$ and $(3)$, $E(\vec{u}_{n}) \rightarrow E_{c}$ from above
and $ \| u_{n} \|_{L_{t}^{\frac{2(d+1)}{d-2}} L_{x}^{\frac{2(d+1)}{d-2}} ([0,T_{+}(u_{n})])} = \infty$. \\
We may assume WLOG that for all $t \in [0, T_{+}(u_{n}))$, $\tilde{d}_{\mathcal{S}} \left( \vec{u}_{n}(t) \right) \geq \delta_{b}^{'} $ with $\delta_{b}^{'} \gg \delta_{b}$ .
If not we may find $\bar{t}_{n} \in I(u)$ such that $ \delta_{b} \leq \tilde{d}_{\mathcal{S}} \left( \vec{u}_{n} (\bar{t}_{n}) \right) < \delta_{b}^{'}$. Hence, if
$\partial_{t} \tilde{d}^{2}_{\mathcal{S}} (\vec{u}_{n})(\bar{t}_{n}) \geq 0 $ (resp. $ \partial_{t}  \tilde{d}^{2}_{\mathcal{S}} (\vec{u})(\bar{t}_{n}) < 0 $)  then we can apply Proposition \ref{Prop:DynEjection} to find $\bar{t}_{n}^{'} > \bar{t}_{n}$ (resp. $\bar{t}_{n}^{'} < \bar{t}_{n}$ ) such that $ \tilde{d}_{\mathcal{S}} \left( u (\bar{t}_{n}^{'}) \right) = \delta_{b}^{'} $  and we replace $\vec{u}_{n}(t)$ with $\vec{u}(t + \bar{t}_{n}^{'})$. \\
We claim that there exists $t_{n} \in [0, T_{+}(u_{n}))$ such that $\| \nabla u_{n}(t_{n}) \|_{L^{2}} \gtrsim 1$. If not $(*)$ would hold, with
$ (*): \, \| \nabla u_{n}(t) \|_{L^{2}} \ll 1, \, t \in [0, T_{+}(u_{n})) $. Now  defining
$ X_{n}(t) :=  \| u_{n} \|_{L_{t}^{\frac{2(d+1)}{d-2}} L_{x}^{\frac{2(d+1)}{d-2}} ([ 0,T_{+}(u_{n}) )) } +  \| u_{n} \|_{ L_{t}^{\frac{2(d+1)}{d-1}} B^{\frac{1}{2}}_{\frac{2(d+1)}{d-1},2} ([ 0,T_{+}(u_{n}) ))} $ and $r$ such that $\frac{1}{r} = \frac{d-1}{2(d+1)} - \frac{1}{2d}$  we get for some $0 < c < 1$

\begin{equation}
\begin{array}{ll}
X_{n}(t) \lesssim \| \vec{u}_{n}(0) \|_{\mathcal{H}} + \| u_{n} \|^{2^{*} -2}_{L_{t}^{\frac{2(d+1)}{d-2}} L_{x}^{\frac{2(d+1)}{d-2}} ([0,T_{+}(u_{n})]) }
\| u_{n} \|_{ L_{t}^{\frac{2(d+1)}{d-1}} B^{\frac{1}{2}}_{\frac{2(d+1)}{d-1},2} ([ 0,T_{+}(u_{n}) ))}   \\
\lesssim \| \vec{u}_{n}(0) \|_{\mathcal{H}} + \| u_{n} \|^{(2^{*} -2)c}_{L_{t}^{\infty} L_{x}^{2^{*}} ([0,T_{+}(u_{n}))}
\| u_{n} \|^{(2^{*}-2)(1-c)}_{L_{t}^{\frac{2(d+1)}{d-1}} L_{x}^{r}} \| u_{n} \|_{ L_{t}^{\frac{2(d+1)}{d-1}} B^{\frac{1}{2}}_{\frac{2(d+1)}{d-1},2} ([ 0,T_{+}(u_{n}) ))}  \\
\lesssim  \| \vec{u}_{n}(0) \|_{\mathcal{H}} + o \left( X^{ 1 + (2^{*}-2)(1-c)}_{n}(t) \right),
\end{array}
\nonumber
\end{equation}
where at the third line we used (\ref{Eqn:SobIneq}), $(*)$, and $H^{\frac{2(d+1)}{d-1}} \hookrightarrow L^{r}$, followed by a continuity argument at the fourth line. This is a contradiction. Hence, by letting $\vec{u}_{n}(t):= \vec{u}_{n}(t+ t_{n})$, we may assume that

\begin{equation}
\begin{array}{ll}
\| \nabla u_{n}(0) \|_{L^{2}} \gtrsim 1
\end{array}
\label{Eqn:LowerBdKinet}
\end{equation}
Unless otherwise specified, let $0 < \alpha \ll 1$, $k \gg 1$, and $n \gg 1$ such that all the statements below are true. Applying
Proposition \ref{Proposition:ProfileDecomp} to $\vec{v}_{n}(0) := \vec{u}_{n}(0):=  \langle D \rangle u_{n}(0) - i \partial_{t} u_{n}(0)$, we get

\begin{equation}
\begin{array}{ll}
\vec{u}_{n}(0) & = \sum_{j=0}^{k} e^{ -i \langle D \rangle t_{n}^{j}} T_{n}^{j} \vec{\phi}^{j} + \vec{w}_{n}^{k} \cdot
\end{array}
\nonumber
\end{equation}
We define the nonlinear concentrating wave associated with $ \left\{ e^{-i \langle D \rangle t_{n}^{j}} T_{n}^{j} \vec{\phi}^{j} \right\} $ in the following fashion

\begin{equation}
\begin{array}{ll}
\vec{u}^{j}_{n} & : = T_{n}^{j} \vec{U}_{n}^{j}  \left( \frac{t - t_{n}^{j}}{h_{n}^{j}} \right)
\end{array}
\nonumber
\end{equation}
with

\begin{equation}
\begin{array}{ll}
\vec{U}_{n}^{j} & = e^{ i t \langle D \rangle^{j}_{n}} \vec{\phi}^{j} - i \int_{\tau_{\infty}^{j}}^{t} e^{i(t-s) \langle D \rangle^{j}_{n}} F
\left( U_{n}^{j}(s) \right) \, ds,
\end{array}
\nonumber
\end{equation}
$F(x) := |x|^{\frac{4}{d-2}} x$, $\tau_{\infty}^{j} := \lim_{n \rightarrow \infty} \tau_{n}^{j} := \lim_{n \rightarrow \infty} -\frac{t_{n}^{j}}{h_{n}^{j}}$
(up to a subsequence) and
$U_{n}^{j} := \Re \left( (\langle D \rangle_{n}^{j})^{-1}  \vec{U}_{n}^{j} \right)$. Observe that
$\vec{U}_{n}^{j}$ satisfies $ (i \partial_{t}  + \langle D \rangle_{n}^{j}) \vec{U}_{n}^{j} = F(U_{n}^{j})$ and, consequently

\begin{equation}
\begin{array}{l}
\vec{U}_{n}^{j} = \langle D \rangle_{n}^{j} U_{n}^{j} - i \partial_{t} U_{n}^{j} \\
\left( \partial_{tt} - \triangle + (h_{n}^{j})^{2} \right) U^{j}_{n} = F \left(  U_{n}^{j}  \right) \\
\end{array}
\nonumber
\end{equation}
Let $u_{n}^{j}:= \Re \left( \langle D \rangle^{-1} \vec{u}_{n}^{j} \right)$. Then observe that
$ u_{n}^{j} =  h_{n}^{j} T_{n}^{j} U_{n}^{j} \left( \frac{t - t_{n}^{j}}{h_{n}^{j}} \right)$ and that $u_{n}^{j}$
is a solution of (\ref{Eqn:NlkgCrit}). We will use $u_{n}^{j}$ when we apply Result \ref{res:PertubationResult}.\\
\\
We now turn to the local existence of $U_{n}^{j}$ around $\tau_{\infty}^{j}$. We claim the following: \\
\\
\underline{Claim}: Let $ 0 < \delta \ll 1$. Then there exists a small open interval $J$ containing $\tau_{\infty}^{j}$
and that does not depend on $n$ such that

\begin{equation}
\begin{array}{ll}
\left\| \cos{( t \langle D \rangle_{n}^{j})} \Re \left( (\langle D \rangle_{n}^{j})^{-1} \vec{\phi}^{j} \right) -
\frac{\sin{ ( t \langle D \rangle_{n}^{j} )}}{\langle D \rangle_{n}^{j}}  \Im  ( \vec{\phi}^{j} )
\right\|_{L_{t}^{\frac{2(d+1)}{d-2}} L_{x}^{\frac{2(d+1)}{d-2}} (J)} & \leq \delta
\end{array}
\label{Eqn:Appdelt}
\end{equation}

\begin{proof}

Assume now that $h_{\infty}^{j}=1$. Then (\ref{Eqn:Appdelt}) clearly follows from the dominated convergence theorem and
(\ref{Eqn:StrichKg0}). \\
\\
Assume now that $h_{\infty}^{j} = 0$. \\
First assume that $\tau_{\infty}^{j} \in \mathbb{R}$. Let $\vec{\phi}^{j}_{\delta}$ be a Schwartz function such that $ \supp \left( \widehat{\vec{\phi}^{j}_{\delta}} \right)  \subset \mathbb{R}^{d} - \{ (0,..,0)  \} $
and $ \| \vec{\phi}^{j}_{\delta} - \vec{\phi}^{j} \| \ll \delta $. Let
$\left( q(x), X_{n,j} (\vec{f}), X_{j}(\vec{f}) \right) :=
\left(  \cos(x), \Re \left( (\langle D \rangle_{n}^{j})^{-1} \vec{f} \right), \Re \left( D^{-1} \vec{f} \right) \right) $ or
$ \left( q(x), X_{n,j}(\vec{f}), X_{j} (\vec{f})   \right) := \left( \frac{\sin(x)}{x}, \Im (\vec{f}), \Im (\vec{f}) \right)$. Then integration by parts (using in particular the formula $e^{i \xi \cdot x} = \frac{x \cdot \nabla (e^{i \xi \cdot x})}{|x|^{2}}$ ) and the Schwartz nature of
$\vec{\phi}^{j}_{\delta}$  allow to use the dominated convergence theorem to conclude that \\
$  \left\| q(t \langle D \rangle_{n}^{j})
X_{n}^{j} (\vec{\phi}^{j}_{\delta}) - q(t D) X_{j} (\vec{\phi}^{j}_{\delta}) \right\|_{L_{t}^{\frac{2(d+1)}{d-2}} L_{x}^{\frac{2(d+1)}{d-2}} (J)} \ll \delta $. We also get from (\ref{Eqn:StrichWave0}) and (\ref{Eqn:StrichWave0Sc}) that
$ \left\| q(t \langle D \rangle_{n}^{j} ) X_{n,j} ( \vec{\phi}^{j} - \vec{\phi}^{j}_{\delta})  \right\|_{L_{t}^{\frac{2(d+1)}{d-2}} L_{x}^{\frac{2(d+1)}{d-2}} (\mathbb{R})} \lesssim \| \vec{\phi}^{j}_{\delta} - \vec{\phi}^{j} \|_{L^{2}} $, and that
$ \left\| q(tD) X_{j} \left( \vec{\phi}^{j}_{\delta}) \right)  \right\|_{L_{t}^{\frac{2(d+1)}{d-2}} L_{x}^{\frac{2(d+1)}{d-2}} (J)} \ll \delta
$. \\
Now assume that $\tau_{\infty}^{j} \in \{ -\infty, \infty \}$. We only prove the claim for $\tau_{\infty}^{j} = \infty $: the proof for $\tau_{\infty}^{j} = -\infty $
is a straightforward modification of the proof for $\tau_{\infty}^{j} = \infty $. It suffices to prove that for $S \gg 1$,
$ \left\| \frac{e^{\pm i t \langle D \rangle_{j}^{n} }}{\langle  D \rangle_{n}^{j}} \phi  \right\|_{L_{t}^{\frac{2(d+1)}{d-2}} L_{x}^{\frac{2(d+1)}{d-2}}  ((S,\infty))} \leq \delta $. Following a similar scheme as the one from (\ref{Eqn:EstDispnPhi}) to (\ref{Eqn:ListEstUnifn}), we get the former estimate, taking into
account that $  \left\|  e^{\pm i t \langle D \rangle } \phi  \right\|_{B^{0}_{\frac{2(d+1)}{d-2},2}} \lesssim  \frac{1}{|t|^{\frac{3(d-1)}{2(d+1)}}}
\| \phi \|_{B^{\frac{3}{2}}_{\left( \frac{2(d+1)}{d-2} \right)^{'}, 2}} $.
\end{proof}
Then, by (\ref{Eqn:StrichWave0Sc}), (\ref{Eqn:NonlinearEstGen}), and (\ref{Eqn:NonlinearEstGenDiff}), we get

\begin{equation}
\begin{array}{l}
\begin{array}{ll}
\| w  \|_{ L_{t}^{\frac{2(d+1)}{d-2}} L_{x}^{\frac{2(d+1)}{d-2}} (J)} & \lesssim
\left\| \cos{( t \langle D \rangle^{j}_{n})} \Re \left( (\langle D \rangle^{j}_{n})^{-1} \vec{\phi}^{j} \right) -
\frac{\sin{ ( t \langle D \rangle^{j}_{n} )}}{\langle D \rangle^{j}_{n}}  \Im  ( \vec{\phi}^{j} )
\right\|_{ L_{t}^{\frac{2(d+1)}{d-2}}  L_{x}^{\frac{2(d+1)}{d-2}} (J)}  \\
 & + \| w \|^{2^{*} -2}_{L_{t}^{\frac{2(d+1)}{d-2}} L_{x}^{\frac{2(d+1)}{d-2}} (J)}
 \| w \|_{ L_{t}^{\frac{2(d+1)}{d-2}} \bar{B}^{\frac{1}{2},h_{n}^{j}}_{\frac{2(d+1)}{d-1},2}  (J) },
\end{array} \\
\\
\begin{array}{l}
\| w \|_{ L_{t}^{\frac{2(d+1)}{d-2}} \bar{B}^{\frac{1}{2},h_{n}^{j}}_{\frac{2(d+1)}{d-1},2}  (J)} +
\| (w,\partial_{t} w) \|_{L_{t}^{\infty} \mathcal{H}_{h_{n}^{j}} ([\tau_{n}^{j}, \tau_{\infty}^{j}])}  \\
\lesssim  \| \vec{\phi}^{j} \|_{L^{2}} + \| w \|^{2^{*}-2}_{L_{t}^{\frac{2(d+1)}{d-2}} L_{x}^{\frac{2(d+1)}{d-2}}(J)  }
\| w \|_{ L_{t}^{\frac{2(d+1)}{d-2}} \bar{B}^{\frac{1}{2},h_{n}^{j}}_{\frac{2(d+1)}{d-1},2}(J) }, \text{and}
\end{array}
\end{array}
\label{Eqn:NonLinProf12}
\end{equation}

%\begin{equation}
%\begin{array}{l}
%\| w \|_{ L_{t}^{\frac{2(d+1)}{d-2}} \bar{B}^{\frac{1}{2},h_{n}^{j}}_{\frac{2(d+1)}{d-1},2}  (J)} +
%\| w \|_{L_{t}^{\infty} L_{x}^{2} ([\tau_{n}^{j}, \tau_{\infty}^{j}])}  \\
%\lesssim  \| \vec{\phi}^{j} \|_{L^{2}} + \| w \|^{2^{*}-2}_{L_{t}^{\frac{2(d+1)}{d-2}} L_{x}^{\frac{2(d+1)}{d-2}}(J)  }
%\| w \|_{ L_{t}^{\frac{2(d+1)}{d-2}} \bar{B}^{\frac{1}{2},h_{n}^{j}}_{\frac{2(d+1)}{d-1},2}(J) },
%\end{array}
%\label{Eqn:NonLinProf2}
%\end{equation}

\begin{equation}
\begin{array}{l}
\| w - \tilde{w} \|_{ L_{t}^{\frac{2(d+1)}{d-2}} \bar{B}^{\frac{1}{2},h_{n}^{j}}_{\frac{2(d+1)}{d-1},2}  (J)} \\
\lesssim
\| w - \tilde{w} \|_{ L_{t}^{\frac{2(d+1)}{d-1}}  \bar{B}^{\frac{1}{2},h_{n}^{j}}_{\frac{2(d+1)}{d-1},2} (J)}
\left(
\begin{array}{l}
\| w \|^{2^{*} - 2}_{L_{t}^{\frac{2(d+1)}{d-2}} L_{x}^{\frac{2(d+1)}{d-2}} (J)}  +
\| \tilde{w} \|^{2^{*} - 2}_{L_{t}^{\frac{2(d+1)}{d-2}} L_{x}^{\frac{2(d+1)}{d-2}} (J)}
\end{array}
\right) \\
+ \| w - \tilde{w} \|_{L_{t}^{\frac{2(d+1)}{d-2}} L_{x}^{\frac{2(d+1)}{d-2}} (J)}
\left(
\| w \|^{2^{*}-3}_{L_{t}^{\frac{2(d+1)}{d-2}} L_{x}^{\frac{2(d+1)}{d-2}} (J)}
+ \| \tilde{w} \|^{2^{*}-3}_{L_{t}^{\frac{2(d+1)}{d-2}} L_{x}^{\frac{2(d+1)}{d-2}} (J)}
\right)  \\
\left(
\| w \|_{L_{t}^{\frac{2(d+1)}{d-1}} \bar{B}^{\frac{1}{2},h_{n}^{j}}_{\frac{2(d+1)}{d-1}} (J) } +
\| \tilde{w} \|_{L_{t}^{\frac{2(d+1)}{d-1}} \bar{B}^{\frac{1}{2},h_{n}^{j}}_{\frac{2(d+1)}{d-1}} (J) }
\right)
\end{array}
\label{Eqn:NonLinProf3}
\end{equation}
From the above claim, (\ref{Eqn:NonLinProf12}), and (\ref{Eqn:NonLinProf3}), we see (see e.g., \cite{kenmer}) that we can construct a
solution on $J$ by using a standard fixed point argument . \\
Let $h_{\infty}^{j} := \lim \limits_{n \rightarrow \infty} h_{n}^{j}$.

\begin{rem}
Assume that $h_{\infty}^{j} =1$. If $ \tau_{\infty}^{j} \neq \{ -\infty, \infty \}$, then

\begin{equation}
\begin{array}{ll}
\left\| \vec{U}_{n}^{j} (t) - e^{i \tau_{\infty}^{j} \langle D \rangle }  \vec{\phi}^{j} \right\|_{L^{2}}
& \lesssim \left\| \vec{U}_{n}^{j} (t) - e^{i t  \langle D \rangle  } \vec{\phi}^{j} \right\|_{L^{2}}
+ \left\| ( e^{i t \langle D \rangle } -  e^{i \tau_{\infty}^{j} \langle D \rangle} ) \vec{\phi}^{j} \right\|_{L^{2}} \\
& \rightarrow_{t \rightarrow \tau_{\infty}^{j}} 0,
\end{array}
\label{Eqn:AsympBehUn}
\end{equation}
since, by construction of $\vec{U}_{n}^{j}$, \\
$ \left\| \vec{U}_{n}^{j} (t) - e^{i t \langle D \rangle} \vec{\phi}^{j} \right\|_{L^{2}}
 \lesssim \| U_{n}^{j} \|^{2^{*} -2}_{ L_{t}^{\frac{2(d+1)}{d-2}} L_{x}^{\frac{2(d+1)}{d-2}} ( [\tau_{n}^{j}, \tau_{\infty}^{j}] )}  \| U_{n}^{j} \|_{  B^{\frac{1}{2}}_{\frac{2(d+1)}{d-1},2} \left([\tau_{n}^{j}, \tau_{\infty}^{j}] \right)},
$
which implies that $ \lim \limits_{t \rightarrow \tau_{\infty}^{j}} \left\| \vec{U}_{n}^{j} (t) - e^{i t \langle D \rangle} \vec{\phi}^{j}
\right\|_{L_{x}^{2}} = 0 $. If $ \tau_{\infty}^{j} \in \{ -\infty,  \infty \}$ then

\begin{equation}
\begin{array}{ll}
 \left\| \vec{U}_{n}^{j} (t) - e^{i t \langle D \rangle} \vec{\phi}^{j} \right\|_{L^{2}} & \rightarrow _{t \rightarrow \tau_{\infty}^{j} } 0
\end{array}
\label{Eqn:UnScatt}
\end{equation}

\label{Rem:UnjProp}

\end{rem}

\begin{rem}

Notice that the existence of  $ \vec{U}_{\infty}^{j} $ can be proved using a similar scheme to prove the existence of
$\vec{U}_{n}^{j}$, replacing (\ref{Eqn:StrichWave0Sc})  with (\ref{Eqn:StrichWave0}). This allows to define $I_{max}(\vec{U}_{\infty}^{j})$. By construction we get the following estimates

\begin{equation}
\begin{array}{ll}
\tau _{\infty}^{j} \neq \{ -\infty, \infty \}:  \left\| \vec{U}_{\infty}^{j} (t) - e^{i \tau_{\infty}^{j} D }  \vec{\phi}^{j} \right\|_{L^{2}}  \rightarrow_{t \rightarrow \tau_{\infty}^{j}} 0 \\
\tau_{\infty}^{j} \in \{ -\infty, \infty \}: \left\| \vec{U}_{\infty}^{j} (t) - e^{i t D }  \vec{\phi}^{j} \right\|_{L^{2}}  \rightarrow_{t \rightarrow \infty} 0
\end{array}
\nonumber
\end{equation}
\label{Rem:UinfjProp}

\end{rem}

Next we write the following theorem that allows to approximate well $\vec{U}_{n}^{j}$ with $\vec{U}_{\infty}^{j}$ in the most difficult case, i.e
$ h_{\infty}^{j} = 0$. The proof is essentially well-known: see Lemma $4.2$ in  \cite{ibramasmnakerr}. We write down a proof for
convenience of the reader in Section \ref{Section:ProofLemmaConv}.

\begin{lem} (see also Lemma $4.2$ in \cite{ibramasmnakerr})
Assume that $h_{\infty}^{j} = 0$. Let $\vec{U}_{\infty}^{j}$  be defined by the following

\begin{equation}
\begin{array}{ll}
\vec{U}_{\infty}^{j} & := e^{ i t D } \vec{\phi}^{j} - i \int_{\tau_{\infty}^{j}}^{t} e^{i(t-s) D }
F \left( U^{j}_{\infty} (s) \right)   \, ds,
\end{array}
\nonumber
\end{equation}
with $ U_{\infty}^{j}:= \Re \left( D^{-1} \vec{U}_{\infty}^{j} \right)$. %Let

%\begin{equation}
%\begin{array}{l}
%\tilde{u}_{n}^{j}:= \Re \left( \langle D \rangle^{-1} T_{n}^{j} \vec{U}_{n}^{j} \left( \frac{t}{h_{n}^{j}} \right) \right) =
%u_{n}^{j} ( t + t_{n}^{j} )
%\end{array}
%\nonumber
%\end{equation}
Assume that $\| U_{\infty}^{j} \|_{L_{t}^{\frac{d+2)}{d-2}}  L_{x}^{\frac{2(d+1)}{d-2}} (\mathbb{R})} < \infty$. Then

\begin{equation}
\begin{array}{l}
K \subset \mathbb{R} \; \text{bounded}: \; \lim_{n \rightarrow \infty} \left\| \vec{U}_{n}^{j} - \vec{U}_{\infty}^{j} \right\|_{L_{t}^{\infty} L_{x}^{2} (K)} +
\left\| U_{n}^{j} - U_{\infty}^{j} \right\|_{L_{t}^{\frac{d+2}{d-2}} L_{x}^{\frac{2(d+2)}{d-2}}(K)}  = 0, \; \text{and} \\
\\
\lim_{n \rightarrow \infty}  \left\| U_{n}^{j} - U_{\infty}^{j} \right\|_{L_{t}^{\frac{d+2}{d-2}} L_{x}^{\frac{2(d+2)}{d-2}}(\mathbb{R})}  = 0  \cdot
\end{array}
\label{Eqn:ConvStrichNormNrj}
\end{equation}

\label{lem:Convhinf}
\end{lem}

\begin{rem}
Notice that if $h_{\infty}^{j}=1$  then all the estimates in this lemma clearly hold, with $\vec{U}_{\infty}^{j}$ substituted for $\vec{U}_{n}^{j}$.
\end{rem}

\begin{rem}

Observe from (\ref{Eqn:StrichWave0Sc}), Lemma \ref{lem:Convhinf}, and the previous remark that

\begin{equation}
\begin{array}{ll}
\left\| \vec{U}_{n}^{j} (\tau_{n}^{j} )  - e^{i \tau_{n}^{j} \langle D \rangle_{n}^{j}} \vec{\phi}^{j} \right\|_{L^{2}} & \lesssim
\| U_{\infty}^{j} \|^{2^{*}-1}_{L_{t}^{\frac{d+2}{d-2}} L_{x}^{\frac{2(d+2)}{d-2}}( \tau_{n}^{j}, \tau_{\infty}^{j} )} \\
& \rightarrow_{ n \rightarrow \infty} 0
\end{array}
\nonumber
\end{equation}
Since $ T_{n}^{j} e^{i \tau_{n}^{j} \langle D \rangle_{n}^{j}} = e^{- i t_{n}^{j} \langle D \rangle } T_{n}^{j} $  and the $L^{2}$ norm
is invariant by application of $T_{n}^{j}$ we get

\begin{equation}
\begin{array}{ll}
\| \vec{u}_{n}^{j}(0) - e^{- i t_{n}^{j} \langle D \rangle } T_{n}^{j} \vec{\phi}^{j}   \|_{L^{2}} & \rightarrow_{n \rightarrow \infty} 0
\end{array}
\label{Eqn:ConvergInit}
\end{equation}

\end{rem}
Let $A_{1}:= \left\{ j \in [0..k]: h_{\infty}^{j} = 0 \right\}$ and $A_{2}:= \left\{ j \in [0..k]: h_{\infty}^{j}= 1 \right\}$. \\
\\
Next we prove the following result:
\begin{res}
if $h_{\infty}^{j} = 1$ (resp. $h_{\infty}^{j} =0$) then $ K( U^{j}_{n}(t) ) > 0$ (resp. $K( U_{\infty}^{j} (t)) > 0$ )
in a neighborhood of $\tau_{\infty}^{j}$. We also have $K(w_{n}^{k}) \geq 0$ (with $ w_{n}^{k}:= \Re \left( \langle  D \rangle^{-1} \vec{w}_{n}^{k}
\right) $).
\label{Res:PositivK}
\end{res}

\begin{proof}

Assume that $\tau_{\infty}^{j} \neq \{ - \infty, \infty \}$. By (\ref{Eqn:RelG}), (\ref{Eqn:LowerBdKinet}), Proposition \ref{Prop:VarEst} and (\ref{Eqn:DecouplKinetNrj})

\begin{equation}
\begin{array}{ll}
E_{wa}(\vec{W}) - \bar{\epsilon}^{2} & \geq E \left( \vec{u}_{n}(0) \right) -  \frac{1}{2^{*}} K(u_{n}(0))  + 2 \bar{\epsilon}^{2}  \\
& \geq  \frac{ \left\| \vec{u}_{n}(0) \right\|^{2}_{L^{2}}}{d} + 2 \bar{\epsilon}^{2} \\
& \geq \sum_{j=0}^{k} \frac{ \left\| \vec{v}_{n}^{j}(0) \right\|^{2}_{L^{2}}}{d} +  \frac{ \| \vec{w}_{n}^{k} \|^{2}_{L^{2}}}{d} + \bar{\epsilon}^{2} \\
& \geq \sum_{j=0}^{k} \frac{ \left\| e^{i \tau_{n}^{j} \langle D \rangle_{n}^{j}  } \vec{\phi}^{j} \right\|^{2}_{L^{2}} }{d} +  \frac{ \| \vec{w}_{n}^{k} \|^{2}_{L^{2}}}{d}
+ \bar{\epsilon}^{2} \\
& \geq \sum \limits_{j \in A_{1}} G \left( \Re \left( D^{-1} e^{i \tau_{\infty}^{j}  D } \vec{\phi}^{j}  \right) \right) +
\sum \limits_{j \in A_{2}}
G \left( \Re \left( \langle D \rangle^{-1}
e^{ i \tau_{\infty}^{j} \langle D \rangle} \vec{\phi}^{j}  \right) \right) + G(w_{n}^{k})
\end{array}
\nonumber
\end{equation}
By (\ref{Eqn:Charac3GdState}) we see that $K(w_{n}^{k}) \geq 0$.

Therefore, if $ \Re \left( D^{-1} e^{i \tau_{\infty}^{j} D } \vec{\phi}^{j}   \right) \neq 0 $
(resp. $ \Re \left( \langle D \rangle ^{-1} e^{i \tau_{\infty}^{j}  \langle D \rangle } \vec{\phi}^{j}   \right) \neq 0 $), then
$ K \left( \Re \left( D^{-1}
e^{i \tau_{\infty}^{j} D} \vec{\phi}^{j}\right) \right) > 0 $ (resp. $  K \left( \Re \left( \langle D \rangle ^{-1} e^{i \tau_{\infty}^{j}  \langle D \rangle } \vec{\phi}^{j}   \right) \right) > 0 $ ) by (\ref{Eqn:Charac3GdState}). Hence we see from Remark \ref{Rem:UnjProp},  Remark \ref{Rem:UinfjProp}, and (\ref{Eqn:AsympBehUn}), that $K (U_{\infty}^{j}(t)) > 0$
(resp. $K(U_{n}^{j}(t)) > 0$. If  $ \Re \left(  D^{-1} e^{i  \tau_{\infty}^{j} D}   \vec{\phi}^{j}   \right) = 0 $ (resp. $ \Re \left( \langle D \rangle^{-1}
e^{i  \tau_{\infty}^{j} \langle D \rangle}   \vec{\phi}^{j}   \right) = 0 $ ), then
 $ \left\| U_{\infty}^{j} (t) \right\|_{\dot{H}^{1}} \rightarrow_{t \rightarrow \tau_{\infty}^{j}} 0 $
(resp. $\left\| U_{n}^{j} (t) \right\|_{H^{1}} \rightarrow_{t \rightarrow \tau_{\infty}^{j}} 0 $)
and, consequently, we see, by (\ref{Eqn:SobIneq}), that the results also holds. \\
If $\tau_{\infty}^{j} = \pm \infty$ then we claim that

\begin{equation}
\begin{array}{ll}
h_{\infty}^{j} =1: & \left\| \Re \left(  \langle D \rangle^{-1} e^{i t \langle D \rangle} \vec{\phi}^{j} \right) \right\|_{L^{2^{*}}}
\rightarrow_{t \rightarrow \pm \infty} 0  \\ \\
h_{\infty}^{j} =0: & \left\| \Re \left( D^{-1} e^{i t D} \vec{\phi}^{j} \right) \right\|_{L^{2^{*}}}
\rightarrow_{t \rightarrow \pm \infty} 0
\end{array}
\label{Eqn:DecayPotZero}
\end{equation}
Assume that $h_{\infty}^{j} = 1$ (resp. $h_{\infty}^{j}=0$ ). Then the claim easily follows for a smooth function $g$: it is enough to use the
standard dispersive estimate (see for example \cite{ginebvelo})

\begin{equation}
\begin{array}{ll}
\left\| \langle D \rangle^{-1} e^{it \langle D \rangle} g \right\|_{B_{2^{*},2}^{0}} & \lesssim  \frac{1}{|t|^{\frac{d-1}{d}}} \| g \|_{B_{(2^{*})^{'},2}^{\frac{1}{d}}},
\end{array}
\label{Eqn:DispEstPr}
\end{equation}
( resp.  $ \left\| D^{-1} e^{it D} g \right\|_{\dot{B}_{2^{*},2}^{0}} \lesssim  \frac{1}{|t|^{\frac{d-1}{d}}} \| g \|_{\dot{B}^{\frac{1}{d}} _{(2^{*})^{'},2}} $  )
and the embedding $ B^{0}_{2^{*},2} \hookrightarrow L^{2^{*}}$ (resp. $\dot{B}^{0}_{2^{*},2} \hookrightarrow L^{2^{*}}$ ). The claim
for $ \vec{\phi}^{j} \in L^{2}$ follows from $L^{2}$ approximation by a smooth function and (\ref{Eqn:SobIneq}), using a standard procedure. \\
From (\ref{Eqn:DecayPotZero}), Remark \ref{Rem:UnjProp}, and Remark \ref{Rem:UinfjProp}, we see that the result holds.

\end{proof}

We prove the following result:

\begin{res}
\begin{equation}
\begin{array}{ll}
E(\vec{u}_{n}) & = \sum \limits_{ j \in A_{1}}  E_{wa}(\vec{U}^{j}_{\infty}) + \sum \limits_{j \in A_{2}} E( \vec{U}^{j}_{n}) +
E (\vec{w}^{k}_{n}) + O(\alpha)
\end{array}
\label{Eqn:DecouplNRJ}
\end{equation}
\label{Res:DecouplNrj}
\end{res}

\begin{proof}

From (\ref{Eqn:DecouplKinetNrj}), Remark \ref{Rem:UnjProp}, and Remark \ref{Rem:UinfjProp} we see that

\begin{equation}
\begin{array}{ll}
\frac{1}{2}  \| \vec{u}_{n}(0) \|^{2}_{L^{2}} & = \frac{1}{2} \sum_{j=0}^{k} \left\| e^{- i \langle D \rangle t_{n}^{j} } T_{n}^{j} \vec{\phi}^{j}  \right\|^{2}_{L^{2}} +  \frac{1}{2} \left\| \vec{w}^{k}_{n} \right\|_{L^{2}}^{2} + O(\alpha) \\
& = \frac{1}{2} \sum_{j=0}^{k} \left\| e^{i \tau_{n}^{j}  \langle D \rangle_{n}^{j} } \vec{\phi}^{j}  \right\|^{2}_{L^{2}}
+ \frac{1}{2} \left\| \vec{w}^{k}_{n} \right\|_{L^{2}}^{2} + O(\alpha) \\
& = \frac{1}{2} \sum_{j \in A_{1}}  \left\| e^{ i \tau_{n}^{j} D } \vec{\phi}^{j}  \right\|^{2}_{L^{2}} +
\frac{1}{2} \sum_{j \in A_{2}} \left\| e^{ i \tau_{n}^{j} \langle D \rangle } \vec{\phi}^{j}  \right\|^{2}_{L^{2}} +
\frac{1}{2} \left\| \vec{w}^{k}_{n} \right\|_{L^{2}}^{2} + O(\alpha) \\
& = \frac{1}{2} \sum_{j \in A_{1}} \left\| \vec{U}^{j}_{\infty} (\tau_{n}^{j}) \right\|^{2}_{L^{2}} + \frac{1}{2} \sum_{j \in A_{2}} \left\| \vec{U}_{n}^{j}(\tau_{n}^{j})  \right\|^{2}_{L^{2}}  + \frac{1}{2} \left\| \vec{w}^{k}_{n} \right\|^{2}_{L^{2}} + O(\alpha) \cdot
\end{array}
\label{Eqn:DecoupKinet}
\end{equation}
It remains to prove the decoupling of the potential part of the energy, i.e

\begin{equation}
\begin{array}{ll}
\tilde{F} \left( u_{n}(0) \right) & = \sum \limits_{j \in A_{1}} \tilde{F} \left( U^{j}_{\infty} (\tau_{n}^{j}) \right)
+ \sum \limits_{j \in A_{2}} \tilde{F} \left( U_{n}^{j} (\tau_{n}^{j} )  \right) + \tilde{F}(w_{n}^{k}) + O(\alpha),
\end{array}
\nonumber
\end{equation}
with $\tilde{F}(g) :=  \frac{\| g \|^{2^{*}}_{L^{2^{*}}}}{2^{*}}$. \\
From  $B_{\infty, \infty}^{-\frac{d}{2} + 1} \cap B^{1}_{2,2} \hookrightarrow B^{0}_{2^{*},1}  $ (see e.g. Proposition $2.22$ in \cite{bahchem}
\footnote{the proposition is stated in \cite{bahchem} in homogeneous Besov spaces but it also clearly holds in inhomogeneous Besov spaces \label{Foot:HomInhom}} ), $ B^{0}_{2^{*},1} \hookrightarrow L^{2^{*}}$, (\ref{Eqn:Asympwnk1}), and (\ref{Eqn:DecouplKinetNrj}), we see that $ \tilde{F}(w_{n}^{k}) = O(\alpha)$.\\
Assume that $\tau_{\infty}^{j} = \pm \infty$. Then in view of (\ref{Eqn:DecayPotZero}, Remark \ref{Rem:UnjProp}, and Remark \ref{Rem:UinfjProp}, we get

\begin{equation}
\begin{array}{ll}
h_{\infty}^{j} =1: & \tilde{F} \left( U_{n}^{j}(\tau_{n}^{j}) \right)  \lesssim
\tilde{F} \left( \Re \left(   \langle D \rangle^{-1} e^{i \tau_{n}^{j} \langle D \rangle} \vec{\phi}^{j} \right) \right) + O(\alpha)  = O (\alpha) \\
h_{\infty}^{j}=0: & \tilde{F} \left( U_{\infty}^{j}(\tau_{n}^{j}) \right)  \lesssim
\tilde{F} \left( \Re \left(  D^{-1} e^{i \tau_{n}^{j} D} \vec{\phi}^{j} \right) \right) + O(\alpha)  = O (\alpha)
\end{array}
\nonumber
\end{equation}
Assume that $h_{\infty}^{j}=0$. We  claim that

\begin{equation}
\begin{array}{l}
\left\| \Re \left( (\langle D \rangle)^{-1} e^{-i t_{n}^{j} \langle D \rangle} T_{n}^{j} \vec{\phi}^{j} \right) \right\|_{L^{2^{*}}}  = O (\alpha) \cdot
\end{array}
\label{Eqn:ResPot}
\end{equation}
Indeed first observe that if $\vec{\phi}^{j}_{\alpha}$ Schwartz function such that $ \|\vec{\phi}^{j} - \vec{\phi}^{j}_{\alpha} \|_{L^{2}} = O(\alpha) $,
then we get from (\ref{Eqn:SobIneq}) and Plancherel theorem that
$ \left\| \Re \left( (\langle D \rangle)^{-1} e^{-i t_{n}^{j} \langle D \rangle} T_{n}^{j} \vec{\phi}^{j} \right) \right\|_{L^{2^{*}}} \lesssim
\|\vec{\phi}^{j} - \vec{\phi}^{j}_{\alpha} \|_{L^{2}} = o(\alpha)$. Hence it suffices to prove that (\ref{Eqn:ResPot}) with
$\vec{\phi}^{j}$ Schwartz. But then, we get from (\ref{Eqn:DispEstPr}) and
$ L^{(2^{*})^{'}} \hookrightarrow B_{(2^{*})^{'},2}^{0} $

\begin{equation}
\begin{array}{ll}
\left\| \Re \left( (\langle D \rangle)^{-1} e^{-i t_{n}^{j} \langle D \rangle} T_{n}^{j} \vec{\phi}^{j} \right) \right\|_{L^{2^{*}}} & \lesssim
\frac{1}{|h_{n}^{j} \tau_{n}^{j}|\frac{d-1}{d}} \| \langle D \rangle^{1 + \frac{1}{d}} T_{n}^{j} \vec{\phi}^{j} \|_{L^{(2^{*})^{'}}} \\
& \lesssim \frac{1}{h_{n}^{j} |\tau_{n}^{j}|^{\frac{d-1}{d}}}
\left\| T_{n}^{j} \left( \langle D \rangle_{n}^{j} \right)^{\frac{1}{d}} \vec{\phi}^{j} \right\|_{L^{(2^{*})^{'}}} \\
& \lesssim \frac{1}{|\tau_{n}^{j}|^{\frac{d-1}{d}}}  \left\| \langle D \rangle^{\frac{1}{d}} \vec{\phi}^{j} \right\|_{L^{(2^{*})^{'}}},
\end{array}
\nonumber
\end{equation}
which is $o(\alpha)$. % \color{red}

%COMMENTS: this was the previous way to estimate $\tilde{F} \left( \Re \left( \langle D \rangle^{-1} e^{- i t_{n}^{j} \langle D \rangle} T_{n}^{j} \vec{\phi}^{j} \right) \right)$ %:

%\begin{equation}
%\begin{array}{ll}
%\tilde{F} \left( \Re \left( \langle D \rangle^{-1} e^{- i t_{n}^{j} \langle D \rangle} T_{n}^{j} \vec{\phi}^{j} \right)   \right)
%&  = \tilde{F} \left( \Re \left( (\langle D \rangle_{n}^{j})^{-1} e^{i \tau_{n}^{j} \langle D \rangle_{n}^{j}} \vec{\phi}^{j} \right) \right) \\
%& \lesssim
%\tilde{F} \left(  \left( (\langle D \rangle_{n}^{j})^{-1} - D^{-1} \right) \Re \left( e^{i \tau_{n}^{j} \langle D \rangle_{n}^{j}} \vec{\phi}^{j}
% \right) \right)  \\
%&  +  \tilde{F} \left( \Re \left( D^{-1} e^{i \tau_{n}^{j} \langle D \rangle_{n}^{j}} \vec{\phi}^{j}   \right) \right) \\
%& \lesssim  \left\| \langle \frac{D}{h_{n}^{j}} \rangle^{-2} \vec{\phi}^{j}  \right\|_{L^{2}}   +  \tilde{F} \left( \Re \left( D^{-1} e^{i \tau_{n}^{j}
%\langle D \rangle_{n}^{j}} \vec{\phi}^{j}   \right) \right) \\
%& = O(\alpha),
%\end{array}
%\label{Eqn:ResPot}
%\end{equation}

%END COMMENTS

%\color{black}

Now assume that $h_{\infty}^{j} =1$. Then (\ref{Eqn:ResPot}) also holds from (\ref{Eqn:DecayPotZero}). \\
Hence (in view of (\ref{Eqn:Decompvn})) it is enough to prove that

\begin{equation}
\begin{array}{ll}
\left|
\begin{array}{l}
\tilde{F} \left( \sum  \limits_{j \in [0,k] : \tau_{\infty}^{j} \in \mathbb{R}  } \Re \left( \langle D \rangle^{-1} T_{n}^{j} e^{i \tau_{n}^{j} \langle D \rangle_{n}^{j}}  \vec{\phi}^{j} \right) \right)
- \sum \limits_{j \in A_{1}: \tau_{\infty}^{j} \in \mathbb{R}} \tilde{F} \left( \Re  (D^{-1} T_{n}^{j} e^{ i \tau_{n}^{j} D } \vec{\phi}^{j} ) \right) \\
- \sum \limits_{j \in A_{2}: \tau_{\infty}^{j} \in \mathbb{R}} \tilde{F} \left( \Re ( \langle D \rangle^{-1} T_{n}^{j} e^{i \tau_{n}^{j} \langle D \rangle}
\vec{\phi}^{j} ) \right)
\end{array}
\right|
& = O(\alpha)
\end{array}
\nonumber
\end{equation}
Assume that $h_{\infty}^{j}=0$. We see from the Plancherel theorem, the equality
$\langle D \rangle^{-1} T_{n}^{j} = h_{n}^{j} T_{n}^{j} ( \langle D \rangle_{n}^{j} )^{-1}$, (\ref{Eqn:SobIneq}),
and the estimate   $ \left\| ( \langle D \rangle_{n}^{j})^{-1} g - D^{-1} g \right\|_{L^{2^{*}}}
\lesssim \left\| \langle \frac{D}{h_{n}^{j}} \rangle^{-2} D^{-1} g \right\|_{L^{2^{*}}} $,
that follows from the Hormander-Mikhlin multiplier that

\begin{equation}
\begin{array}{l}
\tilde{F} \left( \Re \left( \langle D \rangle^{-1} T_{n}^{j} e^{i \tau_{n}^{j} \langle D \rangle_{n}^{j}}  \vec{\phi}^{j}
- h_{n}^{j} T_{n}^{j}  D^{-1} e^{i \tau_{\infty}^{j} D } \vec{\phi}^{j} \right) \right) \\
\lesssim \tilde{F} \left( \Re \left(  \langle D \rangle^{-1} T_{n}^{j}
(e^{i \tau_{n}^{j} \langle D \rangle_{n}^{j}} - e^{i \tau_{\infty}^{j} D })  \vec{\phi}^{j} \right) \right)
+ \tilde{F} \left(  h_{n}^{j} T_{n}^{j} \left(  (\langle D \rangle_{n}^{j})^{-1} - D^{-1}  \right)
\Re \left(  e^{i \tau_{\infty}^{j} D } \vec{\phi}^{j} \right) \right) \\
\lesssim \left\| (e^{i \tau_{n}^{j} \langle D \rangle_{n}^{j}} - e^{i \tau_{\infty}^{j} D  })  \vec{\phi}^{j} \right\|_{L^{2}}
+  \left\| \langle \frac{D}{h_{n}^{j}} \rangle^{-2} \vec{\phi}^{j} \right\|_{L^{2}} \\
= O(\alpha) \cdot
% - (\langle D \rangle^{j}_{n})^{-1} T_{n}^{j} e^{i \tau_{n}^{j} \langle D \rangle_{n}^{j}}  \vec{\phi}^{j}  \right)   \right\|_{L^{2^{*}}}
\end{array}
\label{Eqn:F1}
\end{equation}
We also get from $D^{-1} T_{n}^{j} = h_{n}^{j} T_{n}^{j} D^{-1}$ the estimate below

\begin{equation}
\begin{array}{l}
\left| \tilde{F} \left(  \Re \left( D^{-1} T_{n}^{j} e^{ i \tau_{n}^{j} D } \vec{\phi}^{j} \right) \right) -
\tilde{F} \left(  \Re \left( D^{-1} T_{n}^{j} e^{ i \tau_{\infty}^{j} D } \vec{\phi}^{j} \right)   \right) \right| \\
\lesssim  \tilde{F} \left( h_{n}^{j} T_{n}^{j} D^{-1}  \Re \left(  ( e^{i \tau_{n}^{j}} - e^{i \tau_{\infty}^{j}} ) \vec{\phi}^{j}   \right) \right) \\
\lesssim  \| ( e^{i \tau_{n}^{j} D} - e^{i \tau_{\infty}^{j} D} ) \vec{\phi}^{j} \|_{L^{2}} \\
= O(\alpha)
 \end{array}
\label{Eqn:F2}
\end{equation}
If $h_{\infty}^{j}=1$, then it is easier to see from (\ref{Eqn:SobIneq}) that (\ref{Eqn:F1}) and (\ref{Eqn:F2}) hold: the proof is left to
the reader.\\
Hence letting

\begin{equation}
\tilde{\phi}^{j} :=
\left\{
\begin{array}{l}
\Re \left(  \langle D \rangle^{-1} e^{i \tau_{\infty}^{j} \langle D \rangle } \vec{\phi}^{j}
\right), h_{\infty}^{j}=1 \\
\Re  \left(  D^{-1} e^{i \tau_{\infty}^{j} D  } \vec{\phi}^{j}
\right), h_{\infty}^{j}=0,
\end{array}
\right.
\nonumber
\end{equation}
it is enough to prove

\begin{equation}
\begin{array}{ll}
\left|
F \left( \sum_{j \in [0,k]: \tau_{\infty}^{j} \in \mathbb{R} } h_{n}^{j} T_{n}^{j} \tilde{\phi}^{j}    \right)
- \sum_{j \in [0,k]: \tau_{\infty}^{j} \in \mathbb{R}} F \left( h_{n}^{j} T_{n}^{j} \tilde{\phi}^{j} \right) \right|
& = O(\alpha)
\end{array}
\label{Eqn:DecouplFSimpl}
\end{equation}
We may assume WLOG that $\tilde{\phi}^{j}$ is compactly supported. Let $h_{n}^{j,l}:= \frac{h_{n}^{l}}{h_{n}^{j}}$ and
$x_{n}^{j,l} := \frac{x_{n}^{l} - x_{n}^{j}}{h_{n}^{j}}$. From (\ref{Eqn:OrthParam}) we see that

\begin{equation}
\begin{array}{l}
h_{n}^{j,l} \rightarrow 0, \, or \, h_{n}^{j,l} \rightarrow \infty, \, or \, |x_{n}^{j,l}| \rightarrow \infty,
\end{array}
\label{Eqn:OrthTransl}
\end{equation}
From the elementary estimate  $  \left| \left| \sum \limits_{i=0}^{m} x_{i} \right|^{2^{*}} - \sum  \limits_{i=0}^{m} |x_{i}|^{2^{*}} \right| \lesssim
\sum \limits_{ (i,l) \in [0..m]^{2}: i \neq l} |x_{i}|^{2^{*} -1} |x_{l}| $ (that can be proved by induction on $m$) we see that

\begin{equation}
\begin{array}{l}
\left|
F \left( \sum \limits_{j \in [0,k]: \tau_{\infty}^{j} \in \mathbb{R} } h_{n}^{j} T_{n}^{j} \tilde{\phi}^{j}    \right)
- \sum \limits_{j \in [0,k]: \tau_{\infty}^{j} \in \mathbb{R}} F \left( h_{n}^{j} T_{n}^{j} \tilde{\phi}^{j} \right) \right| \\
\lesssim \\  \sum \limits_{\substack{(j,l) \in [0,k]^{2} \\ j \neq l \\ (\tau_{\infty}^{j} ,\tau_{\infty}^{l}) \in \mathbb{R}^{2}}}
\int_{\mathbb{R}^{d}}  (h_{n}^{j})^{2^{*} -1} h_{n}^{l}  | T_{n}^{j} \tilde{\phi}^{j} |^{2^{*}-1} | T_{n}^{l} \tilde{\phi}^{l}| \; dx \\
\lesssim \sum \limits_{\substack{(j,l) \in [0,k]^{2} \\
j \neq l, (\tau_{\infty}^{j} ,\tau_{\infty}^{l}) \in \mathbb{R}^{2}, h_{n}^{j,l} < 1 } } (h_{n}^{j,l})^{\frac{d-2}{2}}
\int_{\mathbb{R}^{d}}  | \tilde{\phi}^{j} |^{2^{*} -1} (x)  |\tilde{\phi}^{l}| \left( \frac{x - x_{n}^{j,l}}{h_{n}^{j,l}}  \right) \; dx \cdot
\end{array}
\label{Eqn:OrthPotEnergy}
\end{equation}
Clearly the last quantity goes to zero as $n \rightarrow \infty$, in view of (\ref{Eqn:OrthTransl})

%\color{red}
% COMMENTS: Previous version

%We define $\tilde{\phi}_{n}^{j}(x) $ in the following fashion

%\begin{equation}
%\begin{array}{ll}
%\tilde{\phi}_{n}^{j}(x) & := \tilde{\phi}^{j}(x) \times
%\left\{
%\begin{array}{l}
%0, \;  \bigcup_{j \neq l: h_{n}^{j,l} < 1} \left\{  x:  \frac{x- x_{n}^{j,l}}{h_{n}^{j,l}} \in supp \left( \tilde{\phi}^{l} \right) \right\} \\
%1, \; elsewhere
%\end{array}
%\right.,
%\end{array}
%\nonumber
%\end{equation}

%as $n \rightarrow \infty$. The dominated convergence theorem implies that
%$ \tilde{\phi}_{n}^{j} \rightarrow \tilde{\phi}^{j} $  as $ n \rightarrow \infty $ in $L^{2^{*}}$ and, consequently, it is enough
%to prove (\ref{Eqn:DecouplFSimpl}) with $\tilde{\phi}^{j}$ replaced with $\tilde{\phi}_{n}^{j}$.
%Observe that $ \supp(\tilde{\phi}_{n}^{j}) \cap \supp(\tilde{\phi}_{n}^{k}) = \emptyset $ if $j \neq k$. Hence (\ref{Eqn:DecouplFSimpl}) holds.
% END COMMENTS
% \color{black}

\end{proof}

Letting $\alpha \rightarrow 0$ ( hence $k \rightarrow \infty$ and $n \rightarrow \infty)$ in
(\ref{Eqn:DecouplNRJ}) we see that

\begin{equation}
\begin{array}{l}
h_{\infty}^{j}=1: \, E(\vec{U}_{n}^{j}) \leq E_{c} \\
h_{\infty}^{j}=0: \, E_{wa}(\vec{U}_{\infty}^{j}) \leq E_{c}  \\
\end{array}
\nonumber
\end{equation}
If $E(\vec{U}_{\infty}^{j}) < E_{wa}(\vec{W})$ for all $j$ such that $h_{\infty}^{j} =0$ and $E(\vec{U}_{n}^{j}) < E_{wa} (\vec{W})$ for all $j$ such that
$h_{\infty}^{j}=1$ then we conclude from \cite{kenmerwave} and \cite{ibramasmnak}  that $U^{j}_{\infty}$ and $U_{n}^{j}$  exist globally in time
and scatters, i.e $ \| U^{j}_{\infty} \|_{L_{t}^{\frac{2(d+1)}{d-2}} L_{x}^{\frac{2(d+1)}{d-2}} (\mathbb{R}) } < \infty$ and
$ \| U^{j}_{n} \|_{L_{t}^{\frac{2(d+1)}{d-2}} L_{x}^{\frac{2(d+1)}{d-2}} (\mathbb{R}) } < \infty $. Next we prove a
perturbation result:

\begin{res}
Let $w$ be a solution of (\ref{Eqn:NlkgCrit}) and $\tilde{w}$ be a solution of

\begin{equation}
\begin{array}{ll}
\partial_{tt} \tilde{w} - \triangle \tilde{w} + \tilde{w} & = |\tilde{w}|^{\frac{4}{d-2}} \tilde{w} + eq (\tilde{w})
\end{array}
\nonumber
\end{equation}
Let $t_{0} \in  \left( -T_{-}(w), T_{+}(w) \right)$ and $I:= [t_{0}, \cdot] $.
Assume that there exists $C_{1}$ such that $\| \tilde{w} \|_{L_{t}^{\frac{2(d+1)}{d-2}} L_{x}^{\frac{2(d+1)}{d-2}}(I) }  \leq C_{1}$. Then there exists $\epsilon_{0}:= \epsilon_{0}(C_{1})$ such that if $\epsilon^{'} \leq \epsilon_{0}$ and $(a)$ (resp. $(a')$) hold, with

\begin{equation}
\begin{array}{l}
(a): \;  \left\| \cos \left( (t - t_{0}) \langle D \rangle \right) \left( w(t_{0}) - \tilde{w}(t_{0}) \right)
+ \frac{\sin \left( (t-t_{0}) \langle D \rangle  \right)}{ \langle D \rangle}  \left( \partial_{t} w(t_{0}) - \partial_{t} \tilde{w}(t_{0}) \right)
\right\|_{L_{t}^{\frac{d+2}{d-2}} L_{x}^{\frac{2(d+2)}{d-2}} (I)} \leq \epsilon^{'} \, \text{and}  \\
\| eq(\tilde{w}) \|_{L_{t}^{1} L_{x}^{2} (I)}  \leq \epsilon^{'} \\
\\
(a'): \; \left\| \vec{w}(t_{0}) - \vec{\tilde{w}}(t_{0}) \right\|_{L^{2}} \leq \epsilon^{'} \, \text{and} \,  \| eq(\tilde{w}) \|_{L_{t}^{1} L_{x}^{2} (I)}  \leq \epsilon^{'},
\end{array}
\label{Eqn:SmallnessAss}
\end{equation}
then $I \subset (- T_{-}(w),T_{+}(w))$ and $(b)$ (resp. $(b')$) hold, with

\begin{equation}
\begin{array}{l}
(b): \;  \| w - \tilde{w} \|_{L_{t}^{\frac{2(d+1)}{d-2}} L_{x}^{\frac{2(d+1)}{d-2}}(I)}= O (\epsilon^{'})   \\
(b'): \;   \| w - \tilde{w} \|_{L_{t}^{\frac{2(d+1)}{d-2}} L_{x}^{\frac{2(d+1)}{d-2}}(I)}= O (\epsilon^{'}), \; \text{and} \;
\| \vec{w} - \vec{\tilde{w}} \|_{ L_{t}^{\infty} L_{x}^{2} (I)} = O(\epsilon^{'})
\end{array}
\label{Eqn:CloseNormP}
\end{equation}
Here $\vec{w}:= \langle D \rangle w - i \partial_{t} w$ and  $\overrightarrow{\tilde{w}}:= \langle D \rangle \tilde{w} - i \partial_{t} \tilde{w}$.

\label{res:PertubationResult}
\end{res}

\begin{proof}
We only prove the result if $(b)$ holds. Indeed, if $(a)$ holds, then the proof of the result is a slight modification of that
if $(b)$ holds. \\
Let $1 \gg \epsilon^{'}_{0} > 0$ small enough such that all the statements below are true. Let $F(v):= |v|^{\frac{4}{d-2}} v$ and
$e:= w - \tilde{w}$ and $ 0 < \alpha \ll 1$. The proof is made of three steps:

\begin{itemize}

\item $ \| \tilde{w} \|_{ L_{t}^{\frac{d+2}{d-2}} L_{x}^{\frac{2(d+2)}{d-2}}(\mathbb{R})} < \infty$. \\
Indeed, dividing $I$ into subintervals $J=[a,.]$ such that $\| \tilde{w} \|_{L_{t}^{\frac{2(d+1)}{d-2}} L_{x}^{\frac{2(d+1)}{d-2}} (J) } \approx \alpha$, we get
from (\ref{Eqn:StrichKg0}) and (\ref{Eqn:NonlinearEstGen})

\begin{equation}
\begin{array}{l}
\max{  \left( \| \tilde{w} \|_{L_{t}^{\frac{d+2}{d-2}} L_{x}^{\frac{2(d+2)}{d-2}} (J) }, \| \overrightarrow{\tilde{w}} \|_{L_{t}^{\infty} L_{x}^{2} (J)},
\| \tilde{w} \|_{ L_{t}^{\frac{2(d+1)}{d-1}} B^{\frac{1}{2}}_{\frac{2(d+1)}{d-2},2} (J) } \right)} \\
\lesssim  \| \overrightarrow{\tilde{w}}(a) \|_{L^{2}} + \alpha^{2^{*} -2} \| \tilde{w} \|_{L_{t}^{\frac{2(d+1)}{d-1}} B^{\frac{1}{2}}_{\frac{2(d+1)}{d-2},2} (J) } \\
\end{array}
\label{Eqn:EstEachJ}
\end{equation}
By a continuity argument applied to (\ref{Eqn:EstEachJ}) on each $J$ and by iteration, we see that the claim holds.

\item short-time perturbation argument \\
Let $ 1 \gg \alpha > 0$ and $J=[a,] \subset I$ such that $\| \tilde{w} \|_{ L_{t}^{\frac{d+2}{d-2}} L_{x}^{\frac{2(d+2)}{d-2}} (J) } \leq \alpha$.
There exists $ 1 \gg \epsilon^{''}_{0} > 0$ such that if $ \| eq(\tilde{w}) \|_{L_{t}^{1} L_{x}^{2}(J)} \leq \epsilon^{''} $,
$ \| \vec{w}(a) - \overrightarrow{\tilde{w}}(a) \|_{L^{2}} \leq \epsilon^{''} $, and $\epsilon^{''} \leq \epsilon^{''}_{0}$, then

\begin{equation}
\begin{array}{ll}
\| F(w) - F(\tilde{w}) \|_{L_{t}^{1} L_{x}^{2} (J)} \lesssim \epsilon^{''} \cdot
\end{array}
\label{Eqn:DiffFwShort}
\end{equation}
Indeed we see from (\ref{Eqn:StrichKg0}) that

\begin{equation}
\begin{array}{ll}
\| e \|_{L_{t}^{\frac{d+2}{d-2}} L_{x}^{\frac{2(d+2)}{d-2}} (J) } & \lesssim \epsilon^{''}
+  \| F (w) - F( \tilde{w} ) \|_{L_{t}^{1} L_{x}^{2} (J)} \\
& \lesssim \epsilon^{''} +  \| e \|_{ L_{t}^{\frac{d+2}{d-2}}  L_{x}^{\frac{2(d+2)}{d-2}} (J) }
 \left( \| \tilde{w} + e \|^{2^{*} -2}_{L_{t}^{\frac{d+2}{d-2}} L_{x}^{\frac{2(d+2)}{d-2}} (J)}
+ \| \tilde{w} \|^{2^{*} -2}_{L_{t}^{\frac{d+2}{d-2}} L_{x}^{\frac{2(d+2)}{d-2}} (J) }  \right) \\
& \lesssim \epsilon^{''},
\end{array}
\nonumber
\end{equation}
applying a continuity argument at the last line. Hence (\ref{Eqn:DiffFwShort}) holds.

\item long-time perturbation argument \\
We divide $I$ into subintervals $(J_{j} := [a_{j},a_{j+1}])_{1 \leq j \leq \bar{j}}$ such that
$\| \tilde{w} \|_{L_{t}^{\frac{d+2}{d-2}} L_{x}^{\frac{2(d+2)}{d-2}} (J_{j}) } = \alpha$  for $1 \leq j \leq \bar{j}-1$
and $\| \tilde{w} \|_{L_{t}^{\frac{d+2}{d-2}} L_{x}^{\frac{2(d+2)}{d-2}} (J_{\bar{j}}) } \leq \alpha $.
We then prove by induction that if $\epsilon_{0} \ll 1$ then there exists $ 0 \leq C(j) < \infty$ such that
$ (*): \; \| F(w) - F(\tilde{w}) \|_{L_{t}^{1} L_{x}^{2}(J_{j})} \leq C(j) \epsilon^{'}$. Indeed assume that
$(*) $ holds for all $ 1 \leq j \leq k$  with $k \in \{1,...,\bar{j}-1 \}$. Then by iterating (\ref{Eqn:StrichKg0}) on
$j$ we get $\| w - \tilde{w} \|_{L_{t}^{\infty} \mathcal{H} (J_{k}) } \lesssim  \epsilon^{'}
+ \sum \limits_{j=1}^{k} \| F(w) - F(\tilde{w}) \|_{L_{t}^{1} L_{x}^{2}(J_{j})} \lesssim  \epsilon^{'} +
\sum  \limits_{j=1}^{k} C(j) \epsilon^{'} \leq \epsilon^{''}_{0} $. Hence by applying the short-perturbation argument
we see that $(*)$ holds for $j=k+1$. Hence $ \| F(w) - F(\tilde{w}) \|_{L_{t}^{1} L_{x}^{2}(I)} = O (\epsilon^{'})$
and by (\ref{Eqn:StrichKg0}) we get (\ref{Eqn:CloseNormP}).

\end{itemize}

\end{proof}
Let $ w:= u_{n}$ and $ \tilde{w} := \sum \limits_{j=0}^{k} u_{n}^{j}$. \\
We first estimate $ \| \tilde{w} \|_{L_{t}^{\frac{2(d+1)}{d-2}}  L_{x}^{\frac{2(d+1)}{d-2}} (\mathbb{R})} $. In view of (\ref{Eqn:DecouplKinetNrj}),
and the global theory for small data we see that there exists finite set $J$ for which

\begin{equation}
\begin{array}{l}
j \notin J: \; \| u_{n}^{j} \|_{L_{t}^{\frac{d+2}{d-2}} L_{x}^{\frac{2(d+2)}{d-1}}(\mathbb{R})}
+ \| u_{n}^{j} \|_{L_{t}^{\frac{2(d+1)}{d-2}} L_{x}^{\frac{2(d+1)}{d-2}}(\mathbb{R})}  \lesssim
\left\| \left( u_{n}^{j}(0), \partial_{t} u_{n}^{j} (0) \right) \right\|_{\mathcal{H}}
\end{array}
\nonumber
\end{equation}
Using also the elementary estimate $ \left| \left| \sum \limits_{i=0}^{m} x_{i} \right|^{\frac{2(d+1)}{d-2}}
- \sum \limits_{i=0}^{m}  |x_{i}|^{\frac{2(d+1)}{d-2}} \right|  \lesssim \sum \limits_{(i,l) \in [0,..,m]: i \neq l}
|x_{i}|^{\frac{d+4}{d-2}} |x_{l}| $  (that can be proved by induction on $m$), and Lemma \ref{lem:Convhinf} we see that
$ \| \tilde{w} \|^{\frac{2(d+1)}{d-2}}_{L_{t}^{\frac{2(d+1)}{d-2}} L_{x}^{\frac{2(d+1)}{d-2}} (\mathbb{R})}  \lesssim  X_{1} + X_{2} + X_{3} + X_{4}$
with $ X_{1} :=  \sum \limits_{ j \in J \cap  [0,..,k]  \cap A_{1} }
\| u_{n}^{j} \|^{\frac{2(d+1)}{d-2}}_{L_{t}^{\frac{2(d+1)}{d-2}} L_{x}^{\frac{2(d+1)}{d-2}}(\mathbb{R})} $,
 $ X_{2} := \sum \limits_{j \in J \cap  [0,..,k]  \cap A_{2} } \| u_{n}^{j} \|^{\frac{2(d+1)}{d-2}}_{L_{t}^{\frac{2(d+1)}{d-2}} L_{x}^{\frac{2(d+1)}{d-2}}(\mathbb{R})} $,
 $ X_{3} := \sum \limits_{j \notin J} \| u_{n}^{j} \|^{\frac{2(d+1)}{d-2}}_{L_{t}^{\frac{2(d+1)}{d-2}} L_{x}^{\frac{2(d+1)}{d-2}}(\mathbb{R})} $, and
 $ X_{4} := \sum \limits_{\substack{(j,l) \in [0,..,k]^{2} \\ j \neq l}} \int_{\mathbb{R}^{d+1}} |u_{n}^{j}|^{\frac{d+4}{d-2}}  |u_{n}^{l}| \, dx \, dt $.
Clearly $\| u_{n}^{j} \|_{L_{t}^{\frac{2(d+1)}{d-2}} L_{x}^{\frac{2(d+1)}{d-2}}(\mathbb{R})}
= \| U_{n}^{j} \|_{L_{t}^{\frac{2(d+1)}{d-2}} L_{x}^{\frac{2(d+1)}{d-2}}(\mathbb{R})}$. Hence, using also Lemma \ref{lem:Convhinf}, we get
$ X_{1} \lesssim  \sum \limits_{ j \in J \cap  [0,..,k]  \cap A_{1}} \| U_{\infty}^{j} \|^{\frac{2(d+1)}{d-2}}_{L_{t}^{\frac{2(d+1)}{d-2}} L_{x}^{\frac{2(d+1)}{d-2}}(\mathbb{R})}
\lesssim 1$. We also have $ X_{2} \lesssim  1 $. We get from (\ref{Eqn:DecouplKinetNrj}) $ X_{3} \lesssim
\sum \limits_{j \notin J } \| \vec{v}_{n}^{j}(0) \|^{2}_{L^{2}} \ll 1 $. We have

\begin{equation}
\begin{array}{l}
X_{4} \lesssim \sum \limits_{\substack{(j,l) \in [0,..,k]^{2} \\ j \neq l}}   \| u_{n}^{j} u_{n}^{l} \|^{2}_{L_{t}^{\frac{d+1}{d-2}} L_{x}^{\frac{d+1}{d-2}}(J)}
\max \limits_{(j_{1},j_{2}) \in A_{1} \times A_{2}}
\left( \| U_{\infty}^{j_{1}} \|^{\frac{2(d+1)}{d-2} -2}_{L_{t}^{\frac{2(d+1)}{d-2}} L_{x}^{\frac{2(d+1)}{d-2}}(J)},
\| U_{n}^{j_{2}} \|^{\frac{2(d+1)}{d-2}- 2}_{L_{t}^{\frac{2(d+1)}{d-2}} L_{x}^{\frac{2(d+1)}{d-2}}(J)}
\right) \\
\sum \limits_{\substack{(j,l) \in [0,..,k]^{2} \\ j \neq l, h_{n}^{j,l} < 1 }} (h_{n}^{j,l})^{\alpha}
\left( \int_{\mathbb{R}^{d+1}} |U_{n}^{j}(t,x)|^{\frac{d+1}{d-2}} \left| U_{n}^{l} \left( \frac{t}{h_{n}^{j,l}}, \frac{x- x_{n}^{j,l}}{ h_{n}^{j,l}} \right)
\right|^{\frac{d+1}{d-2}} \; dx  \; dt \right)^{\frac{2(d-2)}{d+1}} \\
 \ll 1 \cdot
\end{array}
\end{equation}
In order to prove that the summand approaches zero as $n \rightarrow \infty$, we see from the dominated convergence theorem that we may replace WLOG $U_{n}^{j}$, $U_{n}^{l}$ with $\chi_{R} U_{n}^{j}$, $\chi_{R} U_{n}^{l} $  respectively for $R \gg 1$ (Here  $ \chi_{R} (t,x) := \chi \left( \frac{ \| (x,t) \|} {R} \right) $ with
$\chi: \mathbb{R}^{+} \rightarrow \mathbb{R} $ is a smooth function that satisfies $\chi(s)=1$  if $ s \leq 1$  and  $\chi(s)=0$ if $s \geq 2$); then we see from (\ref{Eqn:OrthTransl}) that the integral involving  $\chi_{R} U_{n}^{j}$ and  $\chi_{R} U_{n}^{l}$ goes to zero as $n \rightarrow \infty$.
Hence $  \| \tilde{w} \|_{L_{t}^{\frac{2(d+1)}{d-2}}  L_{x}^{\frac{2(d+1)}{d-2}} (\mathbb{R})} \lesssim 1 $. \\
Observe from the embedding $L^{\frac{2(d+1)}{d-2}} \hookrightarrow B^{0}_{\frac{2(d+1)}{d-2},2} $ and Proposition 2.22. in \cite{bahchem} \footnote{see Footnote \ref{Foot:HomInhom}} and (\ref{Eqn:Asympwnk1})  that

\begin{equation}
\begin{array}{l}
\left\| \langle D \rangle^{-1} e^{ \pm it \langle D \rangle} \vec{w}_{n}^{k} \right\|_{L_{t}^{\frac{2(d+1)}{d-2}} L_{x}^{\frac{2(d+1)}{d-2}} (\mathbb{R})} \ll
\epsilon_{0} \cdot
\end{array}
\nonumber
\end{equation}
Hence, using also (\ref{Eqn:Decompvn}) and  (\ref{Eqn:StrichKg0}), we see that the first estimate of (\ref{Eqn:SmallnessAss}). \\
 From the elementary estimate above we get for some $\alpha > 0 $

\begin{equation}
\begin{array}{l}
 \left\| F \left( \sum \limits_{j=0}^{k} u_{n}^{j} \right) - \sum \limits_{j=0}^{k}  F (u_{n}^{j}) \right\|_{L_{t}^{1} L_{x}^{2} (\mathbb{R})} \\
\lesssim  \sum \limits_{\substack{(j,l) \in [0,..,k]^{2} \\ j \neq l}} (h_{n}^{j,l})^{\alpha}
\left\|  u_{n}^{j}  u_{n}^{l} \right\|_{L_{t}^{\frac{d+2}{2(d-2)}} L_{x}^{\frac{d+2}{d-2}} (\mathbb{R})}
\max \limits_{(j_{1},j_{2}) \in A_{1} \times A_{2}}
\left(  \| U_{\infty}^{j_{1}} \|^{2^{*} - 3}_{L_{t}^{\frac{d+2}{d-2}} L_{x}^{\frac{2(d+2)}{d-2}} (\mathbb{R})},
\| U_{\infty}^{j_{1}} \|^{2^{*} - 3}_{L_{t}^{\frac{d+2}{d-2}} L_{x}^{\frac{2(d+2)}{d-2}} (\mathbb{R})} \right) \\
\lesssim \sum \limits_{\substack{(j,l) \in [0,..,k]^{2} \\ j \neq l, h_{n}^{j,l} < 1 }} (h_{n}^{j,l})^{\alpha}
\left( \int_{\mathbb{R}}
\left(
\int_{\mathbb{R}^{d}}  |U_{n}^{j}(t,x)|^{\frac{d+2}{d-2}} \left| U_{n}^{l} \left( \frac{t}{h_{n}^{j,l}}, \frac{x- x_{n}^{j,l}}{ h_{n}^{j,l}} \right) \right|^{\frac{d+2}{d-2}} \; dx \right)^{\frac{1}{2}} \; dt \right)^{\frac{2(d-2)}{d+2}} \\ \ll \epsilon_{0} \cdot
\end{array}
\nonumber
\end{equation}
( In order to prove the last estimate we proceed as follows. First observe that
$\| U_{\infty}^{j} \|_{L_{t}^{\frac{d+2}{d-2}} L_{x}^{\frac{2(d+2)}{d-2}} (\mathbb{R})} < \infty $ if $j \in A_{1}$ and that
$\| U_{n}^{j} \|_{L_{t}^{\frac{d+2}{d-2}} L_{x}^{\frac{2(d+2)}{d-2}} (\mathbb{R})} < \infty $  if $j \in A_{2}$
by proceeding similarly as the first stage of the proof of Result \ref{res:PertubationResult}. Then follow the same steps as those to estimate $X_{4}$.)
Hence the second estimate of (\ref{Eqn:SmallnessAss}) holds.

% \color{red}
% COMMENTS: previous version

%Now introduce for $R \gg 1$

%\begin{equation}
%\begin{array}{l}
%\hat{U}^{j}_{n,R} := \chi_{R} (t,x) U_{n}^{j}(t,x) \times
%\left\{
%\begin{array}{l}
%0, \bigcup_{j \neq l: h_{n}^{j,l} < R } \left\{ (t,x): \frac{t - t_{n}^{j,l}}{ h_{n}^{j,l}}  \leq R, \frac{x - x_{n}^{j,l}}{h_{n}^{j,l}} \leq R \right\} \\
%1, \, elsewhere
%\end{array}
%\right.,
%\end{array}
%\nonumber
%\end{equation}
%with $t_{n}^{j,l}:= \frac{t_{n}^{l} - t_{n}^{j}}{h_{n}^{j}}$ and $\chi$ defined below (\ref{Eqn:VirialIdentNLKG}). In view of Lemma \ref{lem:Convhinf}, $\hat{U}^{j}_{n,R}$ is %uniformly bounded in $L_{t}^{\frac{d+2}{d-2}} L_{x}^{\frac{2(d+2)}{d-2}} (\mathbb{R})$, and using also (\ref{Eqn:OrthTransl}), we see that $\hat{U}^{j}_{n,R} - \chi_{R} U_{n}^{j} \rightarrow 0$ in $ L_{t}^{\frac{d+2}{d-2}} L_{x}^{\frac{2(d+2)}{d-2}} (\mathbb{R}) $ as $n \rightarrow \infty$. We also have $ \chi_{R} U_{n}^{j} \rightarrow U_{n}^{j}$ as $R \rightarrow \infty$. Hence
%(invoking again  Lemma \ref{lem:Convhinf}) we may replace $u^{j}_{n}$ with $u^{j}_{n,R}:= h_{n}^{j} T_{n}^{j} \hat{U}^{j}_{n,R} \left( \frac{t -t_{n}^{j}}{h_{n}^{j}} \right)$ in %the second estimate of (\ref{Eqn:SmallnessAss}).
%Observe that $ \supp(u_{n,R}^{j}) \cap \supp (u_{n,R}^{l}) = \emptyset $ if $j \neq l$ and therefore the second estimate of
%(\ref{Eqn:SmallnessAss}) holds. \\ \\
% \color{black}

Hence applying Result \ref{res:PertubationResult}, we see that there is contradiction with the properties of $\vec{u}_{n}$ defined just below Claim \ref{Claim:CritElem}. \\ \\
%NOTE: One has to prove that the Strichartz norms are small for j large since, taking into account that one has to choose k large enough to make the remainder w_{n}^{k} small in a Strichartz norm, one would have to estimate \sum_{j=0..k, j \in A_{1}} \| U_{n}^{j} \|
This means that at least one of the $j$s satisfies  $E(\vec{U}_{\infty}^{j}) \geq E_{wa}(\vec{W})$ if $h_{\infty}^{j} =0$
(or $E(\vec{U}^{j}_{n}) \geq E_{wa}(\vec{W})$ if $h_{\infty}^{j}=1$). Assume without loss of generality that this $j$ is equal to 0. By Result
\ref{Res:PositivK}, Result \ref{Res:DecouplNrj},
(\ref{Eqn:DecouplNRJ}), the local well-posedness theory for small data we see that for the other $j$s
$T_{+}(U^{j}_{\infty})= T_{+} (U_{n}^{j}) = \infty$, $T_{-}(U^{j}_{\infty}) = T_{-} (U_{n}^{j}) = -\infty$ and

\begin{equation}
\begin{array}{ll}
\sum \limits_{j \in A_{1}, j \neq 0} \| \vec{U}_{\infty}^{j} \|^{2}_{L_{t}^{\infty} L_{x}^{2}} + \sum \limits_{j \in A_{2}, j \neq 0}
\| \vec{U}_{n}^{j} \|^{2}_{L_{t}^{\infty} L_{x}^{2}} + \| \vec{w}_{n}^{k} \|^{2}_{L^{2}} & \lesssim \bar{\epsilon}^{2}
\end{array}
\label{Eqn:EstSumL2Norm}
\end{equation}

%NOTICE
%\begin{equation}
%\begin{array}{ll}
%\vec{U}_{\infty}^{j} (- t_{n}^{j})  &  \approx  e^{- i t_{n}^{j} D} \vec{\phi}^{j} \approx T_{n}^{j}
%e^{- i \frac{t_{n}^{j}}{h_{n}^{j}} \langle D \rangle_{n}^{j}} \vec{\phi}^{j} \\
% & \approx  e^{- i t_{n}^{j} \langle D \rangle} T_{n}^{j} \vec{\phi}^{j} \; (Here f \approx g \, means \, \| f \|_{L^{2}} \approx \| g \|_{L^{2}} )
%\end{array}
%\end{equation}

Hence, applying Result \ref{res:PertubationResult} with $\tilde{w}:= u_{n}^{0}$ and $w :=u_{n}$, we see that

\begin{equation}
\begin{array}{ll}
\left\| (u_{n},\partial_{t} u_{n}) - (u^{0}_{n},\partial_{t} u_{n}^{0}) \right\|_{L_{t}^{\infty} \mathcal{H} ( [ 0, T_{+} (u^{0}_{n}) ) } & = O_{\mathcal{H}} (\bar{\epsilon}),
\end{array}
\label{Eqn:Diffunu0}
\end{equation}
Notice that $\tau_{\infty}^{0} \neq \infty$. If not, if $h_{\infty}^{0}=1$ (resp. $h_{\infty}^{0}=0$) then $U_{n}^{0}$ (resp.  $U_{\infty}^{0}$)
would scatter by construction in a neighborhood of $\infty$: see (\ref{Eqn:UnScatt}) and Remark \ref{Rem:UinfjProp}. But then, applying again Remark \ref{Rem:UinfjProp} and Result \ref{res:PertubationResult}, we see that there is contradiction. \\
Notice also that if $h_{\infty}^{1}=1$ (resp. $h_{\infty}^{0}=0$ ) then $\| U_{n}^{0} \|_{ L_{t}^{\frac{d+2}{d-2}} L_{x}^{\frac{2(d+2)}{d-2}} (\tau_{n}^{0}, T_{+}(U_{n}^{0}))} < \infty$ (resp. $ \| U_{\infty}^{0} \|_{ L_{t}^{\frac{d+2}{d-2}}  L_{x}^{\frac{2(d+2)}{d-2}}  (\tau_{n}^{0}, T_{+}(U_{\infty}^{0}))} < \infty $ ) :  if not this would imply by a standard blow-up criterion (the proof is
similar to that in \cite{kenmer}) that $T_{+}(U_{n}^{0}) = \infty$ if $h_{\infty}^{0} =1$ (or $T_{+}(U_{\infty}^{0}) = \infty$ if $h_{\infty}^{0}=0$). But we have already seen that this is not possible.\\
We claim that for $t > \tau_{\infty}^{0}$ we have

\begin{equation}
\begin{array}{l}
h_{\infty}^{0} = 0: \, \tilde{d}_{\mathcal{S}} (\vec{U}_{\infty}^{0}(t))  \geq \delta_{b}; \;  h_{\infty}^{0} =1: \, \tilde{d}_{\mathcal{S}} (\vec{U}_{n}^{0}) \geq \delta_{b}
\end{array}
\label{Eqn:U0Far}
\end{equation}
Assume that $\tau_{\infty}^{0} = -\infty$.  Assume also that $h_{\infty}^{0} =0$. If the claim were wrong, then there would exist $ t_{*} \in I(U_{\infty}^{0})$ such that $\tilde{d}_{\mathcal{S}} \left( \vec{U}_{\infty}^{0}(t_{*}) \right) \leq \delta_{b} $. But then, taking into account (\ref{Eqn:Diffunu0}) and
Lemma \ref{lem:Convhinf} we see that

\begin{equation}
\begin{array}{ll}
\tilde{d}_{\mathcal{S}} \left( \vec{u}_{n} ( h_{n}^{0} (t_{*} - \tau_{n}^{0})  \right) & \lesssim \delta_{b},
\end{array}
\nonumber
\end{equation}
% NOTICE \\
% Letting $\vec{W}_{\sigma,c} :=  \left( \nabla ( \mathcal{T}^{c}  S_{-1}^{\sigma} (W) ), 0 \right) $
% \begin{equation}
% \begin{array}{ll}
% \inf_{(\sigma,c) \in \mathbb{R} \times \mathbb{R}^{N}} \| \vec{U}_{\infty}^{0} (t^{*}) - \vec{W}_{\sigma,c} \|_{L^{2}} & \leq \delta_{bd} \\
% \inf_{(\sigma,c) \in \mathbb{R} \times \mathbb{R}^{N}} \| \vec{U}_{n}^{0} (t^{*}) -  \vec{W}_{\sigma,c} \|_{L^{2}} & \leq 2 \delta_{bd}   \\
% \inf_{(\sigma,c) \in \mathbb{R} \times \mathbb{R}^{N}} \| T^{0}_{n} \vec{U}_{n}^{0} (t^{*}) -  \vec{W}_{\sigma,c} \|_{L^{2}} &  \leq 2 \delta_{bd} \\
% \inf_{(\sigma,c) \in \mathbb{R} \times \mathbb{R}^{N} }
% \left\| \langle D \rangle u_{n}^{0} \left( h_{n}^{0} (t_{*} - \tau_{n}^{0}) \right) -  (W_{\sigma,c})_{1} \right\|_{L^{2}}
% + \left\| \partial_{t} u_{n}^{0} \left( h_{n}^{0} (t_{*} - \tau_{n}^{0}) \right) - 0  \right\|_{L^{2}} & \leq 2 \delta_{bd} \\
% \tilde{d}_{wa,\mathcal{S}}  \left( u_{n} ( h_{n}^{0} (t_{*} - \tau_{n}^{0})  \right) & \leq 2 \delta_{bd}
% \end{array}
% \nonumber
% \end{equation}
% ***** \\
which is a contradiction, since $ t_{*} > \tau_{n}^{0} $. The other case $h_{\infty}^{0} = 1$ is easier and therefore left to the reader. \\
The case $ \tau_{\infty}^{0} \in \mathbb{R}$ is treated similarly to the case $ \tau_{\infty}^{0} = -\infty $ and therefore the proof is omitted. \\
\\
Hance, using also Result \ref{Res:PositivK} and Proposition \ref{Prop:SignProp} (or Proposition \ref{Prop:VarEst}) we see that $K(U_{\infty}(t)) \geq
\min \left( \| \nabla U_{\infty}^{0}(t) \|^{2}_{L^{2}},  k(\delta_{b}) \right) $. \\
\\
Define $\vec{U}_{c}(t) := \vec{U}_{n}^{0}(t+c)$ if $h_{\infty}^{0}=1$ ( resp. $\vec{U}_{c}(t):= \vec{U}_{\infty}^{0}(t + c) $ if $h_{\infty}^{0}=0$
with $ c > \tau_{\infty}^{0}$. Notice from the estimates above, the results in \cite{kriegnakschlagnonrad} for $d \in \{ 3,4 \}$, or the arguments in \cite{kenmerwave} (see Section $5$) for $d \in \{ 3,4,5 \}$, that we cannot have $h_{\infty}^{0} = 0$. If $E(\vec{U}_{c}) < E_{c}$ then, by definition of $E_{c}$, we have $ \| U_{c} \|_{ L_{t}^{\frac{2(d+1)}{d-2}} L_{x}^{\frac{2(d+1)}{d-2}} (0,\infty)} < \infty$ or, equivalently,
$\| U_{n}^{0} \|_{ L_{t}^{\frac{2(d+1)}{d+2}} L_{x}^{\frac{2(d+1)}{d-2}} (c,\infty)} < \infty$. From the construction of the nonlinear profile, we also know that
$\| U_{n}^{0} \|_{ L_{t}^{\frac{2(d+1)}{d-2}} L_{x}^{\frac{2(d+1)}{d-2}} ((\tau_{\infty}^{0},c))} < \infty$. Taking also into account
(\ref{Eqn:EstSumL2Norm}), this leads to a contradiction by applying again Result \ref{res:PertubationResult}, with $\tilde{w}:= u_{n}^{0}$ and $w:=u_{n}$. Hence $E(\vec{U}_{c}) = E_{c}$.

\begin{rem}
Notice that by letting $\alpha \rightarrow 0$ in (\ref{Eqn:DecouplNRJ}) we see that $U_{\infty}^{j} = 0$ if  $h_{\infty}^{j} =0 $ (resp. $U_{n}^{j} =0$
if $h_{n}^{j}=1$ ) for $j \geq 1$ and $\vec{w}_{n}^{k} \rightarrow 0$ in $L^{2}$. Hence

\begin{equation}
\begin{array}{l}
\vec{u}_{n}(0) = e^{- i t_{n}^{0} \langle D \rangle} T_{n}^{0 }\vec{\phi}^{0} + O_{L^{2}} (\alpha)
\end{array}
\nonumber
\end{equation}
\label{rem:Remun0}
\end{rem}

Let $U_{c}:= \Re \left( \langle D \rangle^{-1} \vec{U}_{c} \right)$. We prove the following claim:

\begin{claim}
There exists $c(t) \in \mathbb{R}^{d}$, $t \in [0, T_{+}(U_{c}))$ such that

\begin{equation}
\begin{array}{l}
K  = \left\{
v(t) := \vec{U}_{c} (x-c(t), t), t \in [0, T_{+}(U_{c})) \right\}
\end{array}
\end{equation}
is precompact in $L^{2}$. In fact we have $ T_{+}(U_{c}) = \infty $.
\end{claim}

\begin{proof}

If $K$ were not precompact, then there would exist $\eta > 0$ and a sequence $t_{n} \in [0, T_{+}(U_{c})) $ such that

\begin{equation}
\begin{array}{ll}
\left\| \vec{U}_{c} (t_{n}, \cdot - c_{0}) - \vec{U}_{c}(t_{n^{'}})  \right\|_{L^{2}} \geq \eta
\end{array}
\label{Eqn:CdtionNonComp}
\end{equation}
for $n \neq n'$ and for all $c_{0} \in \mathbb{R}^{d}$. In view of the continuity of the flow, we must have $t_{n} \rightarrow T_{+}(U_{c})$
in order not to violate (\ref{Eqn:CdtionNonComp}). Applying (\ref{Eqn:Decompvn}) to $\vec{v}_{n}(0):= \vec{U}_{c}(t_{n})$  and
following the same steps down to Remark \ref{rem:Remun0}, we see that there exists
$ (\tilde{t}_{n}^{0}, \vec{\tilde{\phi}}_{0}) \in \mathbb{R} \times L^{2}$ such that

\begin{equation}
\begin{array}{ll}
\vec{U}_{c}(t_{n}) & = e^{- i \tilde{t}_{n}^{0} \langle D \rangle} \tilde{T}_{n}^{0} \vec{\tilde{\phi}}_{0} + o_{L^{2}}(1),
\end{array}
\label{Eqn:Uctn}
\end{equation}
with $ \tilde{T}_{n}^{0} \vec{\tilde{\phi}}_{0} := \vec{\tilde{\phi}}_{0}(\cdot - \tilde{x}^{0}_{n}) $ for some $\tilde{x}^{0}_{n} \in \mathbb{R}^{d}$.
Assume that (up to a subsequence) that $\tilde{t}_{n}^{0} \rightarrow -\infty$. Then

\begin{equation}
\begin{array}{ll}
\left\|  \Re  \left(  \langle D \rangle^{-1} e^{i t \langle D \rangle} \vec{U}_{c}(t_{n}) \right)  \right\|_{L_{t}^{\frac{2(d+1)}{d-2}} L_{x}^{\frac{2(d+1)}{d-2}}
((0, \infty)) }  \\
\lesssim \left\| \Re \left( \langle D \rangle^{-1}  e^{ i t \langle D \rangle} \vec{\tilde{\phi}}_{0}  \right)  \right\|_{ L_{t}^{\frac{2(d+1)}{d-2}} L_{x}^{\frac{2(d+1)}{d-2}} ((- \tilde{t}^{0}_{n}, \infty))} + o_{H^{1}}(1) \\
\ll 1 \cdot
\end{array}
\nonumber
\end{equation}
But then, by the local well-posedness theory, this would imply that $  \| U_{c}(t+t_{n}) \|_{ L_{t}^{\frac{2(d+1)}{d-2}} L_{x}^{\frac{2(d+1)}{d-2}} ((0, \infty))} = \| U_{c} \|_{ L_{t}^{\frac{2(d+1)}{d-2}} L_{x}^{\frac{2(d+1)}{d-2}} ((t_{n}, \infty))} < \infty$, which is a contradiction. \\
Assume that (up to a subsequence) that $\tilde{t}^{0}_{n} \rightarrow \infty$. Then

\begin{equation}
\begin{array}{l}
\left\| \Re \left( \langle D \rangle^{-1} e^{it \langle D \rangle} \vec{U}_{c}(t_{n}) \right) \right\|_{ L_{t}^{\frac{2(d+1)}{d-2}} L_{x}^{\frac{2(d+1)}{d-2}}(-\infty,0)} \\ \lesssim \left\| \Re \left( \langle D \rangle^{-1} e^{i t \langle D \rangle}  \vec{\tilde{\phi}}_{0} \right) \right\|_{ L_{t}^{\frac{2(d+1)}{d-2}} L_{x}^{\frac{2(d+1)}{d-2}} ((-\infty, -\tilde{t}^{0}_{n}) } + o_{H^{1}}(1) \\
\ll 1
\end{array}
\nonumber\\\\
\end{equation}
But then, by the local well-posedness theory, this would imply that $  \| U_{c}(t+t_{n}) \|_{ L_{t}^{\frac{2(d+1)}{d-2}} L_{x}^{\frac{2(d+1)}{d-2}} ((-\infty, 0))} =
\| U_{c} \|_{ L_{t}^{\frac{2(d+1)}{d-2}} L_{x}^{\frac{2(d+1)}{d-2}} (  (-\infty, t_{n}) )} \ll 1$, which is a contradiction since $t_{n} \rightarrow T_{+}(U_{c})$. \\
Assume that (up to a subsequence) $\tilde{t}^{0}_{n} \rightarrow \tilde{t}$ for some $\tilde{t} \in \mathbb{R}$. Then, in view of (\ref{Eqn:CdtionNonComp}) and
(\ref{Eqn:Uctn}), we get after an appropriate change of variable

\begin{equation}
\begin{array}{l}
\left\| e^{-i \tilde{t}^{0}_{l} \langle D \rangle}  \vec{\tilde{\phi}}_{0} ( x + \tilde{x}^{0}_{l'} - \tilde{x}^{0}_{l} - c_{0} ) -
e^{- i \tilde{t}^{0}_{l'} \langle D \rangle} \vec{\tilde{\phi}}_{0} (x) \right\|_{L^{2}} \geq \frac{\eta}{2} \cdot
\end{array}
\nonumber
\end{equation}
But then, choosing $c_{0} := \tilde{x}^{0}_{l'} -  \tilde{x}^{0}_{l}$ we see that the estimate above cannot hold for $l \gg 1$ and $l' \geq l$ . \\
Next we claim that $T_{+}(U_{c}) = \infty$. If not let $(t_{l})_{l \geq 1}$ be a sequence such that $t_{l} \rightarrow T_{+}(U_{c})$.  By precompactness there
exists $\vec{V} \in L^{2}$ such that $\vec{U}_{c} \left(t_{l}, \cdot - c(t_{l}) \right) \rightarrow \vec{V}$ in $L^{2}$. Now let $\vec{H}$ (resp. $\vec{H}_{l}$) be the solution of $ (i \partial_{t}  + \langle D \rangle)\vec{H} = F(H)$ (resp. $ (i \partial_{t}  + \langle D \rangle)\vec{H}_{l} = F(H_{l})$ ) with $H:= \Re \left( \langle D \rangle^{-1} \vec{H} \right)$ ( resp. $H_{l}:= \Re \left( \langle D \rangle^{-1} \vec{H}_{l} \right)$) that satisfies
$\vec{H}(T_{+}(U_{c})) = \vec{V} $ (resp. $ \vec{H}_{l}(T_{+}(U_{c})) = \vec{U}_{c}(t_{l}, \cdot -c(t_{l}))  $). Then, by (\ref{Eqn:StrichKg0}) and the dominated
convergence theorem

\begin{equation}
\begin{array}{l}
\left\| \Re \left( \langle D \rangle^{-1} e^{i (t- T_{+}(U_{c})) \langle D \rangle } \vec{H}_{l}(T_{+}(U_{c})) \right)  \right\|_{ L_{t}^{\frac{2(d+1)}{d-2}}
L_{x}^{\frac{2(d+1)}{d-2}} ( T_{+}(U_{c}), T_{+}(U_{c}) + \delta)) } \\
 \lesssim \left\| \vec{H}_{l} \left(T_{+}(U_{c}) \right) - \vec{H} \left( T_{+}(U_{c}) \right)  \right\|_{L^{2}} \\
+ \left\| \Re \left( \langle D \rangle^{-1} e^{i (t- T_{+}(U_{c})) \langle D \rangle } \vec{H}(T_{+}(U_{c}) \right)
\right\|_{ L_{t}^{\frac{2(d+1)}{d-2}} L_{x}^{\frac{2(d+1)}{d-2}} ( T_{+}(U_{c}), T_{+}(U_{c}) + \delta)) }  \\
 \ll 1,
\end{array}
\nonumber
\end{equation}
for some $\delta \ll 1$. Consequently

\begin{equation}
\begin{array}{l}
\left\| \cos{ \left( (t - T_{+}(U_{c})) \langle D \rangle \right)} H_{l} (T_{+}(U_{c})) +
\frac{\sin{ \left( (t- T_{+}(U_{c})) \langle D \rangle \right)}}{ \langle D \rangle} \partial_{t} H_{l}(T_{+}(U_{c})) \right\|_{ L_{t}^{\frac{2(d+1)}{d-2}}
L_{x}^{\frac{2(d+1)}{d-2}} ( T_{+}(U_{c}), T_{+}(U_{c}) + \delta))}  \\
\ll 1.
\end{array}
\nonumber
\end{equation}
Hence by the local well-posedness theory, we see that we can extend $\vec{H}_{l}$ to
$T_{+}(U_{c}) + \delta$. Since $\vec{H}_{l}(t) = \vec{U}_{c} \left( \cdot - c(t_{l}), t + t_{l} - T_{+} (U_{c})  \right) $, this is a contradiction.

\end{proof}
Next we prove the following claim:

\begin{claim}
$\vec{U}_{c}$ does not exist.
\end{claim}

\begin{proof}
We may assume without loss of generality that $c(0)=0$.
A straightforward modification of an argument using the Lorentz transform in \cite{ibramasmnak} (the argument is also used in Section $2.4$ in \cite{nakschlagbook}
for the cubic focusing nonlinear Klein-Gordon equation) shows that $|\vec{P}(\vec{U}_{c})| \lesssim \bar{\epsilon}$, with $\vec{P} (\vec{U}_{c}) := \int_{\mathbb{R}^{d}} \ \nabla U_{c} \partial_{t} U_{c} \; dx  $ being the conserved-in-time momentum.
By the previous claim, we see that there exists $ R_{0}:= R_{0}(\bar{\epsilon}) \gg 1 $ such that

\begin{equation}
\begin{array}{l}
E_{c(t),R_{0}}(\vec{U}_{c}(t)) \leq \bar{\epsilon} E(\vec{U}_{c}), \, t \in [0, \infty)
\end{array}
\label{Eqn:NrjDecayInf}
\end{equation}
From the identity  \footnote{ Here $ X_{R}(\vec{U}_{c}(t)) $ denoting the localized center of energy, i.e
$ X_{R} (\vec{U}_{c}(t)) := \int \chi_{R} (x) x e \left( \vec{U}_{c}(t,x) \right) \,dx $, with $e \left( \vec{U}_{c} (t,x) \right)$ such that
$E \left( \vec{U}_{c}(t) \right) = \int e \left( \vec{U}_{c}(t,x) \right) \, dx$ }

\begin{equation}
\begin{array}{ll}
 \partial_{t} X_{R}(\vec{U}_{c})  & = - d \times \vec{P}(\vec{U}_{c}) + \int_{|x| \geq R} e(\vec{U}_{c}(t)) \, dx,
\end{array}
\nonumber
\end{equation}
and a slight modification of  an argument in \cite{ibramasmnak}, we see that if $ R  \gg R_{0}:= R_{0}(\bar{\epsilon}) $,  then

\begin{equation}
\begin{array}{l}
|c(t) - c(0)| \leq R,
\end{array}
\nonumber
\end{equation}
for $ 0 < t < t_{0} \approx \frac{R}{\bar{\epsilon}}$. Hence

\begin{equation}
\begin{array}{l}
E_{2R,0} (\vec{U}_{c}(t)) \leq \bar{\epsilon} E(\vec{U}_{c}), \, t \in [0, t_{0}]
\end{array}
\nonumber
\end{equation}
By using the previous claim, we can prove (see \cite{ibramasmnak}) that

\begin{equation}
\begin{array}{ll}
\int_{0}^{t_{0}} \| \partial_{t} U_{c}(t) \|^{2}_{L^{2}} + \| U_{c}(t) \|^{2}_{L^{2}} \, dt & \lesssim E(\vec{U}_{c}) + \int_{0}^{t_{0}} \| \nabla U_{c}(t) \|^{2}_{L^{2}} \, dt
\end{array}
\nonumber
\end{equation}
and, consequently

\begin{equation}
\begin{array}{l}
t_{0} E(\vec{U}_{c})  \lesssim E(\vec{U}_{c}) + \int_{0}^{t_{0}} \| \nabla U_{c}(t) \|^{2}_{L^{2}} \, dt
\end{array}
\nonumber
\end{equation}
Next we apply (\ref{Eqn:VirialIdentNLKG}) to $U_{c}$ and $ w(x):= \chi \left( \frac{|x|}{2R} \right)$ ;  we integrate (\ref{Eqn:VirialIdentNLKG}) \footnote
{with $t-t_{2} + m $ substituted with $2R$} on $[0,t_{0})$ and we get (using Proposition \ref{Prop:VarEst})

\begin{equation}
\begin{array}{ll}
R \gtrsim \left| \left[ \langle w \partial_{t} U_{c} \cdot x \cdot \nabla U_{c} + \frac{d}{2} U_{c} \rangle \right]^{t_{0}}_{0} \right|  &
\gtrsim t_{0} - O(1) - \bar{\epsilon}  t_{0} \gtrsim \frac{R}{\bar{\epsilon}}
\end{array}
\nonumber
\end{equation}
This is a contradiction.

\end{proof}

\section{Proof of Lemma \ref{lem:Convhinf}}
\label{Section:ProofLemmaConv}

First observe that the finiteness of $\| U_{\infty}^{j} \|_{L_{t}^{\frac{d+2}{d-2}} L_{x}^{\frac{2(d+2)}{d-2}} (\mathbb{R})}$ implies that
$\| \vec{U}_{\infty}^{j} \|_{L_{t}^{\infty} L^{2} (\mathbb{R})} < \infty$: indeed, partitioning $\mathbb{R}^{+}$ into
a finite number of subintervals $(J_{k}:= [a_{k},b_{k}] )_{1 \leq k \leq l}$ such that $ \| U_{\infty}^{j} \|_{L_{t}^{\frac{d+2}{d-2}} L_{x}^{\frac{2(d+2)}{d-2}} (J_{k})} = \eta$, $1 \leq k < l $, and $ \| U_{\infty}^{j} \|_{L_{t}^{\frac{d+2}{d-2}} L_{x}^{\frac{2(d+2)}{d-2}} (J_{l})} \leq \eta$, we see from
(\ref{Eqn:StrichWave0}) and H\"older inequality that  for all $ K_{k}:= [a_{k},.) \subset J_{k}$ we

\begin{equation}
\begin{array}{l}
\| \vec{U}_{\infty}^{j} \|_{L_{t}^{\infty} L_{x}^{2} (K_{k})} + \| U_{\infty}^{j} \|_{L_{t}^{\frac{d+2}{d-2}} L_{x}^{\frac{2(d+2)}{d-2}} (K_{k})} \\
\lesssim  \| \vec{U}_{\infty}^{j}(a_{k}) \|_{L^{2}} +  \| U_{\infty}^{j} \|^{2^{*}-1}_{L_{t}^{\frac{d+2}{d-2}} L_{x}^{\frac{2(d+2)}{d-2}} (K_{k})}
\end{array}
\nonumber
\end{equation}
Hence a continuity argument and an iteration over $k$ show that $\| \vec{U}_{\infty}^{j} \|_{L_{t}^{\infty} L_{x}^{2} (\mathbb{R}^{+})} < \infty$.
Proceeding similarly we have  $\| \vec{U}_{\infty}^{j} \|_{L_{t}^{\infty} L_{x}^{2} (\mathbb{R}^{-})} < \infty$. \\
\\
If $J$ is an interval let $ X(J) :=  \| \vec{U}_{n}^{j} - \vec{U}_{\infty}^{j} \|_{L_{t}^{\infty} L_{x}^{2} (J)}
+ \| U^{j}_{n} - U^{j}_{\infty} \|_{L_{t}^{\frac{d+2}{d-2}} L_{x}^{\frac{2(d+1)}{d-2}}(J)} $. Let $ 0 <  \epsilon \ll 1 $ and $S \gg 1$ be a large constant. \\
We first estimate $X([-S,S])$. We have

\begin{equation}
\begin{array}{l}
\partial_{tt}  \left( U_{n}^{j} -  U_{\infty}^{j}  \right)  -  \triangle \left( U_{n}^{j} -  U_{\infty}^{j} \right)  = F (U_{n}^{j}) - F (U_{\infty}^{j}) - (h_{n}^{j})^{2}  U_{n}^{j} \cdot
\end{array}
\nonumber
\end{equation}
Let $ 0 < \eta \ll 1$ be a constant small enough such that all the estimate below are true. We partition $[-S,S]$ into
a finite number of subintervals $(J_{k}:= [a_{k},b_{k}] )_{1 \leq k \leq l}$ such that $ \| U_{\infty}^{j} \|_{L_{t}^{\frac{d+2}{d-2}} L_{x}^{\frac{2(d+2)}{d-2}} (J_{k})} = \eta$, $1 \leq k < l $, and $ \| U_{\infty}^{j} \|_{L_{t}^{\frac{d+2}{d-2}} L_{x}^{\frac{2(d+2)}{d-2}} (J_{l})} \leq \eta$. We have

\begin{equation}
\begin{array}{l}
\| F (U_{n}^{j}) - F (U_{\infty}^{j}) \|_{L_{t}^{1} L_{x}^{2}(J_{k})} \\
 \lesssim   \| U_{n}^{j} - U_{\infty}^{j} \|_{L_{t}^{\frac{d+2}{d-2}} L_{x}^{\frac{2(d+2)}{d-2}}(J_{k})}
\left(
\| U_{n}^{j} \|^{2^{*}-2}_{L_{t}^{\frac{d+2}{d-2}} L_{x}^{\frac{2(d+2)}{d-2}}(J_{k})}
+ \| U_{\infty}^{j} \|^{2^{*}-2}_{L_{t}^{\frac{d+2}{d-2}} L_{x}^{\frac{2(d+2)}{d-2}}(J_{k})}
\right)
 \\
\lesssim  \| U_{n}^{j} - U_{\infty}^{j} \|^{2^{*}-1}_{L_{t}^{\frac{d+2}{d-2}} L_{x}^{\frac{2(d+2)}{d-2}}(J_{k})}
+ \| U_{n}^{j} - U_{\infty}^{j} \|_{L_{t}^{\frac{d+2}{d-2}} L_{x}^{\frac{2(d+2)}{d-2}}(J_{k})} \eta^{2^{*}-2}
\end{array}
\cdot
\nonumber
\end{equation}
The H\"older-in-time inequality yields $ \| (h_{n}^{j})^{2} U_{n}^{j} \|_{L_{t}^{1} L_{x}^{2} (J_{k})}  \lesssim
(h_{n}^{j})^{2} S  \| \vec{U}_{n}^{j} \|_{L_{t}^{\infty} L_{x}^{2} (J_{k})} $. \\
Let $ M := \| U_{\infty}^{j} \|_{L_{t}^{\frac{d+2}{d-2}} L_{x}^{\frac{2(d+2)}{d-2}}(\mathbb{R})}$. We see from  (\ref{Eqn:StrichWave0})
that

\begin{equation}
\begin{array}{l}
\| U_{n}^{j} - U_{\infty}^{j} \|_{L_{t}^{\frac{d+2}{d-2}} L_{x}^{\frac{2(d+2)}{d-2}}(J_{k})} \\
\lesssim \left\|  \vec{U}_{n}^{j} (a_{k}) - \vec{U}_{\infty}^{j} (a_{k}) \right\|_{L^{2}} +   \| (h_{n}^{j})^{2} U_{n}^{j} \|_{L_{t}^{1} L_{x}^{2} (J)} + \| F (U_{n}^{j}) - F (U_{\infty}^{j}) \|_{L_{t}^{1} L_{x}^{2}(J_{k})}  \cdot
\end{array}
\nonumber
\end{equation}
Hence, combining the estimates above, iterating over $k$, and recalling that $\vec{U}_{n}^{j} (a_{1}) = \vec{U}_{\infty}^{j} (a_{1})$, we get
$ X([-S,S]) \lesssim \epsilon  $ for $ n \gg 1$. \\
We then have to prove that $  \| U_{n}^{j} - U_{\infty}^{j} \|_{L_{t}^{\frac{d+2}{d-2}} L_{x}^{\frac{2(d+2)}{d-2}} ((S, \infty))} \lesssim \epsilon $ for $n \gg 1$.
We get from (\ref{Eqn:StrichWave0}) $ \left\| \cos{(tD)} \Re \left( D^{-1} \vec{\phi}^{j} \right) + \frac{\sin (tD)}{D} \Im (\vec{\phi}^{j})
\right\|_{L_{t}^{\frac{d+2}{d-2}} L_{x}^{\frac{2(d+2)}{d-2}} (-\infty, \infty)} \lesssim \| \vec{\phi}_{j} \| _{L^{2}} $. This estimate, combined with the dominated convergence theorem, shows that $ \left\| \cos{(tD)} \Re \left( D^{-1} \vec{\phi}^{j} \right) + \frac{\sin (tD)}{D} \Im (\vec{\phi}^{j})
\right\|_{L_{t}^{\frac{d+2}{d-2}} L_{x}^{\frac{2(d+2)}{d-2}} ((S, \infty))} \ll \epsilon $. Let $J: = [S,b]  $ be an interval. Applying  (\ref{Eqn:StrichWave0}) again, we get

\begin{equation}
\begin{array}{l}
\| U_{\infty}^{j} \|_{L_{t}^{\frac{d+2}{d-2}} L_{x}^{\frac{2(d+2)}{d-2}} (J)} \lesssim \epsilon + \| F(U_{\infty}^{j}) \|_{L_{t}^{1} L_{x}^{2} (J)}
\lesssim  \epsilon + \| U_{\infty}^{j} \|^{2^{*}-1}_{L_{t}^{\frac{d+2}{d-2}} L_{x}^{\frac{2(d+2)}{d-2}}(J)} \cdot
 \end{array}
\label{Eqn:EstUinfj}
\end{equation}
Hence a continuity argument shows that $ \| U_{\infty}^{j} \|_{L_{t}^{\frac{d+2}{d-2}} L_{x}^{\frac{2(d+2)}{d-2}} (S, \infty)} \lesssim  \epsilon $.
Let $ \phi \in L^{2} (\mathbb{R}^{d}) $. We claim that for $n \gg 1$

\begin{equation}
\begin{array}{l}
\left\| \frac{e^{ \pm i t \langle D  \rangle_{n}^{j}}}{ \langle D \rangle_{n}^{j}} \phi \right\|_{L_{t}^{\frac{d+2}{d-2}} L_{x}^{\frac{2(d+2)}{d-2}} (S,\infty)} \ll
\epsilon \cdot
\end{array}
\label{Eqn:EstDispnPhi}
\end{equation}
Indeed let $\phi_{\epsilon} $ be a Schwartz function such that  $ \| \phi_{\epsilon} - \phi \|_{L^{2}} \ll  \epsilon $.  From (\ref{Eqn:StrichWave0Sc}) we see that
$
\left\| \frac{e^{\pm i t \langle D  \rangle_{n}^{j}}}{ \langle D \rangle_{n}^{j}} (\phi_{\epsilon} - \phi )
\right\|_{L_{t}^{\frac{d+2}{d-2}} L_{x}^{\frac{2(d+2)}{d-2}} ((S, \infty))} \lesssim \| \phi_{\epsilon} - \phi \|_{L^{2}}
$. Hence it suffice to prove that (\ref{Eqn:EstDispnPhi}) holds with $ \phi $ a Schwartz function. \\
Recall that $  e^{ \pm it h_{n}^{j} \langle D \rangle} T_{n}^{j} = T_{n}^{j} e^{\pm i t \langle D \rangle_{n}^{j}} $ and that
$ \langle D \rangle^{-m} T_{n}^{j} = (h_{n}^{j})^{m} T_{n}^{j}  \left( \langle D  \rangle_{n}^{j} \right)^{-m} $.  Hence
$
\frac{ e^{i t \langle D \rangle_{n}^{j}}} {\langle D \rangle_{n}^{j}} = (h_{n}^{j})^{-1}  (T_{n}^{j})^{-1} \frac{e^{i t h_{n}^{j} \langle D \rangle}}{\langle D \rangle}
T_{n}^{j}
$
and, denoting by $\left( \frac{2(d+1)}{d-2} \right)'$ the number such that $ \frac{1}{\frac{2(d+1)}{d-2}} + \frac{1}{\left( \frac{2(d+1)}{d-2} \right)'} = 1$, we get
from the estimate $ \left\| \langle D \rangle^{-1} e^{i t \langle D \rangle} g \right\|_{ B^{0}_{\frac{2(d+1x)}{d-2},2}} \lesssim
\frac{1}{|t|^{\frac{3(d-1)}{2(d-1)}}} \| g \|_{B^{\frac{1}{2}}_{\left( \frac{2(d+1)}{d-2} \right)^{'},2}} $  and the embedding
 $ B^{\frac{1}{2}}_{\left( \frac{2(d+2)}{d-2} \right)^{'},2} \hookrightarrow  H^{\frac{1}{2}}_{\left( \frac{2(d+2)}{d-2} \right)^{'}}  $

\begin{equation}
\begin{array}{ll}
\left\| \frac{e^{i \pm t \langle D \rangle_{n}^{j}}} {\langle D \rangle_{n}^{j}} \phi \right\|_{L_{t}^{\frac{d+2}{d-2}} L_{x}^{\frac{2(d+2)}{d-2}} ((S, \infty))}
& = \left\| \frac{e^{i \pm t h_{n}^{j} \langle D \rangle}}{\langle D \rangle} T_{n}^{j} \phi  \right\|_{L_{t}^{\frac{d+2}{d-2}} L_{x}^{\frac{2(d+2)}{d-2}} ((S h_{n}^{j}, \infty))} \\
& \lesssim  \left\| \frac{e^{i \pm t  \langle D \rangle }}{\langle D \rangle} T_{n}^{j} \phi  \right\|_{L_{t}^{\frac{d+2}{d-2}} B^{0}_{\frac{2(d+2)}{d-2},2}
((S h_{n}^{j}, \infty))} \\
& \lesssim  \left\| |t|^{-\frac{2(d-1)}{d+2}} \right\|_{L_{t}^{\frac{d+2}{d-2}} (S h_{n}^{j},\infty)} \left\| \langle D \rangle^{\frac{2(d+1)}{d-2}} T_{n}^{j} \phi \right\|_{L^{\left( \frac{2(d+2)}{d-2} \right)^{'}}} \\
& \lesssim \frac{1}{S^{\frac{d}{d-2}}}  \left\| (\langle D \rangle_{n}^{j})^{\frac{1}{2}} \phi \right\|_{L^{\left( \frac{2(d+2)}{d-2} \right)^{'}}} \\
& \lesssim \frac{1}{S^{\frac{d}{d-2}}} \left\| \langle D \rangle^{\frac{1}{2}} \phi \right\|_{L^{\left( \frac{2(d+2)}{d-2} \right)^{'}}} \cdot
\end{array}
\label{Eqn:ListEstUnifn}
\end{equation}
Hence we get (\ref{Eqn:EstUinfj})  from (\ref{Eqn:StrichWave0Sc}) after replacing ``$U_{\infty}^{j}$'' with ``$U_{n}^{j}$''.  Hence $ \| U_{n}^{j} \|_{L_{t}^{\frac{d+2}{d-2}} L_{x}^{\frac{2(d+2)}{d-2}} (S, \infty)} \lesssim  \epsilon $. \\
Similarly $  \| U_{n}^{j} - U_{\infty}^{j} \|_{L_{t}^{\frac{d+2}{d-2}} L_{x}^{\frac{2(d+2)}{d-2}} ((-\infty, -S))} \lesssim \epsilon $ for $n \gg 1$.
Hence (\ref{Eqn:ConvStrichNormNrj}) holds.

\section{Proof of Theorem  \ref{Thm:Main} $(3)$. }
\label{Section:ProofThmMain2}

We apply the decomposition (\ref{Eqn:Decompu}) to the data $\vec{u}(0)$. Hence $ \vec{u}(0)= \vec{S}^{\sigma(0)} \left( \vec{W} + \vec{v}(0) \right) $ with
$\sigma(0)$ to be chosen and $\vec{v}(0) = \vec{S}^{-\sigma(0)} \vec{u}(0) - \vec{W} $. Let $ 0 < \beta \ll \epsilon_{*} $. Let $R \gg 1$ and
$\sigma(0) \gg 1 $  such that

\begin{equation}
\begin{array}{l}
d=3:
\left\{
\begin{array}{l}
\int_{\mathbb{R}^{d}} \left| \nabla \left( \left( 1- \chi \left( \frac{r}{R} \right) \right) W \right) \right|^{2} \, dx \leq \beta^{2} \\
e^{-2 \sigma(0)} \int_{\mathbb{R}^{d}} \left| \chi \left( \frac{r}{R} \right) W \right|^{2} + \rho^{2} \, dx \leq \beta^{2}
\end{array}
\right.
\\
d=4,5: e^{- 2 \sigma(0)} \int_{\mathbb{R}^{d}} W^{2} \, dx  \ll \beta^{2}
\end{array}
\nonumber
\end{equation}
We then apply the decomposition (\ref{Eqn:Decompv}) to $\vec{v}(0)$. Hence $\vec{v}(0) = \vec{\lambda}(0) \rho + \vec{\gamma}(0)$ with

\begin{equation}
\begin{array}{l}
\vec{\lambda}(0):= \beta ( \pm 1, 0 ), \, \beta (0, \pm 1); d=4,5: \vec{\gamma}(0)=0;  \\
d=3: \vec{\gamma}(0)= - \left( 1 - \chi \left( \frac{r}{R} \right) \right)
\vec{W}  + \left\langle  \left( 1- \chi \left( \frac{r}{R} \right) \right) \vec{W} | \vec{\rho} \right\rangle \vec{\rho}
\end{array}
\nonumber
\end{equation}
(Recall that $\chi$ is defined below (\ref{Eqn:VirialIdentNLKG})). Here $\vec{\rho} := \rho (1,1)$. We have $\| \vec{\gamma}(0) \|^{2}_{\dot{\mathcal{H}}} + \| u(0) \|^{2}_{L^{2}} \ll \beta^{2}$. From (\ref{Eqn:Equivdwa}) we see that (with  $\tau(0)=0$, $\partial_{t} \tau = e^{\sigma(t)}$)

\begin{equation}
\begin{array}{ll}
\tilde{d}^{2}_{\mathcal{S}}(\vec{u}(\tau)) \approx  \lambda_{1}^{2}(\tau) + \lambda_{2}^{2}(\tau) + \| \vec{\gamma}(\tau) \|^{2}_{\dot{\mathcal{H}}},
\end{array}
\nonumber
\end{equation}
as long as $ \vec{u}(\tau) \in \mathcal{V}$. Modifying slightly (\ref{Eqn:Ewapar}) and the steps below (\ref{Eqn:Ewapar}), we see that
%Integrating $  \partial_{\tau} ( E(\vec{u}) - E_{wa} (\vec{W} + \vec{v}_{d}) ) $
%we see from (\ref{Eqn:Ewapar}) and (\ref{Eqn:EstRes}) that (with $X$ defined in (\ref{Eqn:DfnX}))

\begin{equation}
\begin{array}{l}
X^{2}(\tau) \lesssim o(\beta^{2}) + \int_{0}^{\tau}
\begin{array}{l}
\left( \tilde{d}^{2}_{\mathcal{S}} (\vec{u}(\tau')) + \tilde{d}^{2^{*}-1}_{\mathcal{S}} (\vec{u}(\tau'))  \right) X(\tau^{'}) +
\left( \tilde{d}_{\mathcal{S}} (\vec{u}(\tau')) + \tilde{d}^{2^{*}-1}_{\mathcal{S}} (\vec{u}(\tau'))  \right) X^{2}(\tau')
\end{array}
d \tau^{'}, \; \text{with}
\end{array}
\nonumber
\end{equation}
$X(\tau) :=  \left\| \vec{\gamma}(\tau) \right\|_{\dot{\mathcal{H}}} +   \left\| u(\tau) \right\|_{L^{2}}$.

Hence (using also (\ref{Eqn:EstDynLambdaPlus}), (\ref{Eqn:EvolPremSigma}) and (\ref{Eqn:EstRes})) we see from a bootstrap argument that for
$|\tau| \in [0,\tau_{f}]$ with $\tau_{f}$ such that $\beta e^{k \tau_{f}} \approx \delta_{f} $,  $ \beta e^{k \tau_{f}} \leq \delta_{f} $,    %\footnote{in particular, this allows to choose $\tau_{f} \gg 1$ %if $\beta << 1$}

\begin{equation}
\begin{array}{l}
\left( \lambda_{1}(0), \lambda_{2}(0) \right) = \beta (\pm 1, 0) \Rightarrow (\lambda_{1}, \lambda_{2}) \approx \beta
\left( \pm \cosh{(k \tau)}, \pm \sinh{(k \tau)} \right), \\
\left( \lambda_{1}(0), \lambda_{2}(0) \right) = \beta (0, \pm 1) \Rightarrow (\lambda_{1}, \lambda_{2}) \approx \beta
\left( \pm \sinh{(k \tau)}, \pm \sinh{(k \tau)} \right), \; \text{and} \\
\left\|  X(\tau) \right\|_{L^{\infty} [- \tau_{f}, \tau_{f}]} - o(\beta) \lesssim  \beta^{2} e^{2 k |\tau|} \cdot
\end{array}
\nonumber
\end{equation}
Hence there exists $\tau_{+, med}$ (resp. $\tau_{-, med}$  ) such that $ \beta e^{k \tau_{+,med}} \sim  \delta_{b} $ (resp. $ \beta e^{- k \tau_{-,med}} \sim \delta_{b} $ )
and $\tilde{d}_{\mathcal{S}} \left( \vec{u}(t_{\pm,med}) \right) \geq \delta_{b}$  with $t_{\pm,med} := \tau^{-1} ( \pm \tau_{med})$. Moreover $E(\vec{u}) <  E_{wa} (\vec{W}) + c_{d}  \tilde{d}^{2}_{\mathcal{S}} \left( \vec{u} ( t_{\pm,med}) \right) $. Hence from Proposition \ref{Prop:SignProp}, Proposition \ref{Prop:FarFromGdPr}, and
(\ref{Eqn:ExpNrjArW}), we see that the following holds

\begin{equation}
\begin{array}{l}
(\lambda_{1}(0),\lambda_{2}(0)) = \beta (1,0), \, \Rightarrow \\
E(\vec{u}) < E_{wa}(\vec{W}), \vec{u}(0) \in A_{-,-}, t \rightarrow \pm T_{\pm}(u):u \, blows-up
\end{array}
\nonumber
\end{equation}

\begin{equation}
\begin{array}{l}
(\lambda_{1}(0),\lambda_{2}(0)) = \beta (-1,0), \, \Rightarrow  \\
E(\vec{u}) < E_{wa}(\vec{W}), \vec{u}(0) \in A_{+,+}, \, t \rightarrow \pm \infty:u \, scatters
\end{array}
\nonumber
\end{equation}

\begin{equation}
\begin{array}{l}
(\lambda_{1} (0), \lambda_{2}(0)) = \beta (0,1), \, \Rightarrow \\
E_{wa} (\vec{W}) < E(\vec{u}) < E_{wa}(\vec{W}) + \epsilon_{*}^{2}, \, \vec{u}(0) \in A_{+,-} \\
t \rightarrow - \infty: u \, scatters, \, t \rightarrow T_{+}(u):u \, blows-up
\end{array}
\nonumber
\end{equation}

\begin{equation}
\begin{array}{l}
(\lambda_{1} (0), \lambda_{2}(0)) = \beta (0, -1), \, \Rightarrow \\
E_{wa} (\vec{W}) < E(\vec{u}) < E_{wa}(\vec{W}) + \epsilon_{*}^{2}, \, \vec{u}(0) \in A_{-,+} \\
t \rightarrow -T_{-}(u): u \, blows-up, \, t \rightarrow \infty: u \, scatters
\end{array}
\nonumber
\end{equation}

%%%%%%%%%%%%%%%%%%%%%%%%%%%%%%%%%%%%%%%%%%%%%%%%%%%%%%%%%%%%%%%%%%%%%%%%%%%%%%%%%%%%%%%%%%%%%%%%%%%%
%%%%%%%%%%%%%%%%%%%%%%%%%Previous Proof of Lemma%%%%%%%%%%%%%%%%%%%%%%%%%%%%%%%%%%%%%%%%%%%%%%%%%%%%%%%%%%%%


\begin{thebibliography}{99}

\bibitem{aubin} Aubin, Thierry, \emph{Equations diff\'erentielles non lin\'eaires et probl\`eme de Yamabe concernant la courbure scalaire}, J. Math. Pures. Appl (9),
 55, 1976, 3, 269-2963

\bibitem{bahchem} H. Bahouri, Jean-Yves Chemin, R. Danchin, \emph{Fourier Analysis and Nonlinear Partial Differential Equations}, Springer, 343, 523 p., 2011, Grundlehren der mathematischen Wissenschaften.

\bibitem{bahger} H. Bahouri and P. Gerard, \emph{High frequency approximation of solutions to critical nonlinear wave equations}, AJM, Vol 121, Number 1, Feb 1999, pp 131-175

\bibitem{duykenmerle} T. Duyckaerts, C. Kenig, and F. Merle, \textit{Universality of the blow-up profile for small type II blow-up solutions
of energy-critical wave equation:  the  non-radial case}, Journal Of The European Mathematical Society, Volume 14, Issue 5, 2012, pp. 1389–1454

\bibitem{bourgjams} J. Bourgain, \emph{Global well-posedness of defocusing $3D$ critical NLS in the radial case}, JAMS 12 (1999), 145-171

\bibitem{bulutlipav} A. Bulut, M. Czubak, D. Li, N. Pavlović, X. Zhang, \emph{Stability and unconditional uniqueness of solutions for
energy critical wave equations in high dimensions}, Comm. Partial Differential Equations 38 (2013), no. 4, 575–607.

\bibitem{cazbook} Cazenave, Thierry, \emph{Semilinear Schr\"odinger equations}, Courant Lecture Notes in Mathematics, 10, New York
University Courant Institute of Mathematical Sciences, New York, 2003

\bibitem{ginebvelo} J. Ginibre and G. Velo, \emph{Time decay of finite energy solutions of the nonlinear Klein-Gordon equation
and Schr\"odinger equations}, Annales de l'I.H.P, section A, volume 43, number 4 (1985), p. 399-442

\bibitem{ginebvelostr} J. Ginibre and G. Velo, \emph{Generalized Strichartz inequalities for the wave equations}, J. Func. Anal., 133(1):50-68, 1995


\bibitem{gril} M. Grillakis, \emph{Analysis of the linearization around a critical point of an infinite-dimensional Hamiltonian system}, Comm. Pure.
Applied. Math., 43(3):299-333, 1990sob

\bibitem{ibramasmnak}  S. Ibrahim, N. Masmoudi, and K. Nakanishi, \emph{Scattering threshold for the focusing nonlinear Klein-Gordon equation}, Analysis
and PDE 4 (2011) no. 3, 405-460

\bibitem{ibramasmnakerr} S. Ibrahim, N. Masmoudi, and K. Nakanishi, \emph{Correction to the article ``Scattering threshold for the focusing nonlinear Klein-Gordon equation''}, Anal. PDE 9 (2016), no. 2, 503–514.

%\bibitem{ibramasmnakerr} S. Ibrahim, N. Masmoudi, and K. Nakanishi, \emph{Errata: Scattering threshold for the focusing nonlinear Klein-Gordon equation},

\bibitem{ibramasmnakshold} S. Ibrahim, N. Masmoudi, and  K. Nakanishi, \emph{Threshold solutions in the case of mass-shift for the critical Klein-Gordon
equation}, Trans. Amer. Math. Soc. 366 (2014), no. 11, 5653–5669.

\bibitem{keeltao} Keel, M., Tao. T., \emph{Endpoints Strichartz estimates}, Amer. J. Math., 120 (1998), 955-980

\bibitem{keraani} S. Keraani, \emph{On the defect of compactness for the Strichartz estimates of the Schr\"odinger equations},
J. Differential Equations, 175, 2001, 2, 353-392

\bibitem{kenmer} C. E. Kenig and F. Merle, \emph{Global well-posedness, scattering, and blow-up for the energy-critical, focusing, non-linear
Schr\"odinger equation in the radial case}, Invent. Math. 166 (2006), no. 3, pp. 645-675

\bibitem{kenmerwave} C. E. Kenig and F. Merle, \emph{Global well-posedness, scattering and blow-up for the energy-critical focusing non-linear
wave equation}, Acta Math 201 (2008), no. 2, 147-212

\bibitem{kriegnakschlagnonrad}  J. Krieger, K. Nakanishi, and W. Schlag, \emph{Global dynamics of the nonradial energy-critical wave
equation above the ground state energy}, Discrete and Continuous Dynamical Systems Volume 33, Issue 6, 2013, 2423-2450

\bibitem{kriegnakschlagrad} J. Krieger, K. Nakanishi, and W. Schlag, \emph{Global dynamics away from the ground state for the energy-critical
  nonlinear wave equation}, American Journal of Mathematics,  Volume 135, Number 4, August 2013, 935-965

\bibitem{linsog} H. Lindblad, C. D Sogge, \emph{On existence and scattering with minimal regularity for semilinear wave equations}, J.
Func.Anal 219 (1995), 227-252

\bibitem{lions}  Lions, P-L, \emph{The concentration-compactness in the calculus of variations. The limit case}, I. Rev. Mat.
Iberoamericana 1 (1985), no. 1, 145-201

\bibitem{merlevega} Merle, F., Vega, L., \emph{Compactness at blow-up time for $L^{2}$ solutions of the critical nonlinear
Schr\"odinger equation in $2D$}, Internat. Math. Res. Notices. 1998, no 8., 399-425

\bibitem{nakashimrn} Nakanishi, K., \emph{Scattering theory for the nonlinear Klein-Gordon equation with critical Sobolev power }, Internat. Math. Res. Notices 1999, no. 1, 31-60

\bibitem{nakschlagkg} Nakanishi, K., Schlag W., \emph{Global dynamics above the ground state energy for the focusing nonlinear Klein-Gordon
equation}, Journal. Diff. Eq. 250 (2011), 2299-22

\bibitem{nakschlagkg2} Nakanishi, K., Schlag W., \emph{Global dynamics above the ground state energy for the focusing nonlinear Klein-Gordon
equation without a radial assumption}, Arch. Rational Mech. Analysis, no 3, 203 (2012), 809-851

\bibitem{nakaschlagschrod} Nakanishi, K., Schlag W., \emph{Global dynamics above the ground state energy for the cubic NLS equation
in 3D}, Calc. Var. and PDE, no 1-2, 44 (2012), 1-45

\bibitem{nakashtriroyschrod} Nakanishi, K., Roy T., \emph{Global dynamics above the ground state for the energy-critical
Schr\"odinger equation with radial data}, Communications On Pure And Applied Analysis, Volume 15, Issue 6, November 2016, 2023-2058.

%\bibitem{nakschlagwave} Nakanishi, K., Schlag W., \emph{Global dynamics of the nonradial energy-critical wave equation above the ground state
%energy}, Discrete and Continuous Dynamical Systems Volume 33, Issue 6, 2013, 2423-2450

\bibitem{nakschlagbook} Nakanishi, K., Schlag, W., \emph{Invariant manifolds and dispersive Hamiltonian evolution equations}, Zurich
Lectures in Advanced Mathematics, EMS, 2011

\bibitem{paynesatt} L. E. Payne and D. H. Sattinger, \emph{Saddle points and instability of nonlinear hyperbolic equations}, Israel J. Math 22 (1975)
, 272-303

\bibitem{pecher1} H. Pecher, \emph{Nonlinear small data scattering for the wave and Klein-Gordon equation}, Math. Z., 185 (1984), pp 261-270

\bibitem{pecher2} H. Pecher, \emph{Low energy scattering for Klein-Gordon equations}, J. Funct. Analysis 63 (1985), pp. 101-122

\bibitem{schlag} W. Schlag, \emph{Spectral theory and nonlinear differential equations: a survey}, Discrete. Cont. Dyn. Syst., 15(3):703-723, 2006

%\bibitem{stafftao} J. Colliander, M. Keel, G. Staffilani, H. Takaoka and T. Tao, \emph{Global existence and scattering in the energy
%space for the critical Scr\"odinger equation on $\mathbb{R}^{3}$}, Annals of Math, 167 (2007), 767-865

%\bibitem{strauss} W.A. Strauss, \emph{Nonlinear scattering theory at low energy}, J. Func. Anal., 41 (1981), pp. 110-133

\bibitem{taovisan} T. Tao and M. Visan, \emph{Stability of energy-critical nonlinear Schr\"odinger equations in high dimensions}, Electron. J. Differential
Equations 118 (2005), 28 pp.

\bibitem{tal} Talenti, Giorgio, \emph{Best Constant In Sobolev Inequality}, Ann. Mat. Pura. Appl., (4) 110 1976, 353-372

\end{thebibliography}
\end{document}